\documentclass{article}
\usepackage{amssymb}
\usepackage{amsmath}

\newcommand{\ie}{i.e.}
\newcommand{\lc}{l.c.}
\newcommand{\nsf}{n.s.f.}
\newcommand{\GNS}{G.N.S.}
\newcommand{\ot}{\otimes}
\newcommand{\de}{\Delta}
\newcommand{\sde}{\delta}
\newcommand{\Mh}{\hat M}
\newcommand{\deh}{\hat\Delta}

\newcommand{\id}{\iota}
\newcommand{\om}{\omega}

\newcommand{\cs}{C^{*}}

\newcommand{\R}{\mathbb{R}}
\newcommand{\Rset}{\mathbb{R}}
\newcommand{\C}{\mathbb{C}}
\newcommand{\Cset}{\mathbb{C}}

\newcommand{\N}{\mathbb{N}}
\newcommand{\Nset}{\mathbb{N}}
\newcommand{\Z}{\mathbb{Z}}

\newcommand{\bh}{\mathcal{B}(H)}

\newcommand{\nfs}{n.s.f.}
\newcommand{\gns}{G.N.S.}
\newcommand{\fie}{\varphi}

\def\11{\mbox1\hspace{-.25em}\text{I}}

\newtheorem{theorem}{Theorem}

\newtheorem{lemma}{Lemma}
\newtheorem{proposition}{Proposition}

\newenvironment{proof}{\textbf{Proof.}}{\hfill $\blacksquare$}

\newenvironment{remark}{\textbf{Remark.}}{}

\begin{document}
\title{Twisting and Rieffel's deformation of locally compact quantum groups.
Deformation of the Haar measure.}
\author{Pierre Fima\footnote{Laboratoire de Math\'ematiques,
Universit\'e de Franche-Comt\'e, 16 route de Gray, 25030 Besancon Cedex, France.
E-mail: fima@math.unicaen.fr} \,\,\,and Leonid Vainerman\footnote{Laboratoire de
Math\'ematiques Nicolas Oresme, Universit\'e de Caen, B.P. 5186, 14032 Caen
Cedex, France. E-mail: vainerman@math.unicaen.fr}}
\date{}
\maketitle

\begin{abstract}
We develop the twisting construction for locally compact quantum groups.
A new feature, in contrast to the previous work of M. Enock and the
second author, is a non-trivial deformation of the Haar measure. Then
we construct Rieffel's deformation of locally compact quantum groups and
show that it is dual to the twisting. This allows to give new interesting
concrete examples of locally compact quantum groups, in particular, deformations
of the classical $az+b$ group and of the Woronowicz' quantum $az+b$ group.
\end{abstract}

\section{Introduction}

The problem of extension of harmonic analysis on abelian locally compact
(l.c.) groups, to non abelian ones, leads to the introduction of more general
objects. Indeed, the set $\hat G$ of characters of an abelian l.c. group
$G$ is again an abelian l.c. group - the dual group of $G$. The
Fourier transform maps functions on $G$ to functions on $\hat G$,
and the Pontrjagin duality theorem claims that $\hat{\hat G}$ is
isomorphic to $G$. If $G$ is not abelian, the set of its
characters is too small, and one should use instead the set $\hat
G$ of (classes of) its unitary irreducible representations and
their matrix coefficients. For compact groups, this leads to the
Peter-Weyl theory and to the Tannaka-Krein duality, where $\hat G$
is not a group, but allows to reconstruct $G$. Such a non-symmetric
duality was established for unimodular groups by W.F. Stinespring,
and for general l.c. groups by P. Eymard and T. Tatsuuma.

In order to restore the symmetry of the duality, G.I. Kac
introduced in 1961 a category of {\it ring groups} which contained
unimodular groups and their duals. The duality constructed by Kac
extended those of Pontrjagin, Tannaka-Krein and Stinespring. This
theory was completed in early 70-s by G.I. Kac and the second author,
and independently by M. Enock and J.-M. Schwartz, in order to cover
all l.c. groups. The objects of this category are called {\it Kac
algebras} \cite{ES}. L.c. groups and their duals can be viewed
respectively as commutative and co-commutative Kac algebras, the
corresponding duality covered all known versions of duality for l.c.
groups.

Quantum groups discovered by V.G. Drinfeld and others gave new
important examples of Hopf algebras obtained by deformation of
universal enveloping algebras and of function algebras on Lie
groups. Their operator algebraic versions did not verify some of
Kac algebra axioms and motivated strong efforts to construct a
more general theory. Important steps in this direction were made
by S.L. Woronowicz with his theory of compact quantum groups and a
series of important concrete examples, S. Baaj and G. Skandalis
with their fundamental concept of a multiplicative unitary and A.
Van Daele who introduced an important notion of a multiplier Hopf
algebra. Finally, the theory of l.c. quantum groups was proposed
by J. Kustermans and S. Vaes \cite{VaesLC}, \cite{KV2}.

A number of "isolated" examples of non-trivial (i.e., non
commutative and non cocommutative) l.c. quantum groups was
constructed by S.L. Woronowicz and other people. They were
formulated in terms of generators of certain Hopf $*$-algebras and
commutation relations between them. It was much harder to
represent them by operators acting on a Hilbert space, to
associate with them an operator algebra and to construct all
ingredients of a l.c. quantum group. There was no general approach
to these highly nontrivial problems, and one must design specific
methods in each specific case (see, for example, \cite{Wor},
\cite{Van}).

In \cite{Enock}, \cite{Vain} M. Enock and the second author
proposed a systematic approach to the construction of non-trivial
Kac algebras by twisting. To illustrate it, consider a
cocommutative Kac algebra structure on the group von Neumann
algebra $M=\mathcal{L}(G)$ of a non commutative l.c. group $G$
with comultiplication $\de(\lambda_g)= \lambda_g\ot\lambda_g$
($\lambda_g$ is the left translation by $g\in G$). Let us
construct on $M$ another (in general, non cocommutative) Kac
algebra structure with comultiplication $\de_\Omega(\cdot) =
\Omega\de(\cdot)\Omega^*$, where $\Omega\in M\ot M$ is
a unitary verifying certain 2-cocycle condition. In order to find
such an $\Omega$, let us, following to M. Rieffel \cite{Rief} and
M. Landstad \cite{Land}, take an inclusion $\alpha :
L^\infty({\hat K})\to M$, where $\hat K$ is the dual to some
abelian subgroup $K$ of $G$ such that $\delta\vert_{K}=1$
($\delta(\cdot)$ is the module of $G$). Then, one lifts a usual
2-cocycle $\Psi$ of $\hat K:\ \Omega=(\alpha\ot\alpha)\Psi$. The
main result of \cite{Enock} is that Haar measure on $\mathcal{L}(G)$
gives also the Haar measure of the deformed object.

Even though a series of non-trivial Kac algebras was constructed
in this way, the above mentioned "unimodularity" condition on $K$ was
restrictive. Here we develop the twisting construction for \lc{}
quantum groups without this condition and compute explicitly the
deformed Haar measure. Thus, we are able to construct l.c. quantum
groups which are not Kac algebras and to deform objects which
are already non-trivial, for example, the $az+b$ quantum group
\cite{Wor}, \cite{Van}.

A dual construction that we call Rieffel's deformation of a \lc{}
group has been proposed in \cite{Rief}, \cite{RiefQG}, and \cite{Land},
where, using a bicharacter on an abelian subgroup, one deforms the
algebra of functions on a group. This construction has been recently
developed by  Kasprzak \cite{Kas} who showed that the dual comultiplication
is exactly the twisted comultiplication of $\mathcal{L}(G)$. Unfortunately,
a trace that he constructed on the deformed algebra is invariant only under
the above mentioned "unimodularity" condition. In this paper we construct
Rieffel's deformation of \lc{} quantum groups without this condition and
compute the corresponding left invariant weight. This proves, in particular,
the existence of invariant weights on the classical Rieffel's deformation.
We also establish the duality between twisting and the Rieffel's deformation.

The structure of the paper is as follows. First, we recall some preliminary
definitions and give our main results. In Section 3 we develop the twisting
construction for l.c. quantum groups. Section 4 is devoted to the Rieffel's
deformations of \lc{} quantum groups and to the proof of the duality theorem.
In Section 5 we present examples obtained by the two constructions: 1) from
group von Neumann algebras $\mathcal{L}(G)$, in particular, when $G$ is the
$az+b$ group; 2) from the $az+b$ quantum group. Some useful technical results
are collected in Appendix.

{\bf Acknowledgment.} We are grateful to Stefaan Vaes who suggested
how twisting can deform the Haar measure and helped us in the proof
of Proposition \ref{IsoDual}.

\section{Preliminaries and main results}

\begin{subsection}{Notations.}

Let us denote by $B(H)$ the algebra of all bounded linear
operators on a Hilbert space $H$, by $\ot$ the tensor product of
Hilbert spaces or von Neumann algebras and by $\Sigma$ (resp.,
$\sigma$) the flip map on it. If $H, K$ and $L$ are Hilbert spaces
and $X \in B(H \ot L)$ (resp., $X \in B(H \ot K), X \in B(K \ot
L)$), we denote by $X_{13}$ (resp., $X_{12},\ X_{23}$) the
operator $(1 \ot \Sigma^*)(X \ot 1)(1 \ot \Sigma)$ (resp., $X\ot
1,\ 1\ot X$) defined on $H \ot K \ot L$. The identity map will be denoted by $\iota$.

Given a {\it normal semi-finite faithful} (n.s.f.) weight $\theta$
on a von Neumann algebra $M$, we denote: ${\mathcal M}^+_\theta =
\{ x \in M^+ \mid \theta(x) < + \infty \},\ {\mathcal N}_\theta =
\{ x \in M \mid x^*x \in M^+_\theta \},$ and ${\mathcal
M}_\theta=span\ {\mathcal M}^+_\theta$ All l.c. groups considered
in this paper are supposed to be second countable, all Hilbert spaces separable and all von Neumann algebras with separable predual.
\end{subsection}

\subsection{Locally compact quantum groups \cite{VaesLC}, \cite{KV2}}

A pair $(M,\de)$ is called a (von Neumann algebraic) l.c.\ quantum
group  when
\begin{itemize}
\item $M$ is a von Neumann algebra and $\de : M \to M \ot M$ is a
normal and unital $*$-homomorphism which is coassociative: $(\de
\ot \id)\de = (\id \ot \de)\de$
\item There exist n.s.f. weights $\varphi$ and
$\psi$ on $M$ such that
\begin{itemize}
\item $\varphi$ is left invariant in the sense that $\varphi
\bigl( (\om \ot \id)\de(x) \bigr) = \varphi(x) \om(1)$ for all $x
\in {\mathcal M}_{\varphi}^+$ and $\om \in M_*^+$, \item $\psi$ is
right invariant in the sense that $\psi \bigl( (\id \ot \om)\de(x)
\bigr) = \psi(x) \om(1)$ for all $x \in {\mathcal M}_{\psi}^+$ and
$\om \in M_*^+$.
\end{itemize}
\end{itemize}
Left and  right invariant weights are unique up to a positive
scalar.

Represent $M$ on the G.N.S. Hilbert space $H$ of $\varphi$ and define a
unitary $W$ on $H \ot H$:
$$
W^* (\Lambda(a) \ot \Lambda(b)) = (\Lambda \ot \Lambda)(\de(b)(a
\ot 1)) \quad\text{for all}\; a,b \in N_{\varphi}\; .
$$
Here, $\Lambda$ denotes the canonical G.N.S.-map for $\varphi$,
$\Lambda \ot \Lambda$ the similar map for $\varphi \ot \varphi$.
One proves that $W$ satisfies the {\it pentagonal equation}:
$W_{12} W_{13} W_{23} = W_{23} W_{12}$, and we say that $W$ is a
{\it multiplicative unitary}. The von Neumann algebra $M$ and the
comultiplication on it can be given in terms of $W$ respectively
as
$$M = \{ (\id \ot \om)(W) \mid \om \in B(H)_* \}^{-\sigma-strong*} \;
$$
and $\de(x) = W^* (1 \ot x) W$, for all $x \in M$. Next, the l.c.\
quantum group $(M,\de)$ has an antipode $S$, which is the unique
$\sigma$-strongly* closed linear map from $M$ to $M$ satisfying
$(\id \ot \om)(W) \in {\mathcal D}(S)$ for all $\om \in B(H)_*$
and $S(\id \ot \om)(W) = (\id \ot \om)(W^*)$ and such that the
elements $(\id \ot \om)(W)$ form a $\sigma$-strong* core for $S$.
$S$ has a polar decomposition $S = R \tau_{-i/2}$, where $R$ is an
anti-automorphism of $M$ and $\tau_t$ is a one-parameter group of
automorphisms of $M$. We call $R$ the
unitary antipode and $\tau_t$ the scaling group of $(M,\de)$. We
have $\sigma (R \ot R) \de = \de R$, so $\varphi R$ is a right
invariant weight on $(M,\de)$, and we take $\psi:= \varphi R$.

There exist a unique number $\nu > 0$ and a unique positif self-adjoint operator $\sde_M$ affiliated to $M$, such that $[D\psi\,:\,D\varphi]_{t}=\nu^{\frac{it^{2}}{2}}\delta_{M}^{it}$. $\nu$ is the scaling constant of $(M,\Delta)$ and $\delta_{M}$ is the modular element of $(M,\Delta)$. The scaling constant
can be characterized as well by the relative invariance property $\varphi
\, \tau_t = \nu^{-t} \, \varphi$.

For the dual l.c.\ quantum group $(\Mh,\deh)$ we have
$$\Mh = \{(\om \ot \id)(W) \mid \om \in B(H)_*
\}^{-\sigma-strong*}$$ and $\deh(x) = \Sigma W (x \ot 1) W^*
\Sigma$ for all $x \in \Mh$. Turn the predual $M_*$ into a Banach
algebra with the product $\om \, \mu = (\om \ot \mu)\de$ and define
$$\lambda: M_* \to \Mh : \lambda(\om) = (\om \ot \id)(W),$$ then
$\lambda$ is a homomorphism and $\lambda(M_*)$ is a
$\sigma$-strongly* dense subalgebra of $\Mh$. A left invariant
n.s.f. weight $\hat \varphi$ on $\Mh$ can be constructed
explicitly. Let $\mathcal{I}=\left\{\omega\in M_{*}\,|\,\exists C\geq 0,\,|\omega(x^{*})|\leq C||\Lambda(x)||\,\forall x\in\mathcal{N}_{\varphi}\right\}$. Then $(H,\iota,\hat{\Lambda})$ is the G.N.S. construction for $\hat{\varphi}$ where $\lambda(\mathcal{I})$ is a $\sigma$-strong-*-norm core for $\hat{\Lambda}$ and $\hat{\Lambda}(\lambda(\omega))$ is the unique vector $\xi(\omega)$ in $H$ such that
$$\langle \xi(\omega),\Lambda(x)\rangle=\omega(x^{*}).$$
The multiplicative unitary of $(\Mh,\deh)$ is $\hat{W} =
\Sigma W^* \Sigma$.

Since $(\Mh,\deh)$ is again a \lc{} quantum group, denote its
antipode by $\hat{S}$, its unitary antipode by $\hat{R}$ and its
scaling group by $\hat{\tau}_t$. Then we can construct the dual of
$(\Mh,\deh)$, starting from the left invariant weight
$\hat\varphi$. The bidual \lc{} quantum group $(\hat\Mh,\hat\deh)$
is isomorphic to $(M,\de)$. Denote by $\hat\sigma_t$ the modular
automorphism group of the weight $\hat\varphi$. The modular
conjugations of the weights $\varphi$ and $\hat\varphi$ will be
denoted by $J$ and $\hat J$ respectively. Let us mention that
$ R(x) =\hat J x^* \hat J$, for all $x \in M$,
and $\hat R(y) = J y^* J$, for all $y \in \Mh$ .

$(M,\de)$ is a {\it Kac algebra} (see \cite{ES}) if and only if
$\tau_t=\id$ and $\sde_M$ is affiliated to the center
of $M$. In particular, $(M,\de)$ is a Kac algebra if $M$ is
commutative. Then $(M,\de)$ is generated by a usual l.c.\ group $G
: M=L^{\infty}(G), (\de_G f)(g,h) = f(gh),$ $(S_Gf)(g) =
f(g^{-1}),\ \varphi_G(f)=\int f(g)\; dg$, where $f\in
L^{\infty}(G),\ g,h\in G$ and we integrate with respect to the
left Haar measure $dg$ on $G$. Then $\psi_G$ is given by
$\psi_G(f) = \int f(g^{-1}) \; dg$ and $\sde_M$ by the strictly
positive function $g \mapsto \sde_G(g)^{-1}$.

$L^\infty(G)$ acts on $H=L^2(G)$ by multiplication and
$(W_G\xi)(g,h)=\xi(g,g^{-1}h),$ for all $\xi\in H\ot H=L^2(G\times
G)$. Then $\Mh=\mathcal L(G)$ is the group von Neumann algebra
generated by the left translations $(\lambda_g)_{g\in G}$ of $G$
and $\deh_G(\lambda_g)= \lambda_g\ot\lambda_g$. Clearly,
$\deh_G^{op}:=\sigma\circ\deh_G=\deh_G$, so $\deh_G$ is
cocommutative. Every cocommutative \lc{} quantum group is obtained
in this way.

\subsection{$q$-commuting pair of operators \cite{WorE2}}

We will use the following notion of commutation relations between unbounded operators. Let $(T,S)$ be a pair of closed operators acting on a Hilbert space $H$. Suppose that $\text{Ker}(T)=\text{Ker}(S)=\{0\}$ and denote by $S=\text{Ph}(S)|S|$ and $T=\text{Ph}(T)|T|$ the polar decompositions. Let $q>0$. We say that $(T,S)$ is a \textit{$q$-commuting pair} and we denote it by $TS=ST$, $TS^{*}=q^{2}S^{*}T$ if the following conditions are satisied
\begin{enumerate}
\item $\text{Ph}(T)\text{Ph}(S)=\text{Ph}(S)\text{Ph}(T)$ and $|T|$ and $|S|$ strongly commute.
\item $\text{Ph}(T)|S|\text{Ph}(T)^{*}=q|S|$ and $\text{Ph}(S)|T|\text{Ph}(S)^{*}=q|T|$.
\end{enumerate}If $T$ and $S$ are q-commuting and normal operators then the product $TS$ is closable and its closure, always denoted by $TS$ has the following polar decomposition $\text{Ph}(TS)=\text{Ph}(T)\text{Ph}(S)$ and $|TS|=q^{-1}|T||S|$.

\subsection{The quantum $az+b$ group \cite{Wor}, \cite{Van}}

Let us describe an explicit example of \lc{} quantum group. Let $s$
and $m$ be two operators defined on the canonical basis
$(e_{k})_{k\in\mathbb{Z}}$ of $l^{2}(\mathbb{Z})$ by $se_{k}=
e_{k+1}$ and $me_{k}=q^{k}e_{k}$ ($0< q <1$). The
G.N.S. space of the quantum $az+b$ group is $H=
l^{2}\left(\mathbb{Z}^{4}\right)$, where we define the
operators
$$a=m\otimes s^{*}\otimes 1\otimes s\quad\text{and}\quad
b=s\otimes m\otimes s\otimes 1$$ with polar decompositions
$a=u|a|$ and $b=v|b|$ given by
$$
\begin{array}{c}
|a|=m\otimes 1\otimes 1\otimes 1\quad\text{and}\quad
u=1\otimes s^{*}\otimes 1\otimes s\\
|b|=1\otimes m\otimes 1\otimes 1\quad\text{and}\quad v=s\otimes
1\otimes s\otimes 1.
\end{array}
$$
Then $u|b|=q|b|u,\ |a|v=qv|a|$, this is the meaning of
the relations $ab=q^{2}ba$ and $ab^{*}=b^{*}a$. Also $\mathrm{Sp}(|a|)=\mathrm{Sp}(|b|)=\mathrm{Sp}(m)=q^{\mathbb{Z}}\cup \{0\},\
\mathrm{Sp}(u)=\mathrm{Sp}(v)=\mathbb{S}^{1}$, where $\mathrm{Sp}$
means the spectrum. Thus,
$\mathrm{Sp}(a)=\mathrm{Sp}(b)=\mathbb{C}^{q}\cup\{0\}$, where $\mathbb{C}^{q}=\left\{z\in\mathbb{C},\,\,\,|z|\in
q^{\mathbb{Z}}\right\}$. The von Neumann algebra of the quantum
$az+b$ group is
$$M:=\left\{\text{finite sums}\\ \sum_{k,l}f_{k,l}(|a|,|b|)v^{k}u^{l}\quad
\text{for}\,\,\,f_{k,l} \in L^{\infty}\left(q^{\mathbb{Z}}\times
q^{\mathbb{Z}}\right)\right\}^{''}.$$ Consider the following
version of the quantum exponential function on $\mathbb{C}^{q}$:
$$
F_{q}(z)=\prod_{k=0}^{+\infty}\frac{1+q^{2k}\overline{z}}{1+q^{2k}z}.
$$
The fundamental unitary of the $az+b$ quantum group is $W=\Sigma V^{*}$ where
$$
V=F_{q}(\hat{b}\otimes b)\chi(\hat{a}\otimes 1, 1\otimes a),
$$
and $\chi(q^{k+i\varphi},q^{l+i\psi})=q^{i(l\varphi+k\psi)}$ is
a bicharacter on $\mathbb{C}^q$. The comultiplication is then
given on generators by
$$W^{*}(1\otimes a)W=a\otimes a\quad\text{and}\quad W^{*}(1\otimes
b)W^{*}=a\otimes b\dot{+}b\otimes 1,$$
where $\dot{+}$ means the closure of the sum. The left invariant weight is
$$\varphi(x)=\sum_{i,j}q^{2(j-i)}f_{0,0}(q^{i},q^{j}),\quad\text{where}\quad
x=\sum_{k,l}f_{k,l}(|a|,|b|)v^{k}u^{l}.$$ The G.N.S. construction
for $\varphi$ is given by $(H,\iota,\Lambda),$ where
$$\Lambda(x)=\sum_{k,l}q^{k+l}\xi_{k,l}\otimes e_{k}\otimes e_{l}\quad
\text{with}\quad\xi_{k,l}(i,j)=q^{j-i}f_{k,l}(q^{i},q^{j}).$$ The
ingredients of the modular theory of $\varphi$ are
$$
\begin{array}{c}
J(e_{r}\otimes e_{s}\otimes e_{k}\otimes e_{l})
=e_{r-k}\otimes e_{s+l}\otimes e_{-k}\otimes e_{-l},\\
\nabla=1\otimes 1\otimes m^{-2}\otimes
m^{-2},
\end{array}
$$ so $\sigma_{t}(a)=q^{-2it}a$ and
$\sigma_{t}^{'}(b)=b$, and the modular element is
$\delta=|a|^{2}$.

The dual von Neumann algebra is
$$\widehat{M}:=\left\{\text{finite sums}\\ \sum_{k,l}f_{k,l}(|\hat{a}|,|\hat{b}|)
\hat{v}^{k}\hat{u}^{l}\quad\text{for}\,\,\,f_{k,l}
\in L^{\infty}\left(q^{\mathbb{Z}}\times
q^{\mathbb{Z}}\right)\right\}^{''}.$$ Here
$\hat{a}=\hat{u}|\hat{a}|$ and $\hat{b}=\hat{v}|\hat{b}|$ are the
polar decompostions of the operators
$$\hat{a}=s^{*}\otimes 1\otimes 1\otimes m,\quad\hat{b}=
s^{*}m\otimes\left(-m^{-1}\otimes m^{-1}s^{*}+m^{-1}s^{*}\otimes
s^{*}\right)\otimes s.$$ The formulas for the dual
comultiplication and the dual left invariant weight are the same,
but this time in terms of $\hat{a}$ and $\hat{b}$.

\subsection{One-parameter groups of automorphisms of von Neumann
algebras}

Consider a von Neumann algebra $M\subset \mathcal{B}(H)$ and a
continuous group homomorphism
$\sigma\,:\,\mathbb{R}\rightarrow\text{Aut}(M)$,
$t\mapsto\sigma_{t}$. There is a standard way to construct, for every
$z\in\mathbb{C}$, a strongly closed densely defined linear
multiplicative in $z$ operator $\sigma_{z}$ in $M$. Let $\mathcal{S}(z)$
be the strip
$\{y\in\mathbb{C}\,|\,\,\text{Im}(y)\in[0,\text{Im}(z)]\}$. Then
we define :
\begin{itemize}
\item The domain $D(\sigma_{z})$ is the set of such elements $x$ in $M$
that the map $t\mapsto\sigma_{t}(x)$ has a strongly continuous
extension to $\mathcal{S}(z)$ analytic on $\mathcal{S}(z)^{0}$.
\item Consider $x$ in $D(\sigma_{z})$ and $f$ the unique extension of the map
$t\mapsto\sigma_{t}(x)$ strongly continuous on $\mathcal{S}(z)$
and analytic on $\mathcal{S}(z)^{0}$. Then, by definition,
$\sigma_{z}(x)=f(z)$.
\end{itemize}
If $x$ is not in $D(\sigma_{z})$, we define an
unbounded operator $\sigma_{z}(x)$ on $H$ as follows:
\begin{itemize}
\item The domain $D(\sigma_{z}(x))$ is the set of such $\xi\in H$
that the map $t\mapsto\sigma_{t}(x)\xi$ has a continuous and bounded
extension to $\mathcal{S}(z)$ analytic on $\mathcal{S}(z)^{0}$.
\item Consider $\xi$ in $D(\sigma_{z}(x))$ and $f$ the unique extension of the map
$t\mapsto\sigma_{t}(x)\xi$ continuous and bounded on $\mathcal{S}(z)$, and
analytic on $\mathcal{S}(z)^{0}$. Then, by definition,
$\sigma_{z}(x)\xi=f(z)$.
\end{itemize}

Let $x$ in $M$, then it is easily seen that the following element
is analytic
$$x(n):=\sqrt{\frac{n}{\pi}}\int_{-\infty}^{+\infty}e^{-nt^{2}}\sigma_{t}(x)dt.$$
The following lemma is a standard exercise:
\begin{lemma}\label{group}
\begin{enumerate}
\item $x(n)\rightarrow x$ $\sigma$-strongly-* and if
$\xi\in\mathcal{D}(\sigma_{z}(x))$ we have
$\sigma_{z}(x(n))\xi\rightarrow\sigma_{z}(x)\xi$. \item Let
$X\subset M$ be a strongly-* dense subspace of $M$ then the set
$\{x(n),\,n\in\mathbb{N},\,x\in X\}$ is a $\sigma$-strong-* core
for $\sigma_{z}$.
\end{enumerate}
\end{lemma}

\begin{proposition}\label{Propext}
Let $A$ be a positive self-adjoint operator affiliated with $M$ and $u$ a unitary in $M$ commuting with $A$ such that $\sigma_{t}(u)=u A^{it}$ for all $t\in\R$,
then $\sigma_{-\frac{i}{2}}(u)$ is a normal operator affiliated with $M$ and its polar decomposition is $ \sigma_{-\frac{i}{2}}(u)=u A^{\frac{1}{2}}.$
\end{proposition}

\begin{proof}
Let $\alpha\in\R$ and
$\mathcal{D}_{\alpha}$ the horizontal strip bounded by $\R$ and $\R-i\alpha$.
Let $\xi\in\mathcal{D}(A^{\frac{1}{2}})$. There exists a continuous bounded
extension $F$ of $t\mapsto A^{it}\xi$ on $\mathcal{D}_{\frac{1}{2}}$ analytic on
$\mathcal{D}_{\frac{1}{2}}^{0}$ (see Lemma 2.3 in \cite{Tak2}). Define
$G(z)=uF(z)$. Then $G(z)$ is continuous and bounded on $\mathcal{S}(-\frac{i}{2})=\mathcal{D}_{\frac{1}{2}}$, and analytic on $\mathcal{S}(-\frac{i}{2})^{0}$. Moreover, $G(t)=u F(t)=u A^{it}\xi=\sigma_{t}(u)\xi$,
so $\xi\in\mathcal{D}(\sigma_{-\frac{i}{2}}(u))$ and $\sigma_{-\frac{i}{2}}(u)\xi=G(-\frac{i}{2})=u A^{\frac{1}{2}}\xi$. Then
$u A^{\frac{1}{2}}\subset\sigma_{-\frac{i}{2}}(u)$.
The other inclusion is proved in the same way.
\end{proof}

\subsection{The Vaes' weight}\label{sectionVaes}

Let $M\subset\mathcal{B}(H)$ be a von Neumann algebra with a
\nsf{} weight $\varphi$ such that $(H,\iota,\Lambda)$ is the \GNS{}
construction for $\varphi$. Let $\nabla$, $\sigma_{t}$ and $J$ be
the objects of the modular theory for $\fie$, and $\delta$ a positive
self-adjoint operator affiliated with $M$ verifying  $\sigma_{t}(\delta^{is})= \lambda^{ist}\delta^{is}$, for all $s,t\in\R$ and some $\lambda>0$.

\begin{lemma}\cite{VaesRN}\label{LemVaesRN}
There exists a sequence of self-adjoint elements $e_{n}\in M$, analytic w.r.t.
$\sigma$ and commuting with any operator that commutes with $\delta$, and such that,
for all $x,z\in\Cset$, $\delta^{x}\sigma_{z}(e_{n})$ is bounded with domain $H$,
analytic w.r.t. $\sigma$ and satisfying $\sigma_{t}(\delta^{x}\sigma_{z}(e_{n}))=\delta^{x}\sigma_{t+z}(e_{n})$, and
$\sigma_{z}(e_{n})$ is a bounded sequence which converges $*$-strongly to $1$,
for all $z\in\Cset$. Moreover, the function $(x,z)\mapsto\delta^{x}\sigma_{z}(e_{n})$
is analytic from $\Cset^{2}$ to $M$.
\end{lemma}
Let $N=\left\{a\in M,\,a\delta^{\frac{1}{2}}\,\,\text{is bounded and}\,\,\overline{a\delta^{\frac{1}{2}}}
\in\mathcal{N}_{\varphi}\right\}$. This is an ideal $\sigma$-strongly$*$ dense in $M$ and the map $a\mapsto\Lambda(\overline{a\delta^{\frac{1}{2}}})$ is $\sigma$-strong$*$-norm closable; its closure will be denoted by $\Lambda_{\delta}$.
\begin{proposition} \cite{VaesRN}
There exists a unique \nsf{} weight $\varphi_{\delta}$ on $M$ such that $(H,\iota,\Lambda_{\delta})$ is a \gns{} construction for $\varphi_{\delta}$. Moreover,
\begin{itemize}
\item the objects of the modular theory of $\fie_{\delta}$ are $J_{\delta}=\lambda^{\frac{i}{4}}J$ and $\nabla_{\delta}=J\delta^{-1}J\delta\nabla$,
\item $[D\fie_{\delta}:D\fie]_{t}=\lambda^{i\frac{t^{2}}{2}}\delta^{it}.$
\end{itemize}
\end{proposition}

\subsection{Main results}\label{SectionResult}

Let $(M,\Delta)$ be a \lc{} quantum group with left and right invariant
weights $\varphi$ and $\psi=\varphi\circ R$, and the corresponding modular
groups $\sigma$ and $\sigma^{'}$. Let $\Omega\in M\otimes M$
be a \textit{2-cocycle}, \ie{}, a unitary such that $(\Omega\otimes 1) (\Delta\otimes\iota)(\Omega)=(1\otimes\Omega)(\iota\otimes\Delta)(\Omega)$.
Then obviously $\Delta_{\Omega}=\Omega\Delta(.)\Omega^{*}$ is
a comultiplication on $M$. If $(M,\Delta)$ is discrete quantum group and
$\Omega$ is any $2$-cocyle on $(M,\Delta)$, then $(M,\Delta_{\Omega})$ is
again a discrete quantum group \cite{BichonVaes}. If $(M,\Delta)$ is not
discrete, it is not known, in general, if $(M,\Delta_{\Omega})$ is a \lc{}
quantum group. Let us consider the following special construction of $\Omega$.
Let $G$ be \lc{} group and $\alpha$ be a unital normal faithful *-homomorphism
from $L^{\infty}(G)$ to $M$ such that $\alpha\otimes\alpha\circ\Delta_{G}=
\Delta\circ\alpha$. In this case we say that $G$ is a \textit{co-subgroup}
of $(M,\Delta)$, and we write $\widehat{G}<(M,\Delta)$. Then the von
Neumann algebraic version of Proposition 5.45 in \cite{VaesLC} gives
$$\tau_{t}\circ\alpha=\alpha\quad\text{and}\quad
R\circ\alpha(F)=\alpha(F(\cdot^{-1})), \,\,\forall\,F\in
L^{\infty}(G).$$

Let $\Psi$ be a continuous bicharacter on $G$. Then $\Omega=(\alpha\otimes\alpha)(\Psi)$
is a $2$-cocycle on $(M,\Delta)$. In \cite{Enock} it was supposed that $\sigma_{t}$ acts trivially on the image of $\alpha$ and it was shown that in this case, $\varphi$ is also $\Delta_{\Omega}$-left invariant. Here we suppose that $\sigma_{t}$ acts by translations,
i.e., that there exists a continuous group homomorphism $t\mapsto \gamma_{t}$ from $\R$
to $G$ such that $\sigma_{t}(\alpha(F))=\alpha(F(.\gamma_{t}^{-1}))$. In this case we say
that the co-subgroup $G$ is \textit{stable}. Then $\sigma_{t}^{'}$ also acts by translations:
\begin{equation}\label{stable}
\sigma_{t}^{'}\circ\alpha(F)=R\circ\sigma_{-t}\circ R\circ\alpha(F)
=\alpha(F(\cdot\gamma_{t}^{-1}))=\sigma_{t}\circ\alpha(F).
\end{equation}
In particular, $\delta^{it}\alpha(F)=\alpha(F)\delta^{it}, \forall\ t\in\mathbb{R},
\ F\in L^{\infty}(G)$. In our case $\varphi$ is not necessarily $\Delta_{\Omega}$-left
invariant, and one has to construct another weight on $M$. Note that $(t,s)\mapsto\Psi(\gamma_{t},\gamma_{s})$ is a bicharacter on $\R$. Thus, there exists $\lambda>0$ such that $\Psi(\gamma_{t},\gamma_{s})=\lambda^{its}$ for all $s,t\in\R$.
Let us define the following unitaries in $M$:
$$u_{t}=\lambda^{i\frac{t^{2}}{2}}\alpha\left(\Psi(\cdot,\gamma_{t}^{-1})\right)\quad
\text{and}\quad
v_{t}=\lambda^{i\frac{t^{2}}{2}}\alpha\left(\Psi(\gamma_{t}^{-1},\cdot)\right).$$
Then equation $(\ref{stable})$ and the definition of a bicharacter imply that $u_{t}$
is a $\sigma$-cocycle and $v_{t}$ is a $\sigma^{'}$-cocycle. The converse of the Connes' Theorem gives then \nsf{} weights $\varphi_{\Omega}$ and $\psi_{\Omega}$ on $M$ such that:
$$u_{t}=[D\varphi_{\Omega}\,:\,D\varphi]_{t}\quad\text{and}\quad
v_{t}=[D\psi_{\Omega}\,:\,D\psi]_{t}.$$
The main result of Section \ref{SectionTwisting} is the following. We denote by $W$ the multiplicative unitary of $(M,\Delta)$, and put $W_{\Omega}^{*}=\Omega(\hat{J}\otimes
J)W\tilde{\Omega}(\hat{J}\otimes J)$.

\begin{theorem}\label{ThmTwist}
$(M,\Delta_{\Omega})$ is a \lc{} quantum group :
\begin{itemize}
\item $\varphi_{\Omega}$ is left invariant, \item $\psi_{\Omega}$
is right invariant, \item $W_{\Omega}$ is the fundamental
multiplicative unitary.
\item The scaling group and the scaling constant are $\tau_{t}^{\Omega}=\tau_{t}$, $\nu_{\Omega}=\nu$.
\end{itemize}
\end{theorem}

If $G$ is abelian, we compute explicitly the modular element and the antipode.

In section \ref{SectionDual} we construct the Rieffel's deformation of a \lc{}
quantum group with an abelian stable co-subgroup $\widehat{G}<(M,\Delta)$ and
prove that this construction is dual to the twisting. Switching to the
additive notations for $G$, define $L_{\gamma}=\alpha(u_{\gamma})$ and $R_{\gamma}
=J L_{\gamma}J$, where $\gamma\in\hat{G},\ u_{\gamma}=\langle\gamma,g\rangle\in
L^{\infty}(G)$, and $J$ is the modular conjugation of $\varphi$. Then Proposition \ref{PropCo-Subgroups} shows that $\widehat{G}^{2}$
acts on $\widehat{M}$ by conjugation by the unitaries $L_{\gamma_{1}}R_{\gamma_{2}}$.
We call this action the \textit{left-right action}. Let $N=\widehat{G}^{2}\ltimes
\widehat{M}$ be the crossed product von Neumann algebra generated by $\lambda_{\gamma_{1},\gamma_{2}}$
and $\pi(x)$, where $\gamma_{i}\in\widehat{G}$ and $x\in\widehat{M}$, and
let $\theta$ be the dual action of $G^{2}$ on $N$. We show that there exists a
unique unital normal *-homomorphism $\Gamma$ from $N$ to $N\otimes N$ such that
$\Gamma(\lambda_{\gamma_{1},\gamma_{2}})=\lambda_{\gamma_{1},0}
\otimes\lambda_{0,\gamma_{2}}$ and $\Gamma(\pi(x))=(\pi\otimes\pi)\hat{\Delta}(x)$.
Let $\Psi$ be a continuous bicharacter
on $G$. Note that, for all $g\in G$, we have $\Psi_{g}\in\widehat{G}$, where $\Psi_{g}(h)=\Psi(h,g)$. We denote by $\theta^{\Psi}$ the \textit{twisted
dual action} of $G^{2}$ on $N$:
\begin{equation}\label{TwistedDualAction}
\theta^{\Psi}_{(g_{1},g_{2})}(x)=\lambda_{\Psi_{g_{1}},\Psi_{g_{2}}}
\theta_{(g_{1},g_{2})}(x)\lambda_{\Psi_{g_{1}},\Psi_{g_{2}}}^{*},
\quad\text{for any}\,\,g_{1},g_{2}\in G,\,\,x\in \widehat{G}^{2}\ltimes\widehat{M},
\end{equation}
and by $N_{\Omega}$ the fixed point algebra under this action (we would
like to point out that $N_\Omega$ is not a deformation of $N$, it is just a
fixed point algebra with respect to the action $\theta^{\Psi}$ related to $\Omega$).
Put $\Upsilon=(\lambda_{R}\otimes\lambda_{L})(\tilde{\Psi}^{*})\in
N\otimes N$, where $\lambda_{R}$ and $\lambda_{L}$ are the unique unital
normal *-homomorphisms from $L^{\infty}(G)$ to $N$ such that $\lambda_{L}(u_{\gamma})=\lambda_{{\gamma},0}$ and $\lambda_{R}(u_{\gamma})=\lambda_{0,{\gamma}}$, and put $\Gamma_{\Omega}(\cdot)=\Upsilon\Gamma(\cdot)\Upsilon^{*}$. Then we show that
$\Gamma_{\Omega}$ is a comultiplication on $N_{\Omega}$ and construct a left
invariant weight on $(N_{\Omega},\Gamma_{\Omega})$.  Because $\theta^{\Psi}_{g_{1},g_{2}}(\lambda_{\gamma_{1},\gamma_{2}})=
\theta{g_{1},g_{2}}(\lambda_{\gamma_{1},\gamma_{2}})=
\overline{\langle\gamma_{1},g_{1}\rangle}\overline
{\langle\gamma_{2},g_{2}\rangle}\lambda_{\gamma_{1},\gamma_{2}}$, we have a canonical
isomorphism $N=\widehat{G}^{2}\ltimes\widehat{M}\to\widehat{G}^{2}\ltimes N_{\Omega}$
intertwining the twisted dual action on $N$ with the dual action on $\widehat{G}^{2}
\ltimes N_{\Omega}$. Denoting by $\widetilde{\hat{\varphi}}$ the dual weight of
$\hat{\varphi}$ on $N$ and by $\widetilde{\hat{\sigma}}$ its modular group, we show
that $w_{t}=\lambda^{-it^{2}}\lambda_{R}(\Psi(-\gamma_{t},.))$ is a $\widetilde{\hat{\sigma}}$-cocycle. This implies the existence of a unique \nsf{} $\tilde{\mu}_{\Omega}$ on $N$ such that $w_{t}=[D\tilde{\mu}_{\Omega}\,:\,D\tilde{\hat{\varphi}}]_{t}$. Moreover, one can
show that $\tilde{\mu}_{\Omega}$ is $\theta^{\Psi}$ invariant. Thus, there exists a
unique \nsf{} $\mu_{\Omega}$ on $N_{\Omega}$ such that the dual weight of $\mu_{\Omega}$
is $\tilde{\mu}_{\Omega}$. In order to formulate the main result of Section \ref{SectionDual},
let us denote by $(\widehat{M}_{\Omega},\hat{\Delta}_{\Omega})$ the dual of $(M,\Delta_{\Omega})$.

\begin{theorem}\label{TheoremRieffel}
$(N_{\Omega},\Gamma_{\Omega})$ is a \lc{} quantum group and $\mu_{\Omega}$ is left
invariant. Moreover there is a canonical isomorphism $(N_{\Omega},\Gamma_{\Omega})\simeq (\widehat{M}_{\Omega},\hat{\Delta}_{\Omega})$.
\end{theorem}

Note that the Rieffel's deformation in the $\cs$-setting was constructed by the first
author in \cite{Fima}, see also Remark \ref{RemarkCstar} and \cite{FimaVain} for an
overview.

In Section \ref{SectionExample} we calculate explicitly two examples. It is known
that if $H$ is an abelian closed subgroup of a \lc{} group $G$, then there is a unique
faithful unital normal *-homomorphism $\alpha$ from $L^{\infty}(\widehat{H})$ to $\mathcal{L}(G)$ such that $\alpha(u_{h})=\lambda_{G}(h)$, for all $h\in\widehat{H}$,
where $\lambda_{G}$ is the left regular representation of $G$, so $H<(\mathcal{L}(G),\hat{\Delta}_{G})$ is a co-subgroup. The left (and right) invariant
weight on $\mathcal{L}(G)$ is the Plancherel weight for which $\sigma_{t}(\lambda_{g})=\delta_{G}^{it}(g)\lambda_{g}$, for all $g\in G$,
where $\delta_{G}$ is the modular function of $G$. Then  $\sigma_{t}\circ\alpha(u_{g})=\alpha(u_{g}(\cdot-\gamma_{t}))$,
where $\gamma_{t}$ is the character on $K$ defined by
$\langle\gamma_{t},g\rangle=\delta_{G}^{-it}(g)$. Because the vector space spanned
by the $u_{h}$ for $h\in H$ is dense in $L^{\infty}(\widehat{H})$, we
have $\sigma_{t}\circ\alpha(F)=\alpha(F(\cdot-\gamma_{t}))$, for all $F\in L^{\infty}(\hat{K})$. Thus, $H<(\mathcal{L}(G),\hat{\Delta}_{G})$ is stable. So,
given a bicharacter $\Psi$ on $\widehat{H}$, we can perform the twisting construction.
The deformation of the Haar weight will be non trivial when $H$ is not in the kernel
of the modular function of $G$.

Let $G=\C^{*}\ltimes\C$ be the $az+b$ group and $H=\mathbb{C}^{*}$ be the abelian
closed subgroup of elements of the form $(z,0)$ with $z\in\C^{*}$. Identifying $\widehat{\mathbb{C}^{*}}$ with $\mathbb{Z}\times\mathbb{R}_{+}^{*}$:
$$
\mathbb{Z}\times\mathbb{R}_{+}^{*}\rightarrow\widehat{\mathbb{C}^{*}},\quad (n,\rho)\mapsto\gamma_{n,\rho}=(re^{i\theta}\mapsto e^{i\ln{r}\ln{\rho}}e^{in\theta}),
$$
let us define, for all $x\in\R$, the following bicharacter on $\mathbb{Z}\times\mathbb{R}_{+}^{*}$:
$$
\Psi_{x}((n,\rho),(k,r))=e^{ix(k\ln{\rho}-n\ln{r})}
$$
and perform the twisting construction. We obtain a family of \lc{} quantum groups $(M_{x},\Delta_{x})$ with trivial scaling group and scaling constant. Moreover, we
show that the antipode is not deformed. The main result of Section \ref{Sectionaz+bclassique} is the following. Let us denote by $\varphi$ the Plancherel weight on $\mathcal{L}(G)$
and by a subscript $x$ the objects associated with $(M_{x},\Delta_{x})$.

\begin{theorem}\label{Theoremaz+bclassique}
We have:
\begin{itemize}
\item  $[D\varphi_{x}\,:\,D\varphi]_{t}=\lambda^{G}_{(e^{itx},0)}$, $\delta_{x}^{it}= \lambda^{G}_{(e^{-2itx},0)}$.

\item $(M_{-x},\Delta_{-x})\simeq (M_{x},\Delta_{x})^{\text{op}}$ and if $x,y \geq 0$ with $x\neq y$ then $(M_{x},\Delta_{x})$ and $(M_{y},\Delta_{y})$ are not isomorphic.
\end{itemize}
The von Neumann algebra of the dual quantum group $(\widehat{M}_{x},\hat{\Delta}_{x})$ is generated by two operators $\hat{\alpha}$ and $\hat{\beta}$ affiliated with it and such that
\begin{itemize}
\item $\hat{\alpha}$ is normal, $\hat{\beta}$ is \textrm{q-normal}, \ie{}, $\hat{\beta}\hat{\beta}^{*}=q\hat{\beta}^{*}\hat{\beta}$,
\item $\hat{\alpha}\hat{\beta}=\hat{\beta}\hat{\alpha}$ and $\hat{\alpha}\hat{\beta}^{*}=q\hat{\beta}^{*}\hat{\alpha}$, with $q=e^{4x}$.
\end{itemize}
The comultiplication is given by  $\hat{\Delta}_{x}(\hat{\alpha})=\hat{\alpha}\otimes\hat{\alpha}$ and $\hat{\Delta}_{x}(\hat{\beta})=\hat{\alpha}\otimes\hat{\beta}\dot{+}\hat{\beta}\otimes 1$.
\end{theorem}

For the dual $(\widehat{M}_{x},\hat{\Delta}_{x})$ we deform, like in the Woronowicz'
quantum $az+b$ group, the commutativity relation between the two coordinate functions,
but the difference is that we also deform the normality of the second coordinate function.

The second example of Section \ref{SectionExample} is the twisting of an
already non trivial object. Consider the Woronowicz' quantum $az+b$ group $(M,\Delta)$
at a fixed parameter $0<q<1$. Let $\alpha:\ L^{\infty}(\mathbb{C}^{q})\rightarrow M$ be
the normal faithful *-homomorphism defined by $\alpha(F)=F(a).$ Because
$\Delta(a)=a\otimes a$, one has $\Delta\circ\alpha= (\alpha\otimes\alpha)\circ\Delta_{\mathbb{C}^{q}}$. Thus, we have a co-subgroup $\widehat{\C^{q}}<(M,\Delta)$ which is stable:
$$\sigma_{t}\circ\alpha(F)=\sigma_{t}(F(a))=F(\sigma_{t}(a))
=F(q^{-2it}a)=\alpha(F(\cdot\gamma_{t}^{-1})),$$ where
$\gamma_{t}=q^{2it}\in\mathbb{C}^{q}$. Performing the twisting
construction with the bicharacters
$$
\Psi_{x}(q^{k+i\varphi},q^{l+i\psi})=q^{ix(k\psi-l\varphi)}, \
\forall x\in\mathbb{Z},
$$
we get the twisted l.c. quantum groups $(M_{x},\Delta_{x})$. The main result of
Section \ref{Sectionaz+bquantique} is the following. Recall that we denote by $a=u|a|$
the polar decomposition of $a$.

\begin{theorem}\label{Theoremaz+bquantique}
One has $\Delta_{x}(a)=a\otimes a$ and $\Delta_{x}(b)=u^{-x+1}|a|^{x+1}\otimes
b\dot{+}b\otimes u^{x}|a|^{-x}.$ The modular element $\delta_{x}=|a|^{4x+2}$, the
antipode is not deformed and we have $[D\varphi_{x}\,:\,D\varphi]_{t}=|a|^{-2ixt}$. Moreover,
for any $x,y\in\N$, one has: if $x\neq y$, then $(M_{x},\Delta_{x})$ and $(M_{y},
\Delta_{y})$ are not isomorphic;
if $x\neq 0$, then $(M_{x},\Delta_{x})$ and $(M_{-x},\Delta_{-x})$ are not isomorphic.
The von Neumann algebra of the dual quantum group $(\widehat{M}_{x},\hat{\Delta}_{x})$
is generated by two operators $\hat{\alpha}$ and $\hat{\beta}$ affiliated with it and
such that
\begin{itemize}
\item $\hat{\alpha}$ is normal, $\hat{\beta}$ is p-normal, \ie{}, $\hat{\beta}\hat{\beta}^{*}=p\hat{\beta}^{*}\hat{\beta}$,
\item $\hat{\alpha}\hat{\beta}=q^{2}\hat{\beta}\hat{\alpha}$ and $\hat{\alpha}\hat{\beta}^{*}=p\hat{\beta}^{*}\hat{\alpha}$, with $p=q^{-4x}$.
\end{itemize}
The comultiplication is given by $\hat{\Delta}_{x}(\hat{\alpha})=\hat{\alpha}\otimes\hat{\alpha}$ and $
\hat{\Delta}_{x}(\hat{\beta})=\hat{\alpha}\otimes\hat{\beta}\dot{+}\hat{\beta}\otimes
1$.
\end{theorem}

We refer to \cite{FimaVain} for the explicit example of the twisting in the
$\cs$-setting of the group $G$ of $2\times 2$ upper triangular matrices of determinant
$1$ with the abelian subgroup of diagonal matrices in $G$. Next subsection contains useful
technical result.

\subsection{Abelian stable co-subgroups}\label{SectionCo-Subgroups}

Let $\widehat{G}<(M,\Delta)$ be a stable co-subgroup with $G$ abelian. For all $\gamma\in\hat{G}$, the map $t\mapsto\langle\gamma,\gamma_{t}\rangle$ is a character
on $\Rset$, so there exists $\lambda(\gamma)>0$ such that $\langle\gamma,\gamma_{t}\rangle=\lambda(\gamma)^{it}$ for all $t\in\Rset$.

\begin{proposition}\label{PropCo-Subgroups}
Let $\widehat{G}<(M,\Delta)$ be a co-subgroup with $G$ abelian. Then:
\begin{equation}\label{EqWL}
1. (1\otimes L_{\gamma})W(1\otimes L_{\gamma}^{*})=W(L_{\gamma}\otimes1)
,\quad
(1\otimes R_{\gamma})W(1\otimes R_{\gamma}^{*})=(L_{-\gamma}\otimes 1)W,
\end{equation}
for all $\gamma\in\widehat{G}$, so we have two commuting actions $\alpha^{L}$ and $\alpha^{R}$ of $\hat{G}$ on $\hat{M}$:
$\alpha^{L}_{\gamma}(x)=L_{\gamma}x L_{\gamma}^{*}$ and $\alpha^{R}_{\gamma}(x)=R_{\gamma}x R_{\gamma}^{*}$. This gives an action of $\hat{G}^{2}$ on $\hat{M}$ $\alpha_{\gamma_{1},\gamma_{2}}=\alpha^{L}_{\gamma_{1}}\circ\alpha^{R}_{\gamma_{2}}$
such that
\begin{equation}\label{eqactionW}
(\iota\otimes\alpha_{\gamma_{1},\gamma_{2}})(W)=(L_{\gamma_{2}}^{*}\otimes 1)W(L_{\gamma_{1}}^{*}\otimes 1).
\end{equation}
2. If $\widehat{G}<(M,\Delta)$ is stable, then, for all
$x\in\mathcal{N}_{\hat{\varphi}}$ and all $\gamma\in\hat{G}$, we have $\alpha^{L}_{\gamma}(x), \alpha^{R}_{\gamma}(x)\in\mathcal{N}_{\hat{\varphi}}$, $L_{\gamma}\hat{\Lambda}(x)=\hat{\Lambda}(\alpha^{L}_{\gamma}(x))$, and $ R_{\gamma}\hat{\Lambda}(x)=\lambda(\gamma)^{-\frac{1}{2}}
\hat{\Lambda}(\alpha^{R}_{\gamma}(x))$.\label{LemRN1}
\end{proposition}

\begin{proof}
Since $\Delta(L_{\gamma})=L_{\gamma}\otimes L_{\gamma}$, $\Delta(x)=W^{*}(1\otimes x)W$ and $(\hat{J}\otimes J)W(\hat{J}\otimes J)=W^{*}$, it is easy to check the first two equalities. The equality for $\alpha$ follows immediately. To prove the second assertion we need the following

\begin{lemma}[\cite{VaesLC}]\label{awb}
Let $\omega\in\mathcal{I}$, $a\in M$, and $b\in\mathcal{D}(\sigma_{-\frac{i}{2}})$,
then $a\omega b\in\mathcal{I}$ and
$$\xi(a\omega b)=aJ\sigma_{-\frac{i}{2}}(b)^{*}J\xi(\omega).$$
\end{lemma}

Let us prove the second assertion. By the first assertion we have $\alpha_{\gamma}^{L}((\omega\otimes\iota)(W))=(L_{\gamma}\omega\otimes\iota)(W)$. Take $\omega\in\mathcal{I}$. By Lemma \ref{awb}, we have $L_{\gamma}\omega\in\mathcal{I}$ and
$$\hat{\Lambda}(\alpha_{\gamma}^{L}(\lambda(\omega)))=\hat{\Lambda}(\lambda(L_{\gamma}\omega))=L_{\gamma}\hat{\Lambda}(\lambda(\omega)).$$
Because $\lambda(\mathcal{I})$ is a core for $\hat{\Lambda}$, for all $x\in\mathcal{N}_{\hat{\varphi}}$, we have $\alpha_{\gamma}^{L}(x)\in\mathcal{N}_{\hat{\varphi}}$ and
$$\hat{\Lambda}(\alpha_{\gamma}^{L}(x))=L_{\gamma}\hat{\Lambda}(x).$$

By the first assertion, we have $\alpha_{\gamma}^{R}((\omega\otimes\iota)(W))=(\omega L_{-\gamma}\otimes\iota)(W)$.
Note that $\sigma_{t}(L_{\gamma})=\lambda(\gamma)^{-it}L_{\gamma}$, thus $L_{\gamma}\in\mathcal{D}(\sigma_{\frac{i}{2}})$ and $\sigma_{\frac{i}{2}}(L_{\gamma})=\lambda(\gamma)^{\frac{1}{2}}L_{\gamma}$. Take $\omega\in\mathcal{I}$. By Lemma \ref{awb}, we have $\omega L_{-\gamma}\in\mathcal{I}$ and
$$\hat{\Lambda}(\alpha_{\gamma}^{R}(\lambda(\omega)))=\hat{\Lambda}(\lambda(\omega L_{-\gamma}))=\lambda(\gamma)^{\frac{1}{2}}R_{\gamma}\hat{\Lambda}(\lambda(\omega)).$$
Because $\lambda(\mathcal{I})$ is a core for $\hat{\Lambda}$, this concludes the proof.
\end{proof}

\section{Twisting of locally compact quantum groups}\label{SectionTwisting}

Let $G$ be a \lc{} group and $(M,\Delta)$ a \lc{} quantum group. Suppose that $\widehat{G}<(M,\Delta)$ is a stable co-subgroup. We keep the notations of Section \ref{SectionResult}. Note that the maps $(t\mapsto \alpha(\Psi(\cdot,\gamma_{t}^{-1})))$
and $(t\mapsto\alpha(\Psi(\gamma_{s}^{-1},\cdot)))$ are unitary representations of $\R$
in $M$. Let $A$ and $B$ be the positive self-adjoint operators affiliated with $M$ such that $A^{it}=\alpha(\Psi(\cdot,\gamma_{t}^{-1}))$ and $B^{is}=\alpha(\Psi(\gamma_{s}^{-1},\cdot))$. We have $\Delta(A)=A\otimes A$, $\Delta(B)=B\otimes B$. Note that $\sigma_{t}(A^{is})=\alpha(\Psi(\cdot\gamma_{t}^{-1},\gamma_{s}^{-1}))
=\lambda^{ist}A^{is}$.  Also we have
$\sigma_{s}^{'}(B^{it})=\sigma_{s}(B^{it})=\lambda^{ist}B^{it}$, so the weights $\varphi_{\Omega}$ and $\psi_{\Omega}$ are the Vaes' weights associated with $\varphi$,
$\lambda$ and $A$, and with $\psi$, $\lambda$ and $B$, respectively. In the sequel, we
denote by $\Lambda_{\Omega}$ the canonical \GNS{} map associated with $\varphi_{\Omega}$,
and by $F\mapsto\tilde{F}$ the *-automorphism of $L^{\infty}(G\times G)$ defined by $\tilde{F}(g,h)=F(g^{-1},gh)$. Theorem \ref{ThmTwist} is in fact a corollary of the
following result.

\begin{theorem}\label{TwThm1}
For all $x,y\in\mathcal{N}_{\varphi_{\Omega}}$, we have
$\Delta_{\Omega}(x)(y\otimes
1)\in\mathcal{N}_{\varphi_{\Omega}\otimes\varphi_{\Omega}}$ and
$$\left(\Lambda_{\Omega}\otimes\Lambda_{\Omega}\right)\left(\Delta_{\Omega}(x)(y\otimes
1)\right)=
W_{\Omega}^{*}(\Lambda_{\Omega}(y)\otimes\Lambda_{\Omega}(x)),$$
where $W_{\Omega}^{*}=\Omega(\hat{J}\otimes
J)W\tilde{\Omega}(\hat{J}\otimes J)$.
\end{theorem}

\begin{proof}
Let us introduce the sets
$$N=\left\{x\in M,\,xA^{\frac{1}{2}}\,\,\text{is bounded
and}\,\,\overline{xA^{\frac{1}{2}}}\in\mathcal{N}_{\varphi}\right\}\quad\text{and}$$
$$L=\left\{x\in
N,\,A^{-\frac{1}{2}}\overline{xA^{\frac{1}{2}}}\,\,\text{is
bounded and}\,\,
\Lambda(\overline{xA^{\frac{1}{2}}})\in\mathcal{D}(A^{-\frac{1}{2}})\right\}.$$
When $y\in L$, we denote the closure of
$A^{-\frac{1}{2}}\overline{xA^{\frac{1}{2}}}$ by
$\overline{A^{-\frac{1}{2}}yA^{\frac{1}{2}}}$. By definition, $N$ is a $\sigma$-strong$*$-norm core for $\Lambda_\Omega$, and Proposition \ref{PropTx} shows that the same is true for $L$. As $\Lambda_\Omega$ is closed in these topologies, it suffices to prove the theorem for elements $x\in N$ and $y\in L$. The first step is as follows.
\begin{lemma}\label{Twisting1}
Let $x\in N$, $y\in L$ and $F\in
(\alpha\otimes\alpha)(L^{\infty}(G\times G))$. Then
$$(\Delta(x)F^{*}y\otimes 1)(A^{\frac{1}{2}}\otimes
A^{\frac{1}{2}})\quad\text{is bounded and}$$
$$\overline{(\Delta(x)F^{*}(y\otimes 1))(A^{\frac{1}{2}}\otimes
A^{\frac{1}{2}})}
=\Delta(\overline{xA^{\frac{1}{2}}})F^{*}\left(\overline{A^{-\frac{1}{2}}yA^{\frac{1}{2}}}\otimes
1\right).$$
\end{lemma}
\begin{proof}
Note that $\Delta(A^{\frac{1}{2}})=A^{\frac{1}{2}}\otimes
A^{\frac{1}{2}}=W^{*}\left(1\otimes A^{\frac{1}{2}}\right)W$. Let
$x\in N$ and $\xi\in\mathcal{D}(A^{\frac{1}{2}}\otimes
A^{\frac{1}{2}})$. Then $W\xi\in\mathcal{D}(1\otimes
A^{\frac{1}{2}})$ and
\begin{eqnarray*}
\Delta(x)(A^{\frac{1}{2}}\otimes A^{\frac{1}{2}})\xi &
= & W^{*}(1\otimes x)WW^{*}(1\otimes A^{\frac{1}{2}})W\xi \\
& = & W^{*}(1\otimes x)(1\otimes A^{\frac{1}{2}})W\xi \\
& = & W^{*}(1\otimes
\overline{xA^{\frac{1}{2}}})W\xi=\Delta(\overline{xA^{\frac{1}{2}}})\xi.\\
\end{eqnarray*}
Thus, $\Delta(x)(A^{\frac{1}{2}}\otimes
A^{\frac{1}{2}})\subset\Delta(\overline{xA^{\frac{1}{2}}})$ and
because it is densely defined, we have shown that, $\forall x\in
N,$ $\Delta(x)(A^{\frac{1}{2}}\otimes
A^{\frac{1}{2}})$ is bounded
and $\overline{\Delta(x)(A^{\frac{1}{2}}\otimes
A^{\frac{1}{2}})}=\Delta(\overline{xA^{\frac{1}{2}}})$. If $x\in
N,\ y\in N^{'}$, the commutativity of
$(\alpha\otimes\alpha)(L^{\infty}(G\times G))$ implies:
\begin{eqnarray*}
(\Delta(x)F^{*}(y\otimes 1))(A^{\frac{1}{2}}\otimes A^{\frac{1}{2}})
& = & \Delta(x)(1\otimes
A^{\frac{1}{2}})F^{*}(yA^{\frac{1}{2}}\otimes 1)
\\
& = & \Delta(x)(A^{\frac{1}{2}}\otimes
A^{\frac{1}{2}})(A^{-\frac{1}{2}}\otimes 1)F^{*}(yA^{\frac{1}{2}}\otimes 1) \\
& = & \Delta(x)(A^{\frac{1}{2}}\otimes
A^{\frac{1}{2}})F^{*}(A^{-\frac{1}{2}}yA^{\frac{1}{2}}\otimes 1) \\
& \subset &
\Delta(\overline{xA^{\frac{1}{2}}})F^{*}
\left(\overline{A^{-\frac{1}{2}}yA^{\frac{1}{2}}}\otimes 1\right).\\
\end{eqnarray*}
Since $(\Delta(x)F^{*}(y\otimes 1))(A^{\frac{1}{2}}\otimes
A^{\frac{1}{2}})$ is densely defined, the
proof is finished.
\end{proof}
In what follows, we identify $L^{\infty}(G)$ with its image
$\alpha(L^{\infty}(G))$. Note that
$$(\iota\otimes\sigma_{t})(\tilde{F})=\widetilde{(\iota\otimes\sigma_{t})(F)},
\quad\text{for all}\,\,\,t\in\mathbb{R},F\in L^{\infty}(G\times G).$$
By analytic continuation, this is also true for $t=z\in
\mathbb{C}$ and $F\in\mathcal{D}(\iota\otimes\sigma_{z})$.

Now we construct a set of certain elements of
$\mathcal{N}_{\varphi_{\Omega}\otimes\varphi_{\Omega}}$ and give
their images by $\Lambda_{\Omega}\otimes\Lambda_{\Omega}$.

\begin{lemma}\label{element}
Let $x\in N$, $y\in L$ and $F\in L^{\infty}(G\times G)$. If
$F\in\mathcal{D}(\iota\otimes\sigma_{-\frac{i}{2}})$ then
$$\Delta(x)F^{*}(y\otimes
1)\in\mathcal{N}_{\varphi_{\Omega}\otimes\varphi_{\Omega}}\quad\text{and}$$
$$\left(\Lambda_{\Omega}\otimes\Lambda_{\Omega}\right)\left(\Delta(x)F^{*}(y\otimes
1)\right)= (\hat{J}\otimes
J)W(\iota\otimes\sigma_{-\frac{i}{2}})(\tilde{F})(\hat{J}\otimes
J)\left(A^{-\frac{1}{2}}
\Lambda_{\Omega}(y)\otimes\Lambda_{\Omega}(x)\right).$$
\end{lemma}
\begin{proof}
According to Proposition \ref{propVaes} and Lemma \ref{Twisting1}, it
suffices to show that
\begin{eqnarray}\label{Eqtwist1}\notag
&\forall F\in\mathcal{D}(\iota\otimes\sigma_{-\frac{i}{2}}),\,\,
\Delta(\overline{xA^{\frac{1}{2}}})F^{*}(\overline{A^{-\frac{1}{2}}yA^{\frac{1}{2}}}\otimes
1)\in\mathcal{N}_{\varphi\otimes\varphi}\quad\text{and,}&\\
&(\Lambda\otimes\Lambda)\left(\Delta(\overline{xA^{\frac{1}{2}}})F^{*}
(\overline{A^{-\frac{1}{2}}yA^{\frac{1}{2}}}\otimes
1)\right)&\\\notag &=(\hat{J}\otimes
J)W(\iota\otimes\sigma_{-\frac{i}{2}})(\tilde{F})(\hat{J}\otimes
J)\left(A^{-\frac{1}{2}}\Lambda_{\Omega}(y)\otimes\Lambda_{\Omega}(x)\right).&\\\notag
\end{eqnarray}
Let $F\in L^{\infty}(G\times G)$. We identify $\sigma$ with its
restriction to $L^{\infty}(G)$. A direct application of Lemma
\ref{group} $(2)$ gives that
$L^{\infty}(G)\odot\mathcal{D}(\sigma_{-\frac{i}{2}})$ is a
$\sigma$-strong* core for $\iota\otimes\sigma_{-\frac{i}{2}}$.
Taking into account the observation preceding this lemma and
because $\Lambda\otimes\Lambda$ is $\sigma$-strong*-norm closed,
it suffices to show (\ref{Eqtwist1}) for $F\in
L^{\infty}(G)\odot\mathcal{D}(\sigma_{-\frac{i}{2}})$. By
linearity, we only have to show (\ref{Eqtwist1}) for $F$ of the
form $F=F_{1}\otimes F_{2}$ with $F_{1}, F_{2}\in L^{\infty}(G)$
and $F_{2}\in\mathcal{D}(\sigma_{-\frac{i}{2}})$. Proposition \ref{PropTx} gives
$\overline{A^{-\frac{1}{2}}yA^{\frac{1}{2}}}\in\mathcal{N}_{\varphi}$, so
$\Delta(\overline{xA^{\frac{1}{2}}})(F_{1}^{*}
\overline{A^{-\frac{1}{2}}yA^{\frac{1}{2}}}\otimes
1)\in\mathcal{N}_{\varphi\otimes\varphi}$, and writing
$$\Delta(\overline{xA^{\frac{1}{2}}})(F_{1}^{*}\otimes
F_{2}^{*})(\overline{A^{-\frac{1}{2}}yA^{\frac{1}{2}}}\otimes 1)
=\Delta(\overline{xA^{\frac{1}{2}}})(F_{1}^{*}\overline
{A^{-\frac{1}{2}}yA^{\frac{1}{2}}}\otimes
1)(1\otimes  F_{2}^{*})$$
with $1\otimes F_{2}\in\mathcal{D}(\iota\otimes\sigma_{-\frac{i}{2}})$, we see
that
$\Delta(\overline{xA^{\frac{1}{2}}})F^{*}(\overline{A^{-\frac{1}{2}}yA^{\frac{1}{2}}}\otimes
1)\in\mathcal{N}_{\varphi\otimes\varphi}$ and
\begin{eqnarray*}
&&(\Lambda\otimes\Lambda)\left(\Delta(\overline{xA^{\frac{1}{2}}})
F^{*}(\overline{A^{-\frac{1}{2}}yA^{\frac{1}{2}}}\otimes
1)\right)\\
&&=\left(1\otimes
J\sigma_{-\frac{i}{2}}(F_{2})J\right)(\Lambda\otimes\Lambda)
\left(\Delta(\overline{xA^{\frac{1}{2}}})(F_{1}^{*}
\overline{A^{-\frac{1}{2}}yA^{\frac{1}{2}}}\otimes
1)\right)\\
&&=\left(1\otimes
J\sigma_{-\frac{i}{2}}(F_{2})J\right)W^{*}\Lambda(F_{(1)}^{*}
\overline{A^{-\frac{1}{2}}yA^{\frac{1}{2}}})\otimes\Lambda(\overline{xA^{\frac{1}{2}}})\\
&&\quad(\text{by definition of }W)\\
&&=\left(1\otimes
J\sigma_{-\frac{i}{2}}(F_{2})J\right)(\hat{J}\otimes
J)W(\hat{J}\otimes J)(F_{1}^{*}\otimes
1)\Lambda(\overline{A^{-\frac{1}{2}}yA^{\frac{1}{2}}})\otimes\Lambda(\overline{xA^{\frac{1}{2}}})\\
&&\quad(\text{because}\quad W^{*}=(\hat{J}\otimes
J)W(\hat{J}\otimes
J))\\
&&=(\hat{J}\otimes
J)\left(1\otimes\sigma_{-\frac{i}{2}}(F_{2})\right)W\left(R(F_{1})\otimes1\right)(\hat{J}\otimes
J)A^{-\frac{1}{2}}\Lambda_{\Omega}(y)\otimes\Lambda_{\Omega}(x)\\
&&\quad(\text{because}\quad
R(x)=\hat{J}x^{*}\hat{J},\quad\text{and}\quad
\Lambda(\overline{A^{-\frac{1}{2}}yA^{\frac{1}{2}}})=A^{-\frac{1}{2}}\Lambda_{\Omega}(y)\,\,\,\text{by
Proposition \ref{PropTx}})\\
&&=(\hat{J}\otimes
J)W\Delta\left(\sigma_{-\frac{i}{2}}(F_{2})\right)\left(R(F_{1})\otimes1\right)(\hat{J}\otimes
J)A^{-\frac{1}{2}}\Lambda_{\Omega}(y)\otimes\Lambda_{\Omega}(x)\\
&&\quad(\text{because}\quad\Delta(x)=W^{*}(1\otimes x)W).\\
\end{eqnarray*}
So we just have to compute:
$$
\Delta\left(\sigma_{-\frac{i}{2}}(F_{2})\right)\left(R(F_{1})\otimes1\right)(g,h)
=F_{1}(g^{-1})\sigma_{-\frac{i}{2}}(F_{2})(gh)$$
$$=(\iota\otimes\sigma_{-\frac{i}{2}})(F)(g^{-1},gh)=
\widetilde{(\iota\otimes\sigma_{-\frac{i}{2}})(F)}(g,h).$$

\end{proof}
The next lemma is necessary to finish the proof of the theorem.
\begin{lemma}\label{Domain}

\item (1) We have $\hat{J}A^{-\frac{1}{2}}=A^{\frac{1}{2}}\hat{J}$.
\item (2) The operator $(\iota\otimes\sigma_{-\frac{i}{2}})(\tilde{\Omega})$ is normal, affiliated with $M\otimes M$, and its polar decomposition is
$$(\iota\otimes\sigma_{-\frac{i}{2}})(\tilde{\Omega})=\tilde{\Omega}(A^{-\frac{1}{2}}\otimes 1).$$

\end{lemma}
\begin{proof}
Let $\alpha\in\R$ and
$\mathcal{D}_{\alpha}$ the horizontal strip bounded by $\R$ and $\R-i\alpha$.

$(1)$ Let $\xi\in\mathcal{D}(A^{-\frac{1}{2}})$. There exists a continuous bounded extension $F$ of $t\mapsto A^{it}\xi$ on $\mathcal{D}_{-\frac{1}{2}}$ which is analytic on
$\mathcal{D}_{-\frac{1}{2}}^{0}$. The function $G(z)=\hat{J}
F(\overline{z})$ is continuous bounded on $\mathcal{D}_{\frac{1}{2}}$ and analytic on $\mathcal{D}_{\frac{1}{2}}^{0}$. Moreover :
$$R(A^{-it})(g)=\Psi(g^{-1},\gamma_{t})=\Psi(g,\gamma_{t}^{-1})=A^{it}(g),\quad\text{for all}\,\,g\in G,t\in\R.$$
Thus, $\hat{J}A^{it}\hat{J}=R(A^{-it})=A^{it}$. We deduce $G(t)=\hat{J}A^{it}\xi=A^{it}\hat{J}\xi$.
This means that $\hat{J}\xi\in\mathcal{D}(A^{\frac{1}{2}})$ and
$A^{\frac{1}{2}}\hat{J}\xi=G(-\frac{i}{2})=\hat{J}F(\frac{i}{2})
=\hat{J}A^{-\frac{1}{2}}\xi$, so
$\hat{J}A^{-\frac{1}{2}}\subset A^{\frac{1}{2}}\hat{J}$. The other inclusion can be
proved in the same way.

$(2)$ Note that
$$
(\iota\otimes\sigma_{t})(\tilde{\Omega})(g,h)=\Psi(g^{-1},gh\gamma_{t}^{-1}) = \Psi(g^{-1},gh)\Psi(g,\gamma_{t})
=\tilde{\Omega}(A^{-it}\otimes 1)(g,h).$$
We conclude the proof applying Proposition \ref{Propext}.
\end{proof}

We can now prove the theorem. Let $x\in N$ and $y\in L$. Put
$\xi=\hat{J}\Lambda_{\Omega}(y)\in\mathcal{D}(A^{\frac{1}{2}})$
and $\eta=J\Lambda_{\Omega}(x)$. By Lemma \ref{Domain} $(2)$,
$A^{\frac{1}{2}}\xi\otimes\eta\in\mathcal{D}\left
((\iota\otimes\sigma_{-\frac{i}{2}})(\tilde{\Omega})\right)$.
Thus, using Lemma \ref{group} $(1)$, there exists
$\tilde{\Omega}_{n}\in L^{\infty}(G\times
G)\cap\mathcal{D}(\iota\otimes\sigma_{-\frac{i}{2}})$ such that
$$
\tilde{\Omega}_{n}\rightarrow\tilde{\Omega}\quad\sigma\text{-strongly}*\quad\text{and}
\quad (\iota\otimes\sigma_{-\frac{i}{2}})(\tilde{\Omega}_{n})
(A^{\frac{1}{2}}\xi\otimes\eta)\rightarrow(\iota\otimes
\sigma_{-\frac{i}{2}})(\tilde{\Omega})(A^{\frac{1}{2}}\xi\otimes\eta).$$
Because $\tilde{\tilde{F}}=F$, we also have
$\Omega_{n}\rightarrow\Omega$ $\sigma$-strongly*, so
$$\Delta(x)\Omega_{n}^{*}(y\otimes
1)\rightarrow\Delta(x)\Omega^{*}y\otimes
1\quad\sigma\text{-strongly}^*.$$ By Lemma \ref{element},
$\Delta(x)\Omega_{n}^{*}(y\otimes
1)\in\mathcal{N}_{\varphi_{\Omega}\otimes\varphi_{\Omega}}$ and
\begin{eqnarray*}
\left(\Lambda_{\Omega}\otimes\Lambda_{\Omega}\right)\left(\Delta(x)\Omega_{n}^{*}(y\otimes
1)\right) &=&(\hat{J}\otimes
J)W(\iota\otimes\sigma_{-\frac{i}{2}})(\tilde{\Omega}_{n})(\hat{J}\otimes
J)\left(A^{-\frac{1}{2}}\Lambda_{\Omega}(y)\otimes\Lambda_{\Omega}(x)\right)\\
&=&(\hat{J}\otimes J)W(\iota\otimes\sigma_{-\frac{i}{2}})
(\tilde{\Omega}_{n})(A^{\frac{1}{2}}\xi\otimes\eta)\quad(\text{by
Lemma \ref{Domain}(1)})\\
&\rightarrow&(\hat{J}\otimes
J)W(\iota\otimes\sigma_{-\frac{i}{2}})(\tilde{\Omega})(A^{\frac{1}{2}}\xi\otimes\eta)\\
&=&(\hat{J}\otimes J)W\tilde{\Omega}(\xi\otimes\eta)\quad(\text{by
Lemma
\ref{Domain}(2)})\\
&=&(\hat{J}\otimes J)W\tilde{\Omega}(\hat{J}\otimes
J)(\Lambda_{\Omega}(y)\otimes\Lambda_{\Omega}(x)).\\
\end{eqnarray*}
Because $\Lambda_{\Omega}\otimes\Lambda_{\Omega}$ is
$\sigma$-strong* - norm closed, we have
$\Delta(x)\Omega^{*}(y\otimes
1)\in\mathcal{N}_{\varphi_{\Omega}\otimes\varphi_{\Omega}}$, so
$\Delta_{\Omega}(x)(y\otimes
1)\in\mathcal{N}_{\varphi_{\Omega}\otimes\varphi_{\Omega}}$ and
\begin{eqnarray*}
\left(\Lambda_{\Omega}\otimes\Lambda_{\Omega}\right)\left(\Delta_{\Omega}(x)(y\otimes
1)\right)
&=&\Omega\left(\Lambda_{\Omega}\otimes\Lambda_{\Omega}\right)\left(\Delta(x)\Omega^{*}(y\otimes
1)\right)\\
&=&\Omega(\hat{J}\otimes J)W\tilde{\Omega}(\hat{J}\otimes
J)(\Lambda_{\Omega}(y)\otimes\Lambda_{\Omega}(x))\\
&=&W_{\Omega}^{*}(\Lambda_{\Omega}(y)\otimes\Lambda_{\Omega}(x)).\\
\end{eqnarray*}
\end{proof}

Let $R_{\Omega}=uR(x)u^{*}$ be the *-anti-automorphism of $M$, where $u=\alpha(\Psi(\cdot^{-1},\cdot))$.

\textit{Proof of Theorem\ref{ThmTwist}}. Let $x,y\in\mathcal{N}_{\varphi_{\Omega}}$.
By Theorem\ref{TwThm1}, we have
\begin{eqnarray}\notag
&||\left(\Lambda_{\Omega}\otimes\Lambda_{\Omega}\right)\left(\Delta_{\Omega}(x)y\otimes
1\right)||^{2}
=||\Lambda_{\Omega}(y)\otimes\Lambda_{\Omega}(x)||^{2}&\\\label{EqTh1}
&\Leftrightarrow\left(\omega_{\Lambda_{\Omega}(y)}\otimes\varphi_{\Omega}\right)
(\Delta_{\Omega}(x^{*}x))=\omega_{\Lambda_{\Omega}(y)}(1)\varphi_{\Omega}(x^{*}x).&\\\notag
\end{eqnarray}
Let $\omega\in M_{*}$, $\omega\geq 0$. The inclusion
$M\subset\mathcal{B}(H)$ is standard, so there is $\xi\in H$ such
that $\omega=\omega_{\xi}$. Let $a_{i}\in M$ such that
$\Lambda_{\Omega}(a_{i})\rightarrow \xi$. Then
\begin{equation}\label{EqTh2}
\omega_{\Lambda_{\Omega}(a_{i})}(x)\rightarrow\omega(x),\quad\text{
for all}\quad x\in M.
\end{equation}
To show that $\varphi_{\Omega}$ is left invariant, it suffices to
show that
$\Delta_{\Omega}(x^{*}x)\in\mathcal{N}_{\iota\otimes\varphi_{\Omega}}$
when $x\in\mathcal{N}_{\varphi_{\Omega}}$. Indeed, in this case we
have, using (\ref{EqTh2}),
\begin{eqnarray*}
&\omega_{\Lambda_{\Omega}(a_{i})}(1)\varphi_{\Omega}(x^{*}x)
\rightarrow\omega(1)\varphi_{\Omega}(x^{*}x)\quad\text{and,
}&\\
&\left(\omega_{\Lambda_{\Omega}(a_{i})}\otimes
\varphi_{\Omega}\right)(\Delta_{\Omega}(x^{*}x))
\rightarrow\left(\omega\otimes\varphi_{\Omega}\right)(\Delta_{\Omega}(x^{*}x)).&\\
\end{eqnarray*}
This implies, using (\ref{EqTh1}), that for all $\omega\in
M_{*}^{+}$ and $x\in\mathcal{N}_{\varphi_{\Omega}}$,
$$
\left(\omega\otimes\varphi_{\Omega}\right)
(\Delta_{\Omega}(x^{*}x))=\omega(1)\varphi_{\Omega}(x^{*}x),
$$
\ie{}, $\varphi_{\Omega}$ is left invariant. Let us show that
$\Delta_{\Omega}(x^{*}x)\in\mathcal{N}_{\iota\otimes\varphi_{\Omega}}$.
We put
$$m=(\iota\otimes\varphi_{\Omega})(\Delta_{\Omega}(x^{*}x))\in
M^{\text{ext}}_{+}.$$ The spectral decomposition of $m$ is $m=\int_{0}^{\infty}\lambda d e_{\lambda} + \infty \, .\,p$. From (\ref{EqTh1}) we see that, for all
$y\in\mathcal{N}_{\varphi_{\Omega}}$,
$m(\omega_{\Lambda_{\Omega}(y)})<\infty$. Thus, the set
$\{\omega\in M_{*}^{+}\,\,|\,\,m(\omega)<\infty\}$ is dense in
$M_{*}^{+}$. This implies $p=0$ and $m=m_{T}$, where $T$ is the
positive operator affiliated with $M$ defined by
$$T=\int_{0}^{\infty}\lambda d e_{\lambda}.$$
So, we only have to show that $T$ is a bounded operator. Using
again (\ref{EqTh1}) and the definition of $m_{T}$, we see that,
for all $y\in\mathcal{N}_{\varphi_{\Omega}}$,
$\Lambda_{\Omega}(y)\in\mathcal{D}(A^{\frac{1}{2}})$ and
$$||T^{\frac{1}{2}}\Lambda_{\Omega}(y)||^{2}=
\varphi_{\Omega}(x^{*}x)||\Lambda_{\Omega}(y)||^{2}.$$
Thus, $T$ is a bounded operator.

It is easy to check (see \cite{Vain}) that $\Delta_{\Omega}\circ
R_{\Omega}=\sigma (R_{\Omega}\otimes R_{\Omega})\Delta_{\Omega},$
so the right invariance of $\varphi_{\Omega}\circ R_{\Omega}$ follows. Thus, $(M,\Delta_{\Omega})$ is a \lc{} quantum group and it
follows immediately from Theorem \ref{TwThm1} that $W_{\Omega}$ is
its multiplicative unitary. Our next aim is to show that $\psi_{\Omega}=
\varphi_{\Omega}\circ R_{\Omega}$. We compute:
\begin{eqnarray*}
R\left(u\sigma_{-t}^{\Omega}(u^{*})\right)(g)
&=& u(g^{-1})u^{*}(g^{-1}\gamma_{t})=u(g)\Psi(g^{-1}\gamma_{t},g^{-1}\gamma_{t})\\
&=&u(g)u^{*}(g)\Psi(g^{-1},\gamma_{t})\Psi(\gamma_{t},g^{-1})\Psi(\gamma_{t},\gamma_{t})\\
&=&\lambda^{it^{2}}(A^{it}B^{it})(g).\\
\end{eqnarray*}
This implies
\begin{eqnarray*}
[D\varphi_{\Omega}\circ R_{\Omega}\,:\,D\psi]_{t}
&=&[D\left(\varphi_{\Omega}\right)_{u}\circ R\,:\,D\varphi\circ
R]_{t}=R\left([D\left(\varphi_{\Omega}\right)_{u}\,:\,D\varphi]_{-t}^{*}\right)\\
&=&R\left([D\left(\varphi_{\Omega}\right)_{u}\,:\,D\varphi_{\Omega}]_{-t}^{*}
\right)R\left([D\varphi_{\Omega}\,:\,D\varphi]_{-t}^{*}\right)\\
&=&R\left(u\sigma_{-t}^{\Omega}(u^{*})\right)R\left(\lambda^{-i\frac{t^{2}}{2}}A^{it}\right)\\
&=&\lambda^{it^{2}}A^{it}B^{it}(\lambda^{-i\frac{t^{2}}{2}}A^{-it})\\
&=&\lambda^{i\frac{t^{2}}{2}}B^{it}.\\
\end{eqnarray*}
Thus, $\psi_{\Omega}=\varphi_{\Omega}\circ R_{\Omega}$. In order to finish the proof,
we have to compute the scaling group and the scaling constant. Recall that if
$(M,\Delta)$ is a \lc{} quantum group, then the scaling group
is the unique one-parameter group $\tau_{t}$ on $M$ such that
$\Delta\circ\sigma_{t}=(\tau_{t}\otimes\sigma_{t})\circ\Delta.$
Since $(\iota\otimes\sigma_{t})(\Omega)=\Omega (A^{it}\otimes 1)$,
using $\tau_{t}\circ\alpha=\alpha$, we have
$(\tau_{t}\otimes\sigma_{t})(\Omega)=\Omega (A^{it}\otimes 1)$, which gives:
\begin{eqnarray*}
(\tau_{t}\otimes\sigma_{t}^{\Omega})(\Delta_{\Omega}(x))
&=&(1\otimes
A^{it})(\tau_{t}\otimes\sigma_{t})(\Omega)(\tau_{t}
\otimes\sigma_{t})(\Delta(x))(\tau_{t}\otimes\sigma_{t})(\Omega^{*})(1\otimes
A^{-it})\\
&=&\Omega (A^{it}\otimes
A^{it})(\tau_{t}\otimes\sigma_{t})(\Delta(x))(A^{-it}\otimes A^{-it})\Omega^{*}\\
&=&\Omega\Delta(A^{it})\Delta(\sigma_{t}(x))\Delta(A^{-it})\Omega^{*}\\
&=&\Delta_{\Omega}(\sigma_{t}^{\Omega}(x)).\\
\end{eqnarray*}
This relation characterizes the scaling group of
$(M,\Delta_{\Omega})$. Recall that the scaling constant of $(M,\Delta)$ verifies $\varphi\circ\tau_{t}=\nu^{-t}\varphi$. Because $\tau_{t}(A^{is})=A^{is}$, for all
$t,s\in\R$, we deduce that $\varphi_{\Omega}\circ\tau_{t}^{\Omega}=\varphi_{\Omega}\circ\tau_{t}=\nu^{-t}\varphi_{\Omega}$. Thus, $\nu^{\Omega}=\nu$.\quad\quad

Let us denote by $X$ and $Y$ the operators
$$X=\Omega^{*}\quad\text{and}\quad Y=(\hat{J}\otimes J)(u^{*}\otimes
1)\Omega(\hat{J}\otimes J).$$ Note that
$\tilde{\Psi}^{*}(g,h)=\Psi^{*}(g^{-1},g)\Psi(g,h)$, so
$\tilde{\Omega}^{*}=(u^{*}\otimes 1)\Omega$ and
$$
W_{\Omega} =(\hat{J}\otimes
J)\tilde{\Omega}^{*}W^{*}(\hat{J}\otimes J){\Omega}^{*}
=(\hat{J}\otimes J)\tilde{\Omega}^{*}(\hat{J}\otimes
J)W{\Omega}^{*} =YWX.$$ From now on we suppose that $G$ is
abelian, we switch to the additive notations for its operations
and denote by $\hat G$ its dual. Recall that the notations $u_{\gamma}$, $L_{\gamma}$ and $R_{\gamma}$ where introduced in Section \ref{SectionResult}. Note that $R(L_{\gamma})=L^{*}_{\gamma}=L_{-\gamma}$.

\begin{proposition}\label{TwistedIng}
$R_{\Omega}$ is the unitary antipode of $(M,\Delta_{\Omega})$.
Moreover,
\begin{itemize}
\item $\delta_{\Omega}=\delta A^{-1}B$,
\item $\mathcal{D}(S_{\Omega})=\mathcal{D}(S)$ and, for all
$x\in\mathcal{D}(S),\ S_{\Omega}(x)=uS(x)u^{*}.$
\end{itemize}
\end{proposition}
\begin{proof}
If $(M,\Delta)$ is a \lc{} quantum group, then the unitary
antipode is the unique *-anti-automorphism $R$ of $M$ such that
$R\left((\iota\otimes\omega_{\xi,\eta})(W)\right)=(\iota\otimes\omega_{J\eta,J\xi})(W)$.
Let us define two *-homomorphisms by
\begin{eqnarray*}
&\pi^{'}\,\,:\,\,L^{\infty}(G\times G)\rightarrow M\otimes M^{'}
:\ \pi^{'}(F)=(\hat{J}\otimes
J)(\alpha\otimes\alpha)(F)^{*}(\hat{J}\otimes J),\\
&\pi\,\,:\,\,L^{\infty}(G\times G)\rightarrow M\otimes M :\
\pi(F)=(\alpha\otimes\alpha)(F).\\
\end{eqnarray*}
We want to prove that, for all $F,G\in L^{\infty}(G\times G)$ and
$\xi,\eta\in H$,
\begin{equation}\label{EqR1}
R\left((\iota\otimes\omega_{\xi,\eta})\left(\pi^{'}(F)W\pi(G)\right)\right)
=(\iota\otimes\omega_{J\eta,J\xi})\left(\pi^{'}(G)W\pi(F)\right).
\end{equation}
By linearity and continuity, it suffices to prove $(\ref{EqR1})$
for $F=u_{\gamma_{1}}\otimes u_{\gamma_{2}}$ and
$G=u_{\gamma_{2}}\otimes u_{\gamma_{4}}$ with
$\gamma_{i}\in\hat{G}$. We have
$$\pi^{'}(u_{\gamma_{1}}\otimes u_{\gamma_{2}})=L_{-\gamma_{1}}\otimes
R_{-\gamma_{2}}\quad\text{and}\quad \pi(u_{\gamma_{3}}\otimes
u_{\gamma_{4}})=L_{\gamma_{3}}\otimes L_{\gamma_{4}},\quad\text{so}$$
\begin{eqnarray*}
&&R\left((\iota\otimes\omega_{\xi,\eta})\left(\pi^{'}(u_{\gamma_{1}}\otimes
u_{\gamma_{2}})W\pi(u_{\gamma_{3}}\otimes
u_{\gamma_{4}})\right)\right)
\\
=&&R\left((\iota\otimes\omega_{\xi,\eta})\left(L_{-\gamma_{1}}\otimes
R_{-\gamma_{2}}WL_{\gamma_{3}}\otimes L_{\gamma_{4}}\right)\right)\\
=&&R\left(L_{-\gamma_{1}}(\iota\otimes
L_{\gamma_{4}}.\omega_{\xi,\eta}.R_{-\gamma_{2}})
\left(W\right)L_{\gamma_{3}}\right)\\
=&&L_{-\gamma_{3}}R\left((\iota\otimes\omega_{L_{\gamma_{4}}
\xi,R_{\gamma_{2}}\eta})\left(W\right)\right)L_{\gamma_{1}}\\
=&&L_{-\gamma_{3}}(\iota\otimes\omega_{JR_{\gamma_{2}}
\eta,JL_{\gamma_{4}}\xi})\left(W\right)L_{\gamma_{1}}\\
=&&L_{-\gamma_{3}}(\iota\otimes\omega_{L_{\gamma_{2}}J
\eta,R_{\gamma_{4}}J\xi})\left(W\right)L_{\gamma_{1}}\\
=&&(\iota\otimes
L_{\gamma_{2}}.\omega_{J\eta,J\xi}.R_{-\gamma_{4}})
\left(L_{-\gamma_{3}}\otimes 1 W L_{\gamma_{1}}\otimes 1\right)\\
=&&(\iota\otimes\omega_{J\eta,J\xi})\left(L_{-\gamma_{3}}\otimes
R_{-\gamma_{4}}W L_{\gamma_{1}}\otimes L_{\gamma_{2}}\right)
\\
=&&(\iota\otimes\omega_{J\eta,J\xi})\left(\pi^{'}(u_{\gamma_{3}}\otimes
u_{\gamma_{4}})W\pi(u_{\gamma_{1}}\otimes u_{\gamma_{2}})\right).\\
\end{eqnarray*}
Note that $Y=\pi^{'}(\tilde{\Psi})$, $X=\pi(\Psi^{*})$ and
$\pi(\tilde{\Psi})(u^{*}\otimes 1)=\tilde{\Omega}(u^{*}\otimes
1)=\Omega^{*}=X.$ Also, using $R(u^{*})=u^{*}$, we have
\begin{eqnarray*}
(u\otimes 1)\pi^{'}(\Psi^{*})
&=&(u\otimes 1)(\hat{J}\otimes J)\Omega(\hat{J}\otimes J)\\
&=&(\hat{J}\otimes J)(R(u^{*})\otimes 1)\Omega(\hat{J}\otimes J)\\
&=&(\hat{J}\otimes J)(u^{*}\otimes 1)\Omega(\hat{J}\otimes J)\\
&=&(\hat{J}\otimes J)\tilde{\Omega}^{*}(\hat{J}\otimes J)=Y.\\
\end{eqnarray*}
Using these remarks and relation $(\ref{EqR1})$, one has
\begin{eqnarray*}
R_{\Omega}\left((\iota\otimes\omega_{\xi,\eta})\left(W_{\Omega}\right)\right)
&=&uR\left((\iota\otimes\omega_{\xi,\eta})\left(\pi^{'}
(\tilde{\Psi})W\pi(\Psi^{*})\right)\right)u^{*}\\
&=&(\iota\otimes\omega_{J\eta,J\xi})\left((u\otimes
1)\pi^{'}(\Psi^{*})W\pi(\tilde{\Psi})(u^{*}\otimes 1)\right)\\
&=&(\iota\otimes\omega_{J\eta,J\xi})\left(YWX\right)\\
&=&(\iota\otimes\omega_{J_{\Omega}\eta,J_{\Omega}\xi})\left(W_{\Omega}\right).\\
\end{eqnarray*}
Where we use, in the last equality, the fact that
$J_{\Omega}=\lambda^{\frac{i}{8}}J$ so
$\omega_{J_{\Omega}\eta,J_{\Omega}\xi}=\omega_{J\eta,J\xi}$. This
relation characterizes the unitary antipode of
$(M,\Delta_{\Omega})$.
We have
\begin{eqnarray*}
[D\psi_{\Omega}\,:\,D\varphi_{\Omega}]_{t}
&=&[D\psi_{\Omega}\,:\,D\psi]_{t}[D\psi\,:\,D\varphi]_{t}[D\varphi\,:\,D\varphi_{\Omega}]_{t}\\
&=&(\lambda^{i\frac{t^{2}}{2}}B^{it})(\nu^{i\frac{t^{2}}{2}}\delta^{it})(\lambda^{-i\frac{t^{2}}{2}}A^{-it})\\
&=&\nu^{i\frac{t^{2}}{2}}(\delta A^{-1}B)^{it}.\\
\end{eqnarray*}
Thus, $\delta_{\Omega}=\delta A^{-1}B$ (because we have seen in the proof of Theorem \ref{ThmTwist} that $\psi_{\Omega}=\varphi_{\Omega}\circ R_{\Omega}$). The
last statement is clear.\quad\quad
\end{proof}
\begin{remark}
If $(M,\Delta)$ is a Kac algebra, $(M,\Delta_{\Omega})$ is
not in general a Kac algebra (see Section \ref{Sectionaz+bclassique}). However, in the case
considered in \cite{Enock}, \cite{Vain}, and \cite{Kas}, when
$\alpha(L^{\infty}(G))$ belongs to the fixed point subalgebra of
$M$ with respect to $\sigma_t$, then $\gamma_{t}$ is trivial and
we have $A^{-1}B=1$, so $(M,\Delta_{\Omega})$ is a Kac algebra .
\end{remark}

\begin{remark}
The map $L^{\infty}(G\times G)\rightarrow\mathcal{B}(H\otimes H)\,\,:\,\,F\mapsto \pi^{'}(\tilde{F})W\pi(F^{*})$ is $\sigma$-strong*-$\sigma$-weak continuous. So, if $(x\mapsto\Psi_{x})$ is $\sigma$-strongly* continuous map from $\R$ to
$L^{\infty}(G\times G)$ such that $\Psi_{x}$ is a continuous bicharacter, for all $x\in\R$, then, denoting by $W_{x}$ the multiplicative unitary of the twisted \lc{} quantum group associated with $\Psi_{x}$, the map $x\mapsto W_{x}$ from $\R$ to the unitaries of $\mathcal{B}(H\otimes H)$ is $\sigma$-weakly continuous. This is the case for the example
of section \ref{Sectionaz+bclassique} and for the examples constructed in \cite{FimaVain}
and \cite{Kas}.
\end{remark}

\section{Rieffel's deformations of locally compact quantum group}\label{SectionDual}
This section is devoted to the proof of Theorem \ref{TheoremRieffel}. We use the hypotheses
and notations from the previous section and from Section \ref{SectionResult}. So let $\widehat{G}<(M,\Delta)$ be a stable co-subgroup with $G$ abelian. Recall that we have (see Section \ref{SectionCo-Subgroups}) two unitary representations of $\hat{G}$ : $\gamma\mapsto L_{\gamma}$ and $\gamma\mapsto R_{\gamma}$ of $\widehat{G}$. This gives two *-homomorphisms from $L^{\infty}(G)$ to $\mathcal{B}(H)$, $\pi_{L}$ and $\pi_{R}$, respectively. We have
$$\pi_{L}=\alpha\quad\text{and}\quad\pi_{R}(F)=
J\alpha(\overline{F(-\cdot)})J=J\widehat{J}\alpha(F)\widehat{J}J.$$
Recall that $W_{\mathrm{\Omega}}=YWX$, where
$X=(\alpha\otimes\alpha)(\Psi^{*})=(\pi_{L}\otimes\pi_{L})(\Psi^{*})$
and $Y=(\widehat{J}\otimes J)\tilde{\Omega}^{*}(\widehat{J}\otimes J)=(\pi_{L}\otimes\pi_{R})(\tilde{\Psi})=(\alpha(u)\otimes 1)(\pi_{L}\otimes\pi_{R})(\Psi^{*})$.
Note that $\widehat{G}<(M_{\Omega},\Delta_{\Omega})$ is also stable (by the preceding
section), so the results of section \ref{SectionCo-Subgroups} can be applied also to $\widehat{G}<(M_{\Omega},\Delta_{\Omega})$. Thus, we have a left-right action $\alpha$ of $\widehat{G}^{2}$ on $\widehat{M}$ and also a left-right action $\beta$ of $\widehat{G}^{2}$ on $\widehat{M}_{\Omega}$. We denote by the same $\pi$ the canonical morphism from $\widehat{M}$ in the crossed product $N=\widehat{G}^{2}\ltimes\widehat{M}$ and from $\widehat{M}_{\mathrm{\Omega}}$ in  $\widehat{G}^{2}\ltimes\widehat{M}_{\mathrm{\Omega}}$. Also we denote by $\lambda_{\gamma_{1},\gamma_{2}}$ the canonical unitaries in the two crossed products and by the same $\theta$ the dual action on $\widehat{G}^{2}\ltimes\widehat{M}$ and $\widehat{G}^{2}\ltimes\widehat{M}_{\mathrm{\Omega}}$. Recall that $\theta$ and $\lambda$ verify
$$
\theta_{g_{1},g_{2}}(\lambda_{\gamma_{1},\gamma_{2}})=\overline{<\gamma_{1},g_{1}>} \overline{<\gamma_{2},g_{2}>}\lambda_{\gamma_{1},\gamma_{2}}.
$$
The unitary representations $\gamma\mapsto\lambda_{(\gamma,0)}$, $\gamma\mapsto\lambda_{(0,\gamma)}$ and $\lambda$ give unital normal *-homomorphisms $\lambda_{L},\lambda_{R}:L^{\infty}(G)\to \widehat{G}^{2}\ltimes\widehat{M}$ and $\lambda: L^{\infty}(G^{2})\to \widehat{G}^{2}\ltimes\widehat{M}$ verifying
$$\lambda_{L}(u_{\gamma})=\lambda_{(\gamma,0)},\quad\lambda_{R}(u_{\gamma})=
\lambda_{(0,\gamma)},\quad
\lambda(u_{(\gamma_{1},\gamma_{2})})=\lambda_{\gamma_{1},\gamma_{2}}.$$
Since $\lambda(f_{1}\otimes
f_{2})=\lambda_{L}(f_{1})\lambda_{R}(f_{2})$, then
\begin{eqnarray}\label{EqDual2}
\theta_{(g_{1},g_{2})}(\lambda_{L}(F))=\lambda_{L}(F(\cdot-g_{1})),
\quad\text{for any}\,\,F\in L^{\infty}(G).\\ \notag
\theta_{(g_{1},g_{2})}(\lambda_{R}(F))=\lambda_{R}(F(\cdot-g_{2})),
\quad\text{for any}\,\,F\in L^{\infty}(G).
\end{eqnarray}
We have for the twisted dual action $\theta^{\Psi}$:
\begin{equation}\label{DualActHatM}
\theta^{\Psi}_{(g_{1},g_{2})}(\pi(x))=\pi(\alpha_{(\Psi_{-g_{1}},\Psi_{g_{2}})}(x)),
\quad\text{for all}\,\, x\in\widehat{M}.
\end{equation}
Considering the following unitaries in $M\otimes N$:
$$
\tilde{X}=(\alpha\otimes\lambda_{L})(\Psi^{*})\,\, ,\
\tilde{Y}=(\alpha\otimes\lambda_{R})(\tilde{\Psi})=(\alpha(u)\otimes 1)(\alpha\otimes\lambda_{R})(\Psi^{*})
,\ \tilde{W}=(\iota\otimes\pi)(W),
$$
we put $\tilde{W}_{\mathrm{\Omega}}=\tilde{Y}\tilde{W}\tilde{X}$. Let $N_{\mathrm{\Omega}}$ be the fixed point subalgebra of $\widehat{G}^{2}\ltimes\widehat{M}$ under the twisted dual action. The step to prove Theorem \ref{TheoremRieffel} is to show that $\widehat{M}_{\mathrm{\Omega}}$ is isomorphic to $N_{\mathrm{\Omega}}$, for this we need a preliminary lemma. Let $\mathcal{B}$ be the von Neumann algebra acting on $H$ and generated
by $\{(\omega\otimes\iota)(W\Omega^{*})\,|\,\omega\in\bh_{*}\}$.

\begin{lemma}\label{TwistedGVN}
We have:
\begin{itemize}
\item $\mathcal{B}\vee\{L_{\gamma}\,|\,\gamma\in\widehat{G}\}^{''}=
    \widehat{M}\vee\{L_{\gamma}\,|\,\gamma\in\widehat{G}\}^{''}$,
\item $\mathcal{B}\vee\{R_{\gamma}\,|\,\gamma\in\widehat{G}\}^{''}=
    \widehat{M}_{\mathrm{\Omega}}\vee\{R_{\gamma}\,|\,\gamma\in\widehat{G}\}^{''}$.
\end{itemize}
\end{lemma}

\begin{proof}
First, take a net in the vector space spanned by elements $u_{\gamma_{1}}\otimes u_{\gamma_{2}}$ such that $\sum c_{i,j} u_{\gamma_{i},\gamma_{j}}\rightarrow \Psi\quad\text{strongly*}$.
Then $(\omega\otimes\iota)(W\Omega^{*})$ is the weak limit of $\sum \overline{c_{i,j}}(L_{-\gamma_{i}}.\omega\otimes\iota)(W)L_{-\gamma_{j}}$,
so $\mathcal{B}\subset\widehat{M}\vee\{L_{\gamma}\,|\,\gamma\in\widehat{G}\}^{''}$. For the converse inclusion note that $(\omega\otimes\iota)(W)=(\omega\otimes\iota)(W\Omega^{*}\Omega)
$. Thus, $(\omega\otimes\iota)(W)$ is the weak limit of $\sum c_{i,j}(L_{\gamma_{i}}.\omega\otimes\iota)(W\Omega^{*})L_{\gamma_{j}}$.
The second assertion can be proved using the same technique.
\end{proof}

\begin{proposition}\label{IsoDual}
There exists a *-isomorphism $\rho\,:\,\widehat{G}^{2}\ltimes\widehat{M}\rightarrow \widehat{G}^{2}\ltimes\widehat{M_{\mathrm{\Omega}}}$ intertwining the actions $\theta^{\Psi}$
on $\widehat{G}^{2}\ltimes\widehat{M}$ and $\theta$ on $\widehat{G}^{2}\ltimes\widehat{M}_{\mathrm{\Omega}}$. Moreover,
$$\rho((\omega\otimes\iota)(\tilde{W}_{\mathrm{\Omega}}))=
\pi((\omega\otimes\iota)(W_{\mathrm{\Omega}})).$$
\end{proposition}

\begin{proof}
Remark that if we put $V=(\mathcal{F}\otimes\mathcal{F})U$ where $\mathcal{F}\,:\,\rightarrow L^{2}(G)$ is the Fourier transform and $U\,:\,L^{2}(\widehat{G}\times\widehat{G})\otimes H\rightarrow L^{2}(\widehat{G}\times\widehat{G})\otimes H$ is the unitary defined by $(U\xi)(\gamma_{1},\gamma_{2})=L_{\gamma_{1}}R_{\gamma_{2}}\xi(\gamma_{1},\gamma_{2})$ then
$$
\left\{
\begin{array}{rcl}
V\pi(x)V^{*} & = & 1\otimes 1\otimes x,\\
V\lambda_{\gamma,0}V^{*}&=&u_{\gamma}\otimes 1\otimes L_{\gamma},\\
V\lambda_{0,\gamma} V^{*}&=&1\otimes u_{\gamma}\otimes R_{\gamma}.\\
\end{array}
\right.
$$

Applying $\alpha\otimes\alpha\otimes\iota$, we conclude that the crossed
products can be defined on $H\otimes H\otimes H$ by:
\begin{eqnarray*}
\widehat{G}^{2}\ltimes\widehat{M} &=&
\{L_{\gamma}\otimes 1\otimes L_{\gamma}\,|\,\gamma\in\widehat{G}\}^{''}\vee
\{1\otimes L_{\gamma}\otimes R_{\gamma}\,|\,\gamma\in\widehat{G}\}^{''}\vee
1\otimes 1\otimes\widehat{M},\\
\widehat{G}^{2}\ltimes\widehat{M_{\mathrm{\Omega}}}&=&
\{L_{\gamma}\otimes 1\otimes L_{\gamma}\,|\,\gamma\in\widehat{G}\}^{''}\vee
\{1\otimes L_{\gamma}\otimes R_{\gamma}\,|\,\gamma\in\widehat{G}\}^{''}\vee
1\otimes 1\otimes\widehat{M_{\mathrm{\Omega}}}.\\
\end{eqnarray*}

Put $\mathcal{W}=(\widehat{J}\otimes\widehat{J})W(\widehat{J}\otimes\widehat{J})$.
Then $\mathcal{W}^{*}(1\otimes x)\mathcal{W}=\Delta^{\text{op}}(x)$, for all $x\in M$
and $[\mathcal{W}, 1\otimes y]=0$, for all $y\in\widehat{M}$. We have also $\mathcal{W}_{\Omega}=(\widehat{J_{\mathrm{\Omega}}}\otimes
\widehat{J_{\mathrm{\Omega}}})W{\mathrm{\Omega}}(\widehat{J_{\mathrm{\Omega}}}
\otimes\widehat{J_{\mathrm{\Omega}}})$ with similar properties.

In the following computation we use the relations $\mathcal{W}^{*}(1\otimes L_{\gamma})\mathcal{W}=L_{\gamma}\otimes L_{\gamma}$, $\mathcal{W}(1\otimes R_{\gamma})\mathcal{W}^{*}=L_{\gamma}\otimes R_{\gamma}$ and similar relations
with $\mathcal{W}_{\Omega}$. We use also the equality $[\mathcal{W}_{13}\Omega^{*}_{31},W_{23}\Omega^{*}_{23}]=0$ implying $[\mathcal{W}\Omega^{*}_{21},1\otimes y]=0$, for all $y\in\mathcal{B}$. Finally,
using Lemma \ref{TwistedGVN}, we have:
$$
\widehat{G}^{2}\ltimes\widehat{M} =
\{L_{\gamma}\otimes 1\otimes L_{\gamma}\,|\,\gamma\in\widehat{G}\}^{''}\vee
\{1\otimes L_{\gamma}\otimes R_{\gamma}\,|\,\gamma\in\widehat{G}\}^{''}\vee
1\otimes 1\otimes\widehat{M}
$$
$$
\downarrow\,\,\text{Ad}(\mathcal{W}_{13})
$$
$$
\begin{array}{lcr}
 &\{1\otimes 1\otimes L_{\gamma}\}^{''}\vee\{L_{\gamma}\otimes L_{\gamma}\otimes R_{\gamma}\}^{''}\vee1\otimes 1\otimes\widehat{M}&\\
 = &\{1\otimes 1\otimes L_{\gamma}\}^{''}\vee\{L_{\gamma}\otimes L_{\gamma}\otimes R_{\gamma}\}^{''}\vee1\otimes 1\otimes\widehat{B}&:=L_{1}\\
\end{array}
$$
$$
\downarrow\,\,\text{Ad}(\Omega_{32}\mathcal{W}^{*}_{23}\Omega_{31}\mathcal{W}^{*}_{13})
$$
$$
\begin{array}{lcr}
 &\{L_{\gamma}\otimes L_{\gamma}\otimes L_{\gamma}\}^{''}\vee\{1\otimes 1\otimes R_{\gamma}\}^{''}\vee 1\otimes 1\otimes\widehat{B}&\\
 = &\{L_{\gamma}\otimes L_{\gamma}\otimes L_{\gamma}\}^{''}\vee\{1\otimes 1\otimes R_{\gamma}\}^{''}\vee1\otimes 1\otimes\widehat{M_{\mathrm{\Omega}}}&:=L_{2}\\
\end{array}
$$
$$
\downarrow\,\,\text{Ad}((\mathcal{W}_{\Omega})_{23})
$$
$$
\begin{array}{cr}
\{L_{\gamma}\otimes 1\otimes L_{\gamma}\,|\,\gamma\in\widehat{G}\}^{''}\vee
\{1\otimes L_{\gamma}\otimes R_{\gamma}\,|\,\gamma\in\widehat{G}\}^{''}\vee
1\otimes 1\otimes\widehat{M_{\mathrm{\Omega}}}&=\widehat{G}^{2}\ltimes\widehat{M_{\mathrm{\Omega}}}.
\end{array}
$$
Define $\rho:=\rho_{2}\circ\varPhi\circ\rho_{1}$, where $\rho_{1}$, $\varPhi$ and $\rho_{2}$ are the isomorphisms from $\widehat{G}^{2}\ltimes\widehat{M}$ to $L_{1}$, from $L_{1}$ to $L_{2}$, and from $L_{2}$ to $\widehat{G}^{2}\ltimes\widehat{M_{\mathrm{\Omega}}}$, respectively. Then one can check that $\rho\circ\theta_{g_{1},g_{2}}^{\Psi}(x)=\theta_{g_{1},g_{2}}\circ\rho(x)$, for all $g_{1}, g_{2}\in G$ and for all $x$ of the form $\lambda_{\gamma_{1},\gamma_{2}}$ (or $L_{\gamma}\otimes 1\otimes L_{\gamma}$ and $1\otimes L_{\gamma}\otimes R_{\gamma}$ in our description of the crossed products). Thus, to finish the proof we only have to show that $(\omega\otimes\iota)(\tilde{W}_{\mathrm{\Omega}})\in N_{\Omega}$ and $\rho((\omega\otimes\iota)(\tilde{W}_{\mathrm{\Omega}}))=
\pi((\omega\otimes\iota)(W_{\mathrm{\Omega}}))$.
Using $(\ref{EqDual2})$, one computes
\begin{eqnarray}\notag
(\iota\otimes\theta_{(g_{1},g_{2})}^{\Psi})(\tilde{X})
&=&(\iota\otimes\theta_{(g_{1},g_{2})})(\tilde{X})
=(\alpha\otimes\lambda_{L})(\Psi^{*}(.,.-g_{1}))\\\notag
&=&(\alpha\otimes\lambda_{L})(\Psi_{g_{1}}
\otimes1)(\alpha\otimes\lambda_{L})(\Psi^{*})\\\label{EqDual4}
&=&(\alpha(\Psi_{g_{1}})\otimes 1)\tilde{X}.\\\notag
\end{eqnarray}
Similarly
\begin{eqnarray}\notag
(\iota\otimes\theta_{(g_{1},g_{2})}^{\Psi})(\tilde{Y})
&=&(\iota\otimes\theta_{(g_{1},g_{2})})(\tilde{Y})
=(\iota\otimes\theta_{(g_{1},g_{2})})((\alpha\otimes\lambda_{R})
(\tilde{\Psi}))\\\label{EqDual5}
&=&(\alpha\otimes\lambda_{R})
(\tilde{\Psi}(.,.-g_{2}))
=\tilde{Y}(\alpha(\Psi_{g_{2}})\otimes 1).
\end{eqnarray}
And, using (\ref{DualActHatM}) and (\ref{eqactionW}), one has
\begin{eqnarray}\notag
(\iota\otimes\theta_{(g_{1},g_{2})}^{\Psi})(\tilde{W})
&=&(\iota\otimes\pi)\left((L_{\Psi_{-g_{2}}}\otimes 1)W(L_{\Psi_{-g_{1}}}\otimes 1)\right)\\\label{EqDual6}
&=&(\alpha(\Psi^{*}_{g_{2}})\otimes 1)\tilde{W}(\alpha(\Psi^{*}_{g_{1}})\otimes 1).
\end{eqnarray}
Now $(\ref{EqDual4})$, $(\ref{EqDual5})$, and
$(\ref{EqDual6})$ imply
$
(\iota\otimes\theta_{(g_{1},g_{2})}^{\Psi})(\tilde{W}_{\mathrm{\Omega}})
=\tilde{W}_{\mathrm{\Omega}}
$, so $(\omega\otimes\iota)(\tilde{W}_{\mathrm{\Omega}})\in N_{\mathrm{\Omega}}$.
Now we want to show that $\rho((\omega\otimes\iota)(\tilde{W}_{\mathrm{\Omega}}))=
\pi((\omega\otimes\iota)(W_{\mathrm{\Omega}}))$. We take a net in the vector space spanned by elements $u_{\gamma_{1}}\otimes u_{\gamma_{2}}$ such that $\sum c_{i,j} (u_{\gamma_{i}}\otimes u_{\gamma_{j}})\rightarrow \Psi\quad\text{strongly*}$, so
$\sum\bar{c}_{i,j}(L_{-\gamma_{i}}\otimes\lambda_{-\gamma_{j},0})\rightarrow\tilde{X}$ and $
\sum\bar{c}_{i,j}(\alpha(u)\otimes 1)(L_{-\gamma_{i}}\otimes\lambda_{0,-\gamma_{j}})\rightarrow\tilde{Y}$
strongly*. This implies
$$
\sum\bar{c}_{i,j}\bar{c}_{k,l}\lambda_{0,-\gamma_{j}}
\pi((L_{-\gamma_{k}}.\omega.L_{-\gamma_{i}}.\alpha(u)\otimes\iota)(W))\lambda_{-\gamma_{l},0}
\rightarrow(\omega\otimes\iota)(\tilde{W}_{\mathrm{\Omega}})\quad\text{weakly}.
$$
Thus $\rho_{1}((\omega\otimes\iota)(\tilde{W}_{\mathrm{\Omega}}))$ is the weak limit of the net
\begin{eqnarray*}
&&\sum\bar{c}_{i,j}\bar{c}_{k,l}(L_{-\gamma_{j}}\otimes L_{-\gamma_{j}}\otimes R_{-\gamma_{j}})(1\otimes 1\otimes (L_{-\gamma_{k}}.\omega.L_{-\gamma_{i}}.\alpha(u)\otimes\iota)(W))(1\otimes 1\otimes L_{-\gamma_{l}})\\
&&=\sum_{i,j}\bar{c}_{i,j}(L_{-\gamma_{j}}\otimes L_{-\gamma_{j}}\otimes R_{-\gamma_{j}})
(1\otimes 1\otimes (\omega.L_{-\gamma_{i}}.\alpha(u)\otimes\iota)(W\sum_{k,l}\bar{c}_{k,l}L_{-\gamma_{k}}\otimes L_{-\gamma_{l}}))\\
&&\longrightarrow_{k,l}\sum_{i,j}\bar{c}_{i,j}(L_{-\gamma_{j}}\otimes L_{-\gamma_{j}}\otimes R_{-\gamma_{j}})
(1\otimes 1\otimes (\omega.L_{-\gamma_{i}}.\alpha(u)\otimes\iota)(W\Omega^{*})).
\end{eqnarray*}
and $\varPhi\circ\rho_{1}((\omega\otimes\iota)(\tilde{W}_{\mathrm{\Omega}}))$ is the weak limit of the net
\begin{eqnarray*}
&&\sum\bar{c}_{i,j}(1\otimes 1\otimes R_{0,-\gamma_{j}})(1\otimes 1\otimes (\omega.L_{-\gamma_{i}}.\alpha(u)\otimes\iota)(W\Omega^{*}))\\
&&=1\otimes 1\otimes (\omega\otimes\iota)(\sum\bar{c}_{i,j}(\alpha(u)\otimes 1)(L_{-\gamma_{i}}\otimes R_{-\gamma_{j}})W\Omega^{*}).
\end{eqnarray*}
Because $\sum\bar{c}_{i,j}(\alpha(u)\otimes 1)(L_{-\gamma_{i}}\otimes R_{-\gamma_{j}})\rightarrow Y$ weakly, we have
$$\varPhi\circ\rho_{1}((\omega\otimes\iota)(\tilde{W}_{\mathrm{\Omega}}))
=1\otimes 1\otimes (\omega\otimes\iota)(W_{\mathrm{\Omega}}).$$
This concludes the proof.
\end{proof}

In particular, Proposition \ref{IsoDual} implies that $N_{\Omega}=\{(\omega\otimes\iota)(\tilde{W}_{\Omega})\,|\,\omega\in\bh_{*}\}^{''}$ and that $\rho$ is a *-isomorphism from $N_{\Omega}$ to $\widehat{M}_{\Omega}$ which sends $(\omega\otimes\iota)(\tilde{W}_{\Omega})$ to $(\omega\otimes\iota)(W_{\Omega})$. Thus, we can transport the \lc{} quantum group structure from $\widehat{M}_{\Omega}$ to $N_{\Omega}$. First, we show that the comultiplication introduced in Section \ref{SectionResult} is the good one. For this we need
\begin{proposition}
There exists a unique unital normal *-homomorphism $\Gamma\,:\, N\rightarrow N\otimes N$ such that
$$\Gamma(\lambda_{\gamma_{1},\gamma_{2}})=\lambda_{\gamma_{1},0}\otimes
\lambda_{0,\gamma_{2}}\quad\text{and}\quad\Gamma(\pi(x))=(\pi\otimes\pi)\hat{\Delta}(x).$$
\end{proposition}

\begin{proof}
Like in the begining of the proof of Proposition \ref{IsoDual} define the crossed product
$$
\widehat{G}^{2}\ltimes\widehat{M} =
\{L_{\gamma}\otimes 1\otimes L_{\gamma}\,|\,\gamma\in\widehat{G}\}^{''}\vee
\{1\otimes L_{\gamma}\otimes R_{\gamma}\,|\,\gamma\in\widehat{G}\}^{''}\vee
1\otimes 1\otimes\widehat{M}.
$$
Let $\mathcal{W}$ be the operator defined in the proof of Proposition \ref{IsoDual} and $Q$
be the unitary on $H\otimes H\otimes H\otimes H\otimes H\otimes H$ such that $Q^{*}=\Sigma_{45}\Sigma_{35}\mathcal{W}^{*}_{15}\hat{W}^{*}_{56}\Sigma_{45}$.
We define $\Gamma(x)=Q^{*}(1\otimes x)Q$. Using that $\hat{W}^{*}(L_{\gamma}\otimes L_{\gamma})\hat{W}=L_{\gamma}\otimes 1$, $\hat{\Delta}(x)=\hat{W}^{*}(1\otimes x)\hat{W}$, for all $x\in\hat{M}$, $[\hat{W},1\otimes y]$=0, for all $y\in M^{'}$, $\mathcal{W}^{*}(1\otimes L_{\gamma})\mathcal{W}= L_{\gamma}\otimes L_{\gamma}$ and $[\mathcal{W},1\otimes y]=0$, for all $y\in\hat{M}$, one can check that the needed properties of $\Gamma$.
\end{proof}

The unitary $\Upsilon=(\lambda_{R}\otimes\lambda_{L})(\tilde{\Psi}^{*})\in
N\otimes N$ allows to define the unital normal *-homomorphism $\Gamma_{\Omega}(x)=\Upsilon\Gamma(x)\Upsilon^{*}:N\to N\otimes N$ which is
a comultiplication on $N_{\mathrm{\Omega}}$:

\begin{proposition}\label{comultiplication}
For all $x\in N_{\mathrm{\Omega}}$, we have $\Gamma_{\mathrm{\Omega}}(x)\in N_{\mathrm{\Omega}}\otimes N_{\mathrm{\Omega}}$ and
$$(\rho\otimes\rho)(\Gamma_{\mathrm{\Omega}}(x))=\hat{\Delta}_{\mathrm{\Omega}}(\rho(x)).$$
\end{proposition}

\begin{proof}
It suffices to show that $(\iota\otimes\rho\otimes\rho)(\iota\otimes\Gamma_{\Omega})(\tilde{W}_{\Omega})=
(W_{\Omega})_{13}(W_{\Omega})_{12}$. By the definition of $\Gamma$, one has, for any $F\in L^{\infty}(G)$,
$$
\Gamma(\lambda_{L}(F))=\lambda_{L}(F)\otimes 1\quad\text{and}\quad\Gamma(\lambda_{R}(F))=1\otimes\lambda_{R}(F),
$$
and since $1\otimes\Upsilon$ commutes with $\tilde{X}_{12}$ and with
$\tilde{Y}_{13}$, one gets:
$$
(\iota\otimes\Gamma_{\Omega})(\tilde{X})=\tilde{X}_{12}\quad\text{and}\quad
(\iota\otimes\Gamma_{\Omega})(\tilde{Y})=\tilde{Y}_{13}.
$$
Moreover,
\begin{eqnarray*}
(\iota\otimes\Gamma_{\Omega})(\tilde{W})
&=&(1\otimes\Upsilon)(\iota\otimes\Gamma\circ\pi)(W)(1\otimes\Upsilon^{*})\\
&=&(1\otimes\Upsilon)(\iota\otimes\pi\otimes\pi)\left((\iota\otimes
\hat{\Delta})(W)\right)(1\otimes\Upsilon^{*})\\
&=&(1\otimes\Upsilon)(\iota\otimes\pi\otimes\pi)(W_{13}W_{12})(1\otimes\Upsilon^{*})\\
&=&\Upsilon_{23}\tilde{W}_{13}\tilde{W}_{12}\Upsilon_{23}^{*}.
\end{eqnarray*}
Using (\ref{EqWL}), we can check the following relations on the generators $u_{\gamma}$ of $L^{\infty}(G)$:
\begin{eqnarray*}
&&W(1\otimes\pi_{R}(F))W^{*}=(\pi_{L}\otimes\pi_{R})(\Delta_{G}(F))\,,\\
&&W^{*}(1\otimes\pi_{L}(F))W=(\pi_{L}\otimes\pi_{L})(\Delta_{G}(F))\,,
\text{for any}\,F\in L^{\infty}(G).\\
\end{eqnarray*}
Then
\begin{eqnarray*}
&&W_{12}(\pi_{R}\otimes\pi_{L})(\tilde{\Psi}^{*})_{23}^{*}W_{12}^{*}=(\pi_{L}\otimes\pi_{R}\otimes\pi_{L})\left((\Delta_{G}\otimes\iota)(\tilde{\Psi})\right)\,\,,\\
&&W_{13}^{*}(\pi_{R}\otimes\pi_{L})(\tilde{\Psi}^{*})_{23}W_{13}=(\pi_{L}\otimes\pi_{R}\otimes\pi_{L})
\left((\sigma\otimes\iota)\left((\iota\otimes\Delta_{G})(\tilde{\Psi}^{*})\right)\right).
\end{eqnarray*}
Let us define
$$V=(\pi_{L}\otimes\pi_{R}\otimes\pi_{L})\left((\sigma\otimes\iota)\left((\iota\otimes\Delta_{G})(\tilde{\Psi}^{*})\right)\right)
(\pi_{L}\otimes\pi_{R}\otimes\pi_{L})\left((\Delta_{G}\otimes\iota)(\tilde{\Psi})\right),$$
then we have
\begin{eqnarray*}
(\iota\otimes\rho\otimes\rho)(\iota\otimes\Gamma_{\Omega})(\tilde{W}_{\Omega})
&=&(\iota\otimes\rho\otimes\rho)(\tilde{Y}_{13}\Upsilon_{23}\tilde{W}_{13}\tilde{W}_{12}\Upsilon_{23}^{*}\tilde{X}_{12})\\
&=&Y_{13}W_{13}VW_{12}X_{12},
\end{eqnarray*}
so it remains to calculate:
\begin{eqnarray*}
&&(\sigma\otimes\iota)\left((\iota\otimes\Delta_{G})(\tilde{\Psi}^{*})\right)(g,h,t)
(\Delta_{G}\otimes\iota)(\tilde{\Psi})(g,h,t)\\
&&=(\iota\otimes\Delta_{G})(\tilde{\Psi}^{*})(h,g,t)(\Delta_{G}\otimes\iota)(\tilde{\Psi})(g,h,t)\\
&&=\Psi^{*}(-h,h+g+t)\Psi(-g-h,g+h+t)=\Psi(-g,g+h+t)\\
&&=\Psi^{*}(g,t)\tilde{\Psi}(g,h).\\
\end{eqnarray*}
Thus, $V=X_{13}Y_{12}$, and this concludes the proof.
\end{proof}

\begin{remark}\label{RemarkCstar}
One can show that $\alpha$ and $\beta$ are actions of $\widehat{G}^{2}$ on the reduced dual $\cs$-algebras $\widehat{A}$ and $\widehat{A}_{\mathrm{\Omega}}$. Moreover, the *-isomorphism $\rho$ gives a *-isomorphism between the reduced crossed products $\widehat{G}^{2}\ltimes\widehat{A}$ and $\widehat{G}^{2}\ltimes\widehat{A}_{\mathrm{\Omega}}$. So $\hat{A}$ is nuclear if and only if $\hat{A}_{\mathrm{\Omega}}$ is nuclear. Moreover, the twisted dual action $\theta^{\Psi}$ gives a deformed $\widehat{G}^{2}$-product structure on $\widehat{G}^{2}\ltimes\widehat{A}$ and the Landstad algebra for this $\widehat{G}^{2}$-product is $[(\omega\otimes\iota)(\tilde{W}_{\mathrm{\Omega}})]$, and it is isomorphic to $\widehat{A}_{\Omega}$. These results can be obtained directly from the universality property of crossed products (see \cite{Fima}).
\end{remark}

The rest of this section is devoted to the computation of the left invariant weight on $(N_{\Omega},\Gamma_{\Omega})$. Since $\rho\,:\,N=\hat{G}^{2}\ltimes\hat{M}\rightarrow\hat{G}^{2}\ltimes
\hat{M}_{\Omega}$ is a $*$-isomorphism, one can consider two natural weights on $N$, $\varphi_{1}=\tilde{\hat{\varphi}}$, the dual weight of $\hat{\varphi}$ on $N$, and $\varphi_{2}=\tilde{\hat{\varphi}}_{\Omega}\circ\rho$, where $\tilde{\hat{\varphi}}_{\Omega}$ is the dual weight of $\hat{\varphi}_{\Omega}$ on $\hat{G}^{2}\ltimes
\hat{M}_{\Omega}$.

\begin{lemma}\label{LemRN}
We have:
\begin{enumerate}
\item $[D\hat{\varphi}\circ\alpha_{\gamma_{1},\gamma_{2}}\,:\,D\hat{\varphi}]_{t}=\langle\gamma_{2},\gamma_{t}\rangle
=[D\hat{\varphi}_{\Omega}\circ\beta_{\gamma_{1},\gamma_{2}}\,:\,D\hat{\varphi}_{\Omega}]_{t}$ $\forall t\in\R$, $\forall \gamma_{1},\gamma_{2}\in\widehat{G}$.\label{LemRN3}
\item $[D\varphi_{1}\circ\theta^{\Psi}_{g_{1},g_{2}}\,:\,D\varphi_{1}]_{t}=\Psi(\gamma_{t},g_{2})$, for all $t\in\Rset$ and all $g_{1},g_{2}\in G$.\label{LemRN4}
\item For any \nsf{} weight $\nu$ on $N$, $\nu$ is invariant under the action $\theta^{\Psi}$ if and only if $\theta^{\Psi}_{g_{1},g_{2}}\left([D\nu\,:\,D\varphi_{1}]_{t}\right)=\Psi(\gamma_{t},g_{2})[D\nu\,:\,D\varphi_{1}]_{t}$.\label{LemRN5}
\end{enumerate}
\end{lemma}

\begin{proof}
Using Proposition \ref{PropCo-Subgroups}(2), and because $L_{\gamma}$ and $R_{\gamma}$ are unitaries, we find $\hat{\varphi}\circ\alpha^{L}_{\gamma}=
\hat{\varphi},\quad\hat{\varphi}\circ\alpha^{R}_{\gamma}=\lambda(\gamma)\hat{\varphi}$, so
\begin{eqnarray*}
[D\hat{\varphi}\circ\alpha_{\gamma_{1},\gamma_{2}}\,:\,D\hat{\varphi}]_{t}
&=&[D\hat{\varphi}\circ\alpha_{\gamma_{1}}^{L}\circ\alpha_{\gamma_{2}}^{R}\,:\,D\hat{\varphi}\circ\alpha_{\gamma_{2}}^{R}]_{t}[D\hat{\varphi}\circ\alpha_{\gamma_{2}}^{R}\,:\,D\hat{\varphi}]_{t}\\
&=&\alpha_{-\gamma_{2}}^{R}\left([D\hat{\varphi}\circ\alpha_{\gamma_{1}}^{L}\,:\,D\hat{\varphi}]_{t}\right)
[D\hat{\varphi}\circ\alpha_{\gamma_{2}}^{R}\,:\,D\hat{\varphi}]_{t}\\
&=&\lambda(\gamma_{2})^{it}=\langle\gamma_{2},\gamma_{t}\rangle.
\end{eqnarray*}
The right-hand side of the first equality is obtained by considering the stable co-subgroup $\widehat{G}<(M,\Delta_{\Omega})$. Let us prove the second assertion. Let $g_{1},g_{2}\in G$, define the unitary $v:=\lambda_{(\Psi_{-g_{1}},\Psi_{g_{2}})}$, and denote by $\varphi_{1}|_{v}$ the weight $\varphi_{1}|_{v}(x)=\varphi_{1}(vxv^{*})$. Using the first assertion, we have
$$
[D\varphi_{1}|_{v}\,:\,D\varphi_{1}]_{t}
=v^{*}\sigma_{t}^{1}(v)=v^{*}\langle\Psi_{g_{2}},\gamma_{t}\rangle v=\Psi(\gamma_{t},g_{2}).
$$
Note that $\varphi_{1}\circ\theta^{\Psi}_{g_{1},g_{2}}=\varphi_{1}|_{v}\circ\theta_{g_{1},g_{2}}$, so
\begin{eqnarray*}
[D\varphi_{1}\circ\theta^{\Psi}_{g_{1},g_{2}}\,:\,D\varphi_{1}]_{t}
&=&[D\varphi_{1}|_{v}\circ\theta_{g_{1},g_{2}}\,:\,D\varphi_{1}\circ\theta_{g_{1},g_{2}}]_{t}
[D\varphi_{1}\circ\theta_{g_{1},g_{2}}\,:\,D\varphi_{1}]_{t}\\
&=&\theta_{-g_{1},-g_{2}}\left([D\varphi_{1}|_{v}\,:\,D\varphi_{1}]_{t}\right)
=\Psi(\gamma_{t},g_{2}).
\end{eqnarray*}
Putting $u_{t}=[D\nu\,:\,D\varphi_{1}]_{t}$ and using the second assertion, one has
$$
[D\nu\circ\theta^{\Psi}_{g_{1},g_{2}}\,:\,D\nu]_{t}
=\theta^{\Psi}_{-g_{1},-g_{2}}(u_{t})[D\varphi_{1}\circ\theta^{\Psi}_{g_{1},g_{2}}\,:\,D\varphi_{1}]_{t}
u_{t}^{*}
=\theta^{\Psi}_{-g_{1},-g_{2}}(u_{t})\Psi(\gamma_{t},g_{2})u_{t}^{*}.
$$
This concludes the proof.
\end{proof}

Note that, using Lemma \ref{LemRN} $(\ref{LemRN3})$, we have, for all $t\in\Rset,F\in L^{\infty}(G^{2})$,
\begin{equation}\label{EqMod1}
\sigma_{t}^{1}(\lambda(F))=\lambda(F(\cdot,\cdot+\gamma_{t}))=\sigma_{t}^{2}(\lambda(F)).
\end{equation}
Let $T$ be the strictly positive operator affiliated with $N$ and such that $T^{it}=\lambda_{R}(\Psi(-\gamma_{t},.))$. Using $(\ref{EqMod1})$, we find $\sigma_{t}^{1}(T^{is})=\lambda^{-its}T^{is}$, so one can consider the Vaes' weight $\tilde{\mu}_{\Omega}$ associated with $T$ and $\lambda^{-1}$. This is the unique \nsf{}
weight on $N$ such that $[D\tilde{\mu}_{\Omega}\,:\,D\varphi_{1}]_{t}=\lambda^{\frac{-it^{2}}{2}}T^{it}$.
From $(\ref{EqDual2})$ we have $\theta^{\Psi}_{g_{1},g_{2}}(T^{it})=\lambda_{R}(\Psi(-\gamma_{t},.-g_{2})=
\Psi(\gamma_{t},g_{2})T^{it}$. By Lemma $\ref{LemRN}$(3), $\tilde{\mu}_{\Omega}$ is invariant under $\theta^{\Psi}$,so the image $\tilde{\mu}_{\Omega}\circ\rho^{-1}$ of $\tilde{\mu}_{\Omega}$ in $\hat{G}^{2}\ltimes\hat{M}_{\Omega}$ is invariant under the dual action. Thus, $\tilde{\mu}_{\Omega}\circ\rho^{-1}$ is the dual weight of some weight $\mu_{\Omega}$ on $\hat{M}_{\Omega}$. To finish the proof of Theorem \ref{TheoremRieffel},
we will show in Theorem \ref{dualweight} that $\mu_{\Omega}=\hat{\varphi}_{\Omega}$, for
which we need

\begin{proposition}
For all $t\in\Rset$ and all $x\in N$, we have
$$\sigma_{t}^{2}(x)=T^{it}\sigma_{t}^{1}(x) T^{-it}.$$
\end{proposition}

\begin{proof}
By $(\ref{EqMod1})$, it suffices to prove this equality for elements of the form $(\omega\otimes\iota)(\tilde{W}_{\Omega})$. By definition of $\sigma_{t}^{2}$, we have
$$\sigma_{t}^{2}\left((\omega\otimes\iota)(\tilde{W}_{\Omega})\right)=
(\rho_{t}^{\Omega}(\omega)\otimes\iota)(\tilde{W}_{\Omega}),$$
where $\rho_{t}^{\Omega}(\omega)(x)=\omega(\delta_{\Omega}^{-it}\tau^{\Omega}_{-t}(x))$. Proposition \ref{TwistedIng} gives
$$
\rho_{t}^{\Omega}(x)=\omega\left(\delta^{-it}A^{it}B^{-it}\tau_{t}(x)\right).
$$
On the other hand, using $(\ref{EqMod1})$, one has
$$(\iota\otimes\sigma_{t}^{1})(\tilde{X})=\tilde{X},\quad (\iota\otimes\sigma_{t}^{2})(\tilde{Y})=(A^{it}\otimes 1)\tilde{Y},$$
which implies
\begin{eqnarray*}
(\iota\otimes\sigma_{t}^{1})(\tilde{W}_{\Omega}) &=&
(A^{it}\otimes 1)\tilde{Y}(\rho_{t}(\omega)\otimes\iota)(\tilde{W})\tilde{X}\\
&=&(A^{it}\otimes 1)\tilde{Y}(\delta^{-it}\otimes 1)(\tau_{-t}\otimes\iota)(\tilde{W})\tilde{X}\\
&=&(B^{it}\otimes 1)(\delta^{-it}A^{it}B^{-it}\otimes 1)(\tau_{-t}\otimes\iota)(\tilde{W}_{\Omega})\\
& & \text{because}\,\,\delta^{it}\alpha(.)\delta^{-it}=\alpha(.)\,\,\text{and}\,\,\tau_{t}\circ\alpha=\tau_{t}\\
&=&(B^{it}\otimes 1)(\iota\otimes\sigma_{t}^{2})(\tilde{W}_{\Omega}).
\end{eqnarray*}
Next, using $(\ref{EqWL})$ with the character $\chi_{t}(g)=\Psi(\gamma_{t},g)$, we find
\begin{eqnarray}\notag
(B^{it}\otimes 1)\tilde{W}_{\Omega} &=& \tilde{Y}(L_{-\chi_{t}}\otimes 1)\tilde{W}\tilde{X}
= \tilde{Y}(\iota\otimes\pi)\left((L_{-\chi_{t}}\otimes 1)W\right)\tilde{X} \\\notag
&=& \tilde{Y}(\iota\otimes\pi)\left((1\otimes R_{\chi_{t}})W(1\otimes R_{\chi_{t}}^{*})\right)\tilde{X}
= (1\otimes\lambda_{R}(\chi_{t}))\tilde{W}_{\Omega}(1\otimes\lambda_{R}(\chi_{t})^{*})\\
&=&(1\otimes T^{-it})\tilde{W}_{\Omega}(1\otimes T^{it}).\label{WOmegaT}
\end{eqnarray}
Thus, for all $t\in\Rset,\omega\in M_{*}$, one has
\begin{eqnarray*}
\sigma_{t}^{1}((\omega\otimes\iota)(\tilde{W}_{\Omega}))
&=& (\omega\otimes\iota)((B^{it}\otimes 1)(\iota\otimes\sigma_{t}^{2})(\tilde{W}_{\Omega}))\\
&=& (\omega\otimes\iota)((\iota\otimes\sigma_{t}^{2})((1\otimes T^{-it})\tilde{W}_{\Omega}(1\otimes T^{it})))\\
&=& T^{-it}\sigma_{t}^{2}((\omega\otimes\iota)(\tilde{W}_{\Omega}))T^{it},
\end{eqnarray*}
where we used, in the last equation, $\sigma_{t}^{2}(T^{is})=\lambda^{-ist}T^{is}$. 
\end{proof}

\begin{theorem}\label{dualweight}
We have $\mu_{\Omega}=\hat{\varphi}_{\Omega}$.
\end{theorem}

Let us denote by $\varphi_{P}$ the Plancherel weight on $\mathcal{L}(\hat{G}^{2})$, by $\Lambda_{P}$ its canonical \GNS{} map, by $\lambda_{L}^{\hat{G}^{2}}$ and $\lambda_{R}^{\hat{G}^{2}}$ the $*$-homorphisms $L^{\infty}(G)\rightarrow \mathcal{L}(\hat{G}^{2})$ coming from the representations $(\gamma\mapsto\lambda^{\hat{G}^{2}}_{(\gamma,0)})$ and $(\gamma\mapsto\lambda^{\hat{G}^{2}}_{(0,\gamma)})$, respectively, and by $T_{1}^{it}=\lambda_{R}^{\hat{G}^{2}}(\Psi(-\gamma_{t},\cdot)$. Thus, $T=T_{1}\otimes 1$.
We also introduce the $*$-homomorphism $\alpha^{'}(F)=J\alpha(F)^{*}J$ and denote by
$F\mapsto F^{\circ}$ the $*$-automorphism of $L^{\infty}(G\times G)$ defined by $F^{\circ}(g,h)=F(h,g+h)$.

The standard \GNS{} construction for $\varphi_{1}$ is $(L^{2}(\hat{G}^{2},H),\iota,\Lambda_{1})$, where a $\sigma$-strong-*-norm core for $\Lambda_{1}$ is given by
$$\mathcal{D}_{1}=\left\{(x\otimes 1)(\omega\otimes\iota)(\tilde{W})\,|\,x\in\mathcal{N}_{\varphi_{P}},\,\omega
\in\mathcal{I}\right\},$$
and, if $x\in\mathcal{N}_{\varphi_{P}},\,\omega\in\mathcal{I}$, we have
$$\Lambda_{1}\left((x\otimes 1)(\omega\otimes\iota)(\tilde{W})\right)=\Lambda_{P}(x)\otimes\xi(\omega).$$
Let $\lambda_{\Omega}(\omega)$, $\mathcal{I}_{\Omega}$, and $\xi_{\Omega}$ be the standard objects associated with $(M,\Delta_{\Omega})$. For $\varphi_{2}$, we take the \GNS{} construction $(L^{2}(\hat{G}^{2},H),\tilde{\rho},\Lambda_{2})$, where a $\sigma$-strong-*-norm core for $\Lambda_{2}$ is
$$\mathcal{D}_{2}=\left\{(x\otimes 1)(\omega\otimes\iota)(\tilde{W}_{\Omega})\,|\,x\in\mathcal{N}_{\varphi_{P}},\,\omega
\in\mathcal{I}_{\Omega}\right\},$$
and, if $x\in\mathcal{N}_{\varphi_{P}},\,\omega\in\mathcal{I}_{\Omega}$, one has
$$\Lambda_{2}\left((x\otimes 1)(\omega\otimes\iota)(\tilde{W}_{\Omega})\right)=\Lambda_{P}(x)\otimes\xi_{\Omega}(\omega).$$

Let us introduce the following sets:
\begin{eqnarray*}
C_{1} & = & \left\{x\in\mathcal{N}_{\varphi_{P}}\,|\,T^{\frac{1}{2}}(x\otimes 1)\,\,\,\text{is bounded}\right\},\\
C_{1}^{0} & = & \left\{x\in C_{1}\,|\,\Lambda_{P}(x)\in\mathcal{D}(T_{1}^{-\frac{1}{2}})\right\},\\
C_{2} & = & \left\{\omega_{\xi,\eta}\in\mathcal{I}_{\Omega}\,|\,\eta\in\mathcal{D}(A^{-\frac{1}{2}})
\cap\mathcal{D}(B^{\frac{1}{2}})\right\}.
\end{eqnarray*}

\begin{lemma}\label{LemCore}
For all $\omega_{\xi,\eta}\in C_{2}$ one has $\omega_{\xi,A^{-\frac{1}{2}}\eta},\omega_{\xi,u^{*}B^{\frac{1}{2}}\eta}\in\mathcal{I}$.
Moreover,
$$\xi_{\Omega}\left(\omega_{\xi,\eta}\right)=\xi\left(\omega_{\xi, A^{-\frac{1}{2}}\eta}\right),\quad
\xi\left(\omega_{\xi,u^{*}B^{\frac{1}{2}}\eta}\right)=
\lambda^{\frac{i}{4}}Ju^{*}J\xi\left(\omega_{\xi,A^{-\frac{1}{2}}\eta}\right).$$
The following set is a $\sigma$-weak-weak core for $\Lambda_{2}$:
$$\mathcal{D}=\left\{(x\otimes 1)(\omega_{\xi,\eta}\otimes\iota)(\tilde{W}_{\Omega})\,|\,x\in C^{0}_{1},\,\omega_{\xi,\eta}\in C_{2}\right\}.$$
Moreover, if $x\in C_{1}$ and $\omega_{\xi,\eta}\in C_{2}$, then
$$\Lambda_{2}((x\otimes 1)(\omega_{\xi,\eta}\otimes\iota)(\tilde{W}_{\Omega}))=
\Lambda_{P}(x)\otimes\xi(\omega_{\xi,A^{-\frac{1}{2}}\eta}).$$
\end{lemma}

\begin{proof}
Let $\omega_{\xi,\eta}\in\mathcal{I}_{\Omega}$ and $\eta\in\mathcal{D}(A^{-\frac{1}{2}})$. Let $e_{n}$ be self-adjoint elements, like in Lemma \ref{LemVaesRN}, for the operator $A$. When $x\in\mathcal{N}_{\varphi}$, we have
\begin{eqnarray*}
|\omega_{\xi,A^{-\frac{1}{2}}\eta}(e_{n}x^{*})| &=& |\langle e_{n}x^{*}\xi,A^{-\frac{1}{2}}\eta\rangle |=
|\langle (A^{-\frac{1}{2}} e_{n})x^{*}\xi,\eta\rangle | \\
 &=& |\langle (xA^{-\frac{1}{2}}e_{n})^{*}\xi,\eta\rangle | \leq C||\Lambda_{\Omega}(xA^{-\frac{1}{2}}e_{n})||,
\end{eqnarray*}
because $xA^{-\frac{1}{2}}e_{n}A^{\frac{1}{2}}$ is bounded and its closure,
which equals to $xe_{n}$, belons to $\mathcal{N}_{\varphi}$. Thus, we obtain
$$|\omega_{\xi,A^{-\frac{1}{2}}\eta}(e_{n}x^{*})|\leq C||\Lambda(xe_{n})||=C||J\sigma_{-\frac{i}{2}}(e_{n})J\Lambda(x)||\rightarrow C||\Lambda(x)||.$$
Since $|\omega_{\xi,A^{-\frac{1}{2}}\eta}(e_{n}x^{*})|\rightarrow |\omega_{\xi,A^{-\frac{1}{2}}\eta}(x^{*})|$, we conclude that $\omega_{\xi, A^{-\frac{1}{2}}\eta}$ is in $\mathcal{I}$. Moreover, for all $x\in\mathcal{N}_{\varphi}$, we have
\begin{eqnarray*}
\langle \xi_{\Omega}(\omega_{\xi,\eta}),J\sigma_{-\frac{i}{2}}(e_{n})J\Lambda(x)\rangle
&=&\langle \xi_{\Omega}(\omega_{\xi,\eta}),\Lambda(xe_{n})\rangle
=\langle \xi_{\Omega}(\omega_{\xi,\eta}),\Lambda_{\Omega}(xA^{-\frac{1}{2}}e_{n})\rangle\\
&=&\omega_{\xi,\eta}(A^{-\frac{1}{2}}e_{n}x^{*})
=\langle e_{n}x^{*}\xi,A^{-\frac{1}{2}}\eta\rangle\\
&=&\omega_{\xi,A^{-\frac{1}{2}}\eta}(e_{n}x^{*})
=\langle \xi(\omega_{\xi,A^{-\frac{1}{2}}\eta}),\Lambda(xe_{n})\rangle\\
&=&\langle \xi(\omega_{\xi,A^{-\frac{1}{2}}\eta}),J\sigma_{-\frac{i}{2}}(e_{n})J\Lambda(x)\rangle.
\end{eqnarray*}
Taking the limit when $n\rightarrow\infty$, we conclude that $\xi_{\Omega}(\omega_{\xi,\eta})=\xi(\omega_{\xi,A^{-\frac{1}{2}}\eta})$.

Suppose that $\eta\in\mathcal{D}(B^{\frac{1}{2}})$. Let $f_{m}$ be self-adjoint elements, like in Lemma \ref{LemVaesRN}, for the operator $B$. Note that $f_{m}$ commute with $e_{n}$ and $u$, also $e_{n}$ commute with $u$. Let us show that $uB^{\frac{1}{2}}f_{m}A^{\frac{1}{2}}e_{n}$ is analytic w.r.t. $\sigma$. We have
$$\Psi(-(g-\gamma_{t}),g-\gamma_{t})=\lambda^{-it^{2}}\Psi(-g,g)
\Psi(g,\gamma_{t})\Psi(\gamma_{t},g),$$
so $\sigma_{t}(u)=\lambda^{-it^{2}}uA^{-it}B^{-it}$ and, using Lemma \ref{LemVaesRN}, we obtain
\begin{eqnarray*}
\sigma_{t}(uB^{\frac{1}{2}}f_{m}A^{\frac{1}{2}}e_{n})
&=& \lambda^{-it^{2}}uA^{-it}B^{-it}B^{\frac{1}{2}}\sigma_{t}(f_{m})A^{\frac{1}{2}}
\sigma_{t}(e_{n})\\
&=& \lambda^{-it^{2}}uB^{\frac{1}{2}-it}\sigma_{t}(f_{m})A^{\frac{1}{2}-it}\sigma_{t}(e_{n}).
\end{eqnarray*}
Define the following function from $\Cset$ to $M$:
$$f(z)=\lambda^{-iz^{2}}uB^{\frac{1}{2}-iz}\sigma_{z}(f_{m})A^{\frac{1}{2}-iz}\sigma_{z}(e_{n}).$$
By Lemma \ref{LemVaesRN}, $f$ is analytic, so $uB^{\frac{1}{2}}f_{m}A^{\frac{1}{2}}e_{n}$ is analytic, and we have
$$\sigma_{-\frac{i}{2}}(uB^{\frac{1}{2}}f_{m}A^{\frac{1}{2}}e_{n})
=\lambda^{\frac{i}{4}}u\sigma_{-\frac{i}{2}}(f_{m})\sigma_{-\frac{i}{2}}(e_{n}).$$

Thus, for $x\in\mathcal{N}_{\varphi}$, $xu^{*}B^{\frac{1}{2}}f_{m}e_{n}A^{\frac{1}{2}}$ is bounded and its closure, which is equal to $xu^{*}B^{\frac{1}{2}}f_{m}A^{\frac{1}{2}}e_{n}$,
belongs to $\mathcal{N}_{\varphi}$. Moreover,
\begin{eqnarray*}
| \omega_{\xi,u^{*}B^{\frac{1}{2}}\eta}(e_{n}f_{m}x^{*})| & = & |\langle e_{n}f_{m}x^{*}\xi,u^{*}B^{\frac{1}{2}}\eta\rangle |
 =  | \langle B^{\frac{1}{2}}f_{m}e_{n}ux^{*}\xi,\eta\rangle |\\
&=& | \langle \left(xu^{*}B^{\frac{1}{2}}f_{m}e_{n}\right)^{*}\xi,\eta\rangle |
\leq C||\Lambda(xu^{*}B^{\frac{1}{2}}f_{m}A^{\frac{1}{2}}e_{n})||\\
&\leq &C||J\lambda^{\frac{i}{4}}u\sigma_{-\frac{i}{2}}(f_{m})\sigma_{-\frac{i}{2}}(e_{n})J\Lambda(x)||.
\end{eqnarray*}
Taking the limit over $m$ and $n$, we get
$$|\omega_{\xi,u^{*}B^{\frac{1}{2}}\eta}(x^{*})|\leq C||JuJ\Lambda(x)||\leq C||u||||\Lambda(x)||.$$
Thus, $\omega_{\xi,u^{*}B^{\frac{1}{2}}\eta}\in\mathcal{I}$. Moreover, for all $x\in\mathcal{N}_{\varphi}$, one has
$$
\begin{array}{l}
\langle \xi(\omega_{\xi,u^{*}B^{\frac{1}{2}}\eta}),J\sigma_{-\frac{i}{2}}(e_{n})\sigma_{-\frac{i}{2}}(f_{m})J\Lambda(x)\rangle
=\langle \xi(\omega_{\xi,u^{*}B^{\frac{1}{2}}\eta}),\Lambda(xe_{n}f_{m})\rangle\\
=\omega_{\xi,u^{*}B^{\frac{1}{2}}\eta}(e_{n}f_{m}x^{*})
=\langle e_{n}f_{m}x^{*}\xi,u^{*}B^{\frac{1}{2}}\eta\rangle
=\langle x^{*}\xi,u^{*}B^{\frac{1}{2}}f_{m}A^{\frac{1}{2}}e_{n}A^{-\frac{1}{2}}\eta\rangle\\
=\langle uB^{\frac{1}{2}}f_{m}A^{\frac{1}{2}}e_{n}x^{*}\xi,A^{-\frac{1}{2}}\eta\rangle
=\omega_{\xi,A^{-\frac{1}{2}}\eta}\left((xu^{*}B^{\frac{1}{2}}f_{m}A^{\frac{1}{2}}e_{n})^{*}\right)\\
=\langle\xi(\omega_{\xi,A^{-\frac{1}{2}}\eta}),\Lambda(xu^{*}B^{\frac{1}{2}}f_{m}A^{\frac{1}{2}}e_{n})\rangle
=\langle\xi(\omega_{\xi,A^{-\frac{1}{2}}\eta}),J\lambda^{\frac{i}{4}}u\sigma_{-\frac{i}{2}}(f_{m})\sigma_{-\frac{i}{2}}(e_{n})J\Lambda(x)\rangle\\
=\langle\lambda^{\frac{i}{4}}Ju^{*}J\xi(\omega_{\xi,A^{-\frac{1}{2}}\eta}),J\sigma_{-\frac{i}{2}}(e_{n})\sigma_{-\frac{i}{2}}(f_{m})J\Lambda(x)\rangle.
\end{array}
$$

Taking the limit over $m$ and $n$, we get $\xi(\omega_{\xi,u^{*}B^{\frac{1}{2}}\eta})=
\lambda^{\frac{i}{4}}Ju^{*}J\xi(\omega_{\xi,A^{-\frac{1}{2}}\eta})$. Now we want to prove that $\mathcal{D}$ is a $\sigma$-weak-weak core for $\Lambda_{2}$. Because $T=T_{1}\otimes 1$, we know that $T^{\frac{1}{2}}(x\otimes1)$ is bounded if and only if $T_{1}^{\frac{1}{2}}x$ is bounded. Thus, by Proposition \ref{PropTx}, $C_{1}^{0}$ is a $\sigma$-strong*-norm core for $\Lambda_{P}$, and, by Proposition \ref{PropCor}, it suffices to show that the set $\{(\omega\otimes\iota)(W_{\Omega})\,|\,\omega\in C_{2}\}$ is a $\sigma$-strong*-norm core
for $\hat{\Lambda}_{\Omega}$. Let $x=(\omega_{\xi,\eta}\otimes\iota)(W_{\Omega})$ with $\omega_{\xi,\eta}\in\mathcal{I}_{\Omega}$ . Let $L=\Nset\times\Nset$ with the product order and consider the net $x_{(n,m)}=(\omega_{\xi,e_{n}f_{m}\eta}\otimes\iota)(W_{\Omega})$.
Because $e_{n}f_{m}\eta\rightarrow\eta$, we have $x_{(n,m)}\rightarrow x$ in norm. Note that $e_{n}f_{m}\eta\in\mathcal{D}(A^{-\frac{1}{2}})\cap\mathcal{D}(B^{\frac{1}{2}})$. Moreover, using the same techniques, one can show that $\omega_{\xi,e_{n}f_{m}\eta}\in\mathcal{I}_{\Omega}$. Thus, $\omega_{\xi,e_{n}f_{m}\eta}\in C_{2}$ and we have, for all $x\in M$ such that $xA^{\frac{1}{2}}$ is bounded and $\overline{xA^{\frac{1}{2}}}\in\mathcal{N}_{\varphi}$,
\begin{eqnarray*}
\langle\hat{\Lambda}_{\Omega}(x_{(n,m)}),\Lambda_{\Omega}(x)\rangle
&=&\langle\xi_{\Omega}(\omega_{\xi,e_{n}f_{m}\eta}),\Lambda_{\Omega}(x)\rangle
=\langle x^{*}\xi,e_{n}f_{m}\eta\rangle\\
&=&\langle (xe_{n}f_{m})^{*}\xi,\eta\rangle
=\langle\xi_{\Omega}(\omega_{\xi,\eta}),\Lambda_{\Omega}(xe_{n}f_{m})\rangle,
\end{eqnarray*}
because $xe_{n}f_{m}A^{\frac{1}{2}}$ is bounded and $\overline{xe_{n}f_{m}A^{\frac{1}{2}}}=\overline{xA^{\frac{1}{2}}}e_{n}f_{m}\in
\mathcal{N}_{\varphi}$, so
\begin{eqnarray*}
\langle\hat{\Lambda}_{\Omega}(x_{(n,m)}),\Lambda_{\Omega}(x)\rangle
&=&\langle\xi_{\Omega}(\omega_{\xi,\eta}),J\sigma_{-\frac{i}{2}}(e_{n})\sigma_{-\frac{i}{2}}(f_{m})J\Lambda(\overline{xA^{\frac{1}{2}}})\rangle\\
&=&\langle J\sigma_{\frac{i}{2}}(e_{n})\sigma_{\frac{i}{2}}(f_{m})J\xi_{\Omega}(\omega_{\xi,\eta}),\Lambda_{\Omega}(x)\rangle
\end{eqnarray*}
Thus, $\hat{\Lambda}_{\Omega}(x_{(n,m)})= J\sigma_{\frac{i}{2}}(e_{n})\sigma_{\frac{i}{2}}(f_{m}) J\hat{\Lambda}_{\Omega}(x)\rightarrow\hat{\Lambda}_{\Omega}(x)$.
\end{proof}

Next proposition describes the image by $\Lambda_{1}$ of typical elements from $\Lambda_{2}$. Let us define the unitaries
$$U=(\lambda_{L}^{\hat{G}^{2}}\otimes\alpha)(\Psi^{\circ})^{*}\quad\text{and}\quad V=(\lambda_{R}^{\hat{G}^{2}}\otimes\alpha^{'})(\sigma\Psi^{*}).$$

\begin{proposition}\label{PropLambda1}
Let $x\in C_{1}^{0}$ and $\omega\in M_{*}$ be such that $\omega u\in\mathcal{I}$, then $\Lambda_{P}\otimes\xi(\omega u)\in\mathcal{D}(T^{-\frac{1}{2}})$, $ (x\otimes 1)(\omega\otimes\iota)(\tilde{W}_{\Omega})\in\mathcal{N}_{\varphi_{1}}$ and
$$\Lambda_{1}\left((x\otimes 1)(\omega\otimes\iota)(\tilde{W}_{\Omega})\right)=UVT^{-\frac{1}{2}}
\Lambda_{P}(x)\otimes\xi(\omega u).$$
\end{proposition}

First, we need some preliminary results.

\begin{lemma}\label{J1}
Let $J_{1}$ be the modular conjugation associated with $\varphi_{1}$. Then, for all $F\in L^{\infty}(G)$,
$$(\lambda_{L}^{\hat{G}^{2}}\otimes\alpha)(\Delta_{G}(F))=J_{1}\lambda_{L}(F)^{*}J_{1}.$$
\end{lemma}

\begin{proof}
Using Lemma \ref{LemRN} $(\ref{LemRN1})$, we see that $((\gamma_{1},\gamma_{2})\mapsto L_{\gamma_{1}}R_{\gamma_{2}})$ is the standard implementation of the action $\alpha$ on $H=H_{\hat{\varphi}}$, so the operator $J_{1}$ is given by $J_{1}\xi(\gamma_{1},\gamma_{2})=L_{-\gamma_{1}}R_{-\gamma_{2}}\hat{J}\xi(-\gamma_{1},-\gamma_{2})$, for $\xi\in L^{2}(\hat{G}^{2},H)$. It is now easy to check the needed equality for $F=u_{\gamma}$ with $\gamma\in\hat{G}$. Because $(\lambda_{L}^{\hat{G}^{2}}\otimes\alpha)\circ\Delta_{G}$ and $J_{1}\lambda_{L}(.)^{*} J_{1}$ are $*$-homomorphisms, this concludes the proof.
\end{proof}

Define one parameter groups of automorphisms of $L^{\infty}(G): \gamma_{t}(F)(x)=F(x-\gamma_{t})$ and of $M^{'}: \sigma_{t}^{'}(x)=J\sigma_{t}(JxJ)J$. Note that $\sigma_{t}^{'}\circ\alpha^{'}=\alpha^{'}\circ\gamma_{t}$. By analytic continuation, $\alpha^{'}(F)\in\mathcal{D}(\sigma_{z}^{'})$ and $\sigma_{z}^{'}(\alpha^{'}(F))=\alpha^{'}(\gamma_{z}(F))\ \forall z\in\Cset, \ F\in\mathcal{D}(\gamma_{z})$.

\begin{lemma}\label{Lemcalcul}
Let $F\in L^{\infty}(G^{2})$, $x\in\mathcal{N}_{\varphi_{P}}$, and $\omega\in\mathcal{I}$. If $F\in\mathcal{D}(\gamma_{-\frac{i}{2}}\otimes\iota)$, then $(\lambda_{R}^{\hat{G}^{2}}\otimes\alpha^{'})(\sigma F)\in\mathcal{D}(\iota\otimes\sigma_{-\frac{i}{2}}^{'}),\ (x\otimes 1) (\omega\otimes\iota)\left((\alpha\otimes\lambda_{R})(F)\tilde{W}(\alpha\otimes\lambda_{L})(F)\right)
\in\mathcal{N}_{\varphi_{1}}$, and
$$\Lambda_{1}\left((x\otimes 1)(\omega\otimes\iota)\left((\alpha\otimes\lambda_{R})(F)\tilde{W}(\alpha\otimes\lambda_{L})(F)\right)\right)$$
$$=(\lambda_{L}^{\hat{G}^{2}}\otimes\alpha)(F^{\circ})(\iota\otimes\sigma_{-\frac{i}{2}}^{'})\left((\lambda_{R}^{\hat{G}^{2}}\otimes\alpha^{'})(\sigma F)\right)\Lambda_{P}(x)\otimes\xi(\omega).$$
\end{lemma}

\begin{proof}
Because $\mathcal{D}(\gamma_{-\frac{i}{2}})\odot L^{\infty}(G)$ is a $\sigma$-strong* core for $\gamma_{-\frac{i}{2}}\otimes\iota$ and $\Lambda_{1}$ is $\sigma$-weak-weak closed, we can take $F\in\mathcal{D}(\gamma_{-\frac{i}{2}})\odot L^{\infty}(G)$. By linearity, we can take $F=F_{1}\otimes F_{2}$ with $F_{1}\in\mathcal{D}(\gamma_{-\frac{i}{2}})$ and $F_{2}\in L^{\infty}(G)$. If $x\in\mathcal{N}_{\varphi_{P}}$ and $\omega\in\mathcal{I}$, then
\begin{eqnarray*}
&&(x\otimes 1)(\omega\otimes\iota)\left((\alpha\otimes\lambda_{R})(F)\tilde{W}(\alpha\otimes\lambda_{L})(F)\right)\\
&&=\lambda_{R}(F_{2})(x\otimes 1)(\alpha(F_{1}).\omega.\alpha(F_{1})\otimes\iota)(\tilde{W})\lambda_{L}(F_{2}).
\end{eqnarray*}
Because $F_{1}\in\mathcal{D}(\gamma_{-\frac{i}{2}})$, we have $\alpha(F_{1})\in\mathcal{D}(\sigma_{-\frac{i}{2}})$. Lemma \ref{awb} and the definition of $\Lambda_{1}$ imply $(x\otimes 1)(\alpha(F_{1})\cdot\omega\cdot\alpha(F_{1})\otimes\iota)(\tilde{W})\in\mathcal{N}_{\varphi_{1}}$ and
$$
\begin{array}{l}
\Lambda_{1}\left((x\otimes 1)(\alpha(F_{1})\cdot\omega\cdot\alpha(F_{1})\otimes\iota)(\tilde{W})\right)\\
=(1\otimes\alpha(F_{1}))(1\otimes J\sigma_{-\frac{i}{2}}(\alpha(F_{1})^{*}J)(\Lambda_{P}(x)\otimes\xi(\omega))\\
=(1\otimes\alpha(F_{1}))(1\otimes \alpha^{'}(\gamma_{-\frac{i}{2}}(F_{1})))(\Lambda_{P}(x)\otimes\xi(\omega)).
\end{array}
$$
Moreover, $(\ref{EqMod1})$ gives $\lambda_{L}(F_{2})\in N^{\varphi_{1}}$, so
$$
\begin{array}{l}
\lambda_{R}(F_{2})(x\otimes 1)(\alpha(F_{1}).\omega.\alpha(F_{1})\otimes\iota)(\tilde{W})\lambda_{L}(F_{2})
\in\mathcal{N}_{\varphi_{1}}\quad\text{and,}\\
\Lambda_{1}\left(\lambda_{R}(F_{2})(x\otimes 1)(\alpha(F_{1}).\omega.\alpha(F_{1})\otimes\iota)(\tilde{W})\lambda_{L}(F_{2})\right)\\
=J_{1}\lambda_{L}(F_{2})^{*}J_{1}\lambda_{R}(F_{2})(1\otimes\alpha(F_{1}))(1\otimes \alpha^{'}(\gamma_{-\frac{i}{2}}(F_{1})))(\Lambda_{P}(x)\otimes\xi(\omega))\\
=J_{1}\lambda_{L}(F_{2})^{*}J_{1}(1\otimes\alpha(F_{1}))
\left((\lambda_{R}^{\hat{G}^{2}}(F_{2})\otimes1)(1\otimes \alpha^{'}(\gamma_{-\frac{i}{2}}(F_{1})))\right)
(\Lambda_{P}(x)\otimes\xi(\omega))\\
\quad(\text{because}\,\,\,\lambda_{R}(F_{2})=\lambda_{R}^{\hat{G}^{2}}(F_{2})\otimes1\,\,\,\text{commute with}\,\,\,1\otimes\alpha(F_{1}))\\
=(\lambda_{L}^{\hat{G}^{2}}\otimes\alpha)(\Delta_{G}(F_{2})1\otimes F_{1})
(\lambda_{R}^{\hat{G}^{2}}\otimes\alpha^{'})((F_{2}\otimes1)(1\otimes\gamma_{-\frac{i}{2}}(F_{1})))
(\Lambda_{P}(x)\otimes\xi(\omega)),
\end{array}
$$
where we used Lemma \ref{J1} in the last equality. Note that
$$(\Delta_{G}(F_{2})1\otimes F_{1})(g,h)=F_{2}(g+h)F_{1}(h)=F(h,g+h)=F^{\circ}(g,h),$$
and because
\begin{eqnarray*}
(\lambda_{R}^{\hat{G}^{2}}\otimes\alpha^{'})((F_{2}\otimes1)(1\otimes\gamma_{-\frac{i}{2}}(F_{1})))
&=&(\lambda_{R}^{\hat{G}^{2}}\otimes\alpha^{'})((\iota\otimes\gamma_{-\frac{i}{2}})(\sigma F))\\
&=&(\iota\otimes\sigma^{'}_{-\frac{i}{2}})((\lambda_{R}^{\hat{G}^{2}}\otimes\alpha^{'})(\sigma F)),
\end{eqnarray*}
we conclude the proof.
\end{proof}

\begin{lemma}\label{LemVext}
The operator $(\iota\otimes\sigma^{'}_{-\frac{i}{2}})(V)$ is normal, affilated with $\mathcal{L}(\widehat{G}^{2})\otimes M^{'}$, and its polar decomposition is
$(\iota\otimes\sigma^{'}_{-\frac{i}{2}})(V)=V(T_{1}^{-\frac{1}{2}}\otimes 1)=VT^{-\frac{1}{2}}$.
\end{lemma}
\begin{proof}
We have $(\iota\otimes\gamma_{t})(\sigma\Psi^{*})(g,h)=\Psi^{*}(h,g)\Psi^{*}(-\gamma_{t},g)$,
so
$$(\iota\otimes\sigma^{'}_{t})(V)=
(\lambda_{R}^{\hat{G}^{2}}\otimes\alpha^{'})((\iota\otimes\gamma_{t})(\sigma\Psi^{*}))=V(T_{1}^{-it}
\otimes 1).$$
We conclude the proof by applying Proposition \ref{Propext}.
\end{proof}

\textit{Proof of Proposition \ref{PropLambda1}}. Let $x\in\mathcal{C}^{0}_{1}$ and $\omega\in M_{*}$ such that $\omega\cdot u\in\mathcal{I}$. By Lemma \ref{LemVext}, $\Lambda_{P}(x)\otimes\xi(\omega \cdot u)\in\mathcal{D}((\iota\otimes\sigma^{'}_{-\frac{i}{2}})(V))$ and
$$(\iota\otimes\sigma^{'}_{-\frac{i}{2}})(V)\Lambda_{P}(x)\otimes\xi(\omega\cdot u)
=VT^{-\frac{1}{2}}\Lambda_{P}(x)\otimes\xi(\omega\cdot u).$$
By Lemma \ref{group}, $V(n)\rightarrow V$ $\sigma$-strongly* and
$$(\iota\otimes\sigma^{'}_{-\frac{i}{2}})(V(n))\Lambda_{P}(x)\otimes\xi(\omega\cdot u)
\rightarrow VT^{-\frac{1}{2}}\Lambda_{P}(x)\otimes\xi(\omega\cdot u),$$
where
$$V(n)=\sqrt{\frac{n}{\pi}}\int e^{-nt^{2}}(\iota\otimes\sigma^{'}_{t})(V)dt
=(\lambda_{R}^{\hat{G}^{2}}\otimes\alpha^{'})(\sigma\Psi^{*}(n)),$$
with $\Psi^{*}(n)=\sqrt{\frac{n}{\pi}}\int e^{-nt^{2}}(\gamma_{t}\otimes\iota)(\Psi^{*})dt$. So $\Psi^{*}(n)$ is analytic w.r.t. $(t\mapsto\gamma_{t}\otimes\iota)$ and $\Psi^{*}(n)\rightarrow \Psi^{*}$ $\sigma$-strongly*. Now we can apply Lemma \ref{Lemcalcul} to $\Psi^{*}(n)$ and $\omega\cdot u:
(x\otimes 1)(\omega\cdot u\otimes\iota)\left((\alpha\otimes\lambda_{R}) (\Psi^{*}(n))\tilde{W}(\alpha\otimes\lambda_{L})(\Psi^{*}(n))\right)\in\mathcal{N}_{\varphi_{1}}$ and
$$
\Lambda_{1}\left((x\otimes 1)(\omega \cdot u\otimes\iota)\left((\alpha\otimes\lambda_{R}) (\Psi^{*}(n))\tilde{W}(\alpha\otimes\lambda_{L})(\Psi^{*}(n))\right)\right)
$$
$$
=(\lambda_{L}^{\hat{G}^{2}}\otimes\alpha)(\Psi^{*}(n)^{\circ})(\iota\otimes\sigma^{'}_{-\frac{i}{2}})(V(n))
\Lambda_{P}(x)\otimes\xi(\omega\cdot u).$$
Note that
$$(\alpha\otimes\lambda_{R})(\Psi^{*}(n))\tilde{W}(\alpha\otimes\lambda_{L})(\Psi^{*}(n))
\rightarrow(\alpha\otimes\lambda_{R})(\Psi^{*})\tilde{W}(\alpha\otimes\lambda_{L})(\Psi^{*})\quad\sigma -\text{weakly, so},
$$
$$
(x\otimes 1)(\omega \cdot u\otimes\iota)\left((\alpha\otimes\lambda_{R}) (\Psi^{*}(n))\tilde{W}(\alpha\otimes\lambda_{L})(\Psi^{*}(n))\right)\rightarrow (x\otimes 1)(\omega\otimes\iota)(\tilde{W}_{\Omega})\quad
\sigma -\text{weakly},$$
and
$$
(\lambda_{L}^{\hat{G}^{2}}\otimes\alpha)(\Psi^{*}(n)^{\circ})(\iota\otimes
\sigma^{'}_{-\frac{i}{2}})(V(n))\Lambda_{P}(x)\otimes\xi(\omega .u)\rightarrow
UVT^{-\frac{1}{2}}\Lambda_{P}(x)\otimes\xi(\omega .u)\quad\text{weakly}.
$$
Because $\Lambda_{1}$ is $\sigma$-weak-weak closed, this concludes the proof.

\begin{lemma}\label{LemT12}
Let $\eta\in\mathcal{D}(B^{\frac{1}{2}})$, $\xi\in H$ and $x\in\mathcal{C}_{1}$. Then
$$
(x\otimes 1)(\omega_{\xi,\eta}\otimes\iota)(\tilde{W}_{\Omega})T^{\frac{1}{2}}\,\,\text{is bounded and its closure is}\,\,\, \overline{T^{\frac{1}{2}}(x\otimes 1)}(\omega_{\xi,B^{\frac{1}{2}}\eta}\otimes\iota)(\tilde{W}_{\Omega}).
$$
\end{lemma}

\begin{proof}
Using $(\ref{WOmegaT})$, for all $t\in\Rset$, we have $\tilde{W}_{\Omega}(1\otimes T^{it})\tilde{W}_{\Omega}^{*}=B^{it}\otimes T^{it}$, so $\tilde{W}_{\Omega}(1\otimes T^{\frac{1}{2}})\tilde{W}_{\Omega}^{*}=B^{\frac{1}{2}}\otimes T^{\frac{1}{2}}$.
Let $\eta\in\mathcal{D}(B^{\frac{1}{2}})$, $\xi\in H$, $x\in\mathcal{C}_{1},\ f\in\mathcal{D}(T^{\frac{1}{2}})$, and $l\in L^{2}(\hat{G}^{2},H)$, then
\begin{eqnarray*}
\langle (x\otimes 1)(\omega_{\xi,\eta}\otimes\iota)(\tilde{W}_{\Omega})T^{\frac{1}{2}}f,l\rangle
&=&\langle\tilde{W}_{\Omega}\xi\otimes T^{\frac{1}{2}}f,\eta\otimes(x\otimes1)^{*}l\rangle\\
&=&\langle(B^{\frac{1}{2}}\otimes T^{\frac{1}{2}})\tilde{W}_{\Omega}\xi\otimes f,\eta\otimes(x\otimes1)^{*}l\rangle\\
&=&\langle (1\otimes\overline{(x\otimes 1)T^{\frac{1}{2}}})\tilde{W}_{\Omega}\xi\otimes f,B^{\frac{1}{2}}\eta\otimes l \rangle\\
&=&\langle \overline{T^{\frac{1}{2}}(x\otimes 1)}(\omega_{\xi,B^{\frac{1}{2}}\eta}\otimes\iota)(\tilde{W}_{\Omega})f,l \rangle.\\
\end{eqnarray*}
Thus, we have $(x\otimes 1)(\omega_{\xi,\eta}\otimes\iota)(\tilde{W}_{\Omega})T^{\frac{1}{2}}
\subset\overline{T^{\frac{1}{2}}(x\otimes 1)}(\omega_{\xi,B^{\frac{1}{2}}\eta}\otimes\iota)(\tilde{W}_{\Omega})$. Because $\mathcal{D}(T^{\frac{1}{2}})$ is dense, this concludes the proof.
\end{proof}

\begin{proposition}\label{lambdamu2}
Let $x\in\mathcal{C}^{0}_{1}$ and $\omega_{\xi,\eta}\in C_{2}$. Then
$$
(x\otimes 1)(\omega_{\xi,\eta}\otimes\iota)(\tilde{W}_{\Omega}) \in\mathcal{N}_{\mu}\cap\mathcal{N}_{\varphi_{2}}\quad\text{and}
$$
$$\Lambda_{\mu}\left((x\otimes 1)(\omega_{\xi,\eta}\otimes\iota)(\tilde{W}_{\Omega})\right)
=\lambda^{\frac{i}{4}}UV(1\otimes Ju^{*}J)\Lambda_{2}\left((x\otimes 1)(\omega_{\xi,\eta}\otimes\iota)(\tilde{W}_{\Omega})\right).
$$
\end{proposition}
\begin{proof}
Let $x\in\mathcal{C}^{0}_{1}$ and $\omega_{\xi,\eta}\in C_{2}$. By Lemma \ref{LemT12}, $(x\otimes 1)(\omega_{\xi,\eta}\otimes\iota)(\tilde{W}_{\Omega})T^{\frac{1}{2}}$ is bounded and its closure is $\overline{T^{\frac{1}{2}}(x\otimes 1)}(\omega_{\xi,B^{\frac{1}{2}}\eta}\otimes\iota)(\tilde{W}_{\Omega})$. Moreover, by Lemma \ref{LemCore}, $\omega_{\xi,B^{\frac{1}{2}}\eta}\cdot u=\omega_{\xi,u^{*}B^{\frac{1}{2}}\eta}\in\mathcal{I}$, so we can apply Proposition \ref{PropLambda1}, and we find that $(x\otimes 1)(\omega_{\xi,B^{\frac{1}{2}}\eta}\otimes\iota)(\tilde{W}_{\Omega})\in\mathcal{N}_{\varphi_{1}}$ and
$$
\Lambda_{1}\left((x\otimes 1)(\omega_{\xi,B^{\frac{1}{2}}\eta}\otimes\iota)(\tilde{W}_{\Omega})\right)
=UVT^{-\frac{1}{2}}\Lambda_{P}(x)\otimes\xi(\omega_{\xi,u^{*}B^{\frac{1}{2}}\eta}).
$$
Finally, using Proposition \ref{PropTx}, and because $T^{\frac{1}{2}}$ commutes with $UV$, we find that $\overline{T^{\frac{1}{2}}(x\otimes 1)}(\omega_{\xi,B^{\frac{1}{2}}\eta}\otimes\iota) (\tilde{W}_{\Omega})\in\mathcal{N}_{\varphi_{1}}$ and
$$
\Lambda_{1}\left(\overline{T^{\frac{1}{2}}(x\otimes 1)}(\omega_{\xi,B^{\frac{1}{2}}\eta}\otimes\iota)(\tilde{W}_{\Omega})\right)
=UV\Lambda_{P}(x)\otimes\xi(\omega_{\xi,u^{*}B^{\frac{1}{2}}\eta}).
$$
By Lemma \ref{LemCore}, $\xi(\omega_{\xi,u^{*}B^{\frac{1}{2}}\eta})= \lambda^{\frac{i}{4}}Ju^{*}J\xi_{\Omega}(\omega_{\xi,\eta})$, so
$$
\begin{array}{l}
\Lambda_{1}\left(\overline{T^{\frac{1}{2}}(x\otimes 1)}(\omega_{\xi,B^{\frac{1}{2}}\eta}\otimes\iota)(\tilde{W}_{\Omega})\right)
=\lambda^{\frac{i}{4}}UV(1\otimes Ju^{*}J)\Lambda_{P}(x)\otimes\xi_{\Omega}(\omega_{\xi,\eta})\\
=\lambda^{\frac{i}{4}}UV(1\otimes Ju^{*}J)\Lambda_{2}\left((x\otimes 1)(\omega_{\xi,\eta}\otimes\iota)(\tilde{W}_{\Omega})\right).
\end{array}
$$

\end{proof}
\textit{Proof of Theorem \ref{dualweight}.}
Let $\mathcal{D}$ be the $\sigma$-weak-weak core for $\Lambda_{2}$ introduced before Lemma \ref{LemCore}. By Proposition \ref{lambdamu2}, $\mathcal{D}\subset\mathcal{N}_{\mu}\cap\mathcal{N}_{\varphi_{2}}$ and there is a unitary $Z$ such that $\Lambda_{2}(x)=Z\Lambda_{\mu}(x)$, for all $x\in\mathcal{D}$. By Proposition \ref{Equality}, $\varphi_{2}=\tilde{\mu}_{\Omega}$, so $\hat{\varphi}_{\Omega}=\mu_{\Omega}$.

\section{Examples}\label{SectionExample}

\subsection{Twisting of the $az+b$ group}\label{Sectionaz+bclassique}

Our aim is to prove Theorem \ref{Theoremaz+bclassique}. According to Section \ref{SectionResult},
if $H$ is a closed abelian subgroup of a \lc{} group $G$, then $H<(\mathcal{L}(G),\hat{\Delta}_{G})$
is an abelian stable co-subgroup. The morphism $\alpha: L^{\infty}(\widehat{H})\to\mathcal{L}(G)$ is given by $\alpha(u_{h})=\lambda_{G}(h)$, and the
morphism $(t\mapsto\gamma_{t}):\R\to\widehat{H}$ by $\langle\gamma_{t},h\rangle=\delta_{G}^{-it}(h)$.
Let $G=\mathbb{C}^{*}\ltimes\mathbb{C}$
and $K\subset G$ be the subgroup $K=\left\{(z,0),\,z\in\mathbb{C}^{*}\right\}$.
The modular function of $G$ is $\delta_{G}(z,w)=|z|^{-1}$, for all $z\in\C^*,\ \omega\in\C$, and $\langle\gamma_{t},z\rangle= |z|^{it}$, for all $z\in\C^{*},t\in\R$. Let us identify $\widehat{\mathbb{C}^{*}}$ with $\mathbb{Z}\times\mathbb{R}_{+}^{*}$:
$$
\mathbb{Z}\times\mathbb{R}_{+}^{*}\rightarrow\widehat{\mathbb{C}^{*}},\quad (n,\rho)\mapsto\gamma_{n,\rho}=(re^{i\theta}\mapsto e^{i\ln{r}\ln{\rho}}e^{in\theta}).
$$
Then $\gamma_{t}=(0,e^{t})\in \mathbb{Z}\times\mathbb{R}_{+}^{*}$. For any $x\in\R$, there is a bicharacter on $\mathbb{Z}\times\mathbb{R}_{+}^{*}:
\Psi_{x}((n,\rho),(k,r))=e^{ix(k\ln{\rho}-n\ln{r})}$. Let $(M_{x},\Delta_{x})$ be the l.c. quantum group obtained by twisting. Then
$\Psi_{x}((n,\rho),\gamma_{t}^{-1})=e^{ixtn}=u_{e^{ixt}}((n,\rho))$, and we
get the operator $A_{x}$ deforming the Plancherel weight $\varphi$:
$$
A_{x}^{it}=\alpha(u_{e^{2ixt}})=\lambda^{G}_{(e^{itx},0)}.
$$
Since $\Psi_{x}(\gamma_{t},\gamma_{s})=1$, for all $s,t\in\R$, the twisted left-invariant weight $\varphi_{x}$ satisfies $[D\varphi_{x}\,:\,D\varphi]_{t}=A_{x}^{it}=\lambda^{G}_{(e^{itx},0)}$.
The modular element of the twisted quantum group is $$
\delta_{x}^{it}=\alpha(\Psi_{x}(\cdot,\gamma_{t})\Psi_{x}(-\gamma_{t},\cdot))=\lambda^{G}_{(e^{-2itx},0)},
$$
so $\delta_{x}$ is not affiliated with the center of $\mathcal{L}(G)$, and the twisted quantum group is not a Kac algebra. Let us look if $(M_{x},\Delta_{x})$ is isomorphic for different values of $x$. Since $\Psi_x$ is antisymmetric, $\Psi_{-x}=\Psi_{x}^{*}$, and $\Delta$ is cocommutative, we have $\Delta_{-x}=\sigma\Delta_{x}$. Thus, $(M_{-x},\Delta_{-x})\simeq (M_{x},\Delta_{x})^{\text{op}}$. Moreover, using the Fourier transformation in the first variable, one has immediately $\text{Sp}(\delta_{x})=q_{x}^{\Z}\cup\{0\},$ where $q_{x}=e^{-2x}$. Thus, if $x\neq y,\ x>0,y>0$, one
has $q_{x}^{\Z}\neq q_{y}^{\Z}$ and, consequently, $(M_{x},\Delta_{x})$ and $(M_{y},\Delta_{y})$ are
not isomorphic.

In order to finish the proof of Theorem \ref{Theoremaz+bclassique}, we must compute
the dual \lc{} quantum group. The action of $K^{2}$ on $L^{\infty}(G)$ can be lifted
to its Lie algebra $\mathbb{C}^{2}$ which does not change the result of
deformation (see \cite{Kas}) but simplifies
calculations. The group $\mathbb{C}$ is self-dual with the duality
$(z_{1},z_{2})\mapsto \exp\left(i\text{Im}(z_{1}z_{2})\right)$.
Let $x\in\mathbb{R}$. The lifted bicharacter on $\mathbb{C}$ is $\Psi_{x}(z_{1},z_{2})=\exp\left(ix\text{Im}(z_{1}\overline{z}_{2})\right)$. The action
$\rho$ of $\mathbb{C}^{2}$ on $L^{\infty}(G)$ is
\begin{equation}\label{actionclassique}
\rho_{z_{1},z_{2}}(f)(w_{1},w_{2})=f(e^{z_{2}-z_{1}}w_{1},e^{-z_{1}}w_{2}).
\end{equation}
Let $N=\C^{2}\ltimes L^{\infty}(G)$ and $\theta$ be the dual action of $\C^{2}$ on $N$.
One has, for all $z,w\in\C$, $\Psi_{x}(w,z)=u_{x\overline{z}}(w)$. So, the twisted dual
action is
\begin{equation}\label{actiondualedeformee}
\theta_{z_{1},z_{2}}^{\Psi_{x}}=
\lambda_{-x\overline{z}_{1},x\overline{z}_{2}}\theta_{z_{1},z_{2}}(\cdot)
\lambda^{*}_{-x\overline{z}_{1},x\overline{z}_{2}}.
\end{equation}
Let $\widehat{M}_{x}$ be the fixed point algebra. We will construct two
operators affiliated with $\widehat{M}_{x}$ which generate $\widehat{M}_{x}$.
Let $a$ and $b$ be the coordinate functions on $G$, and
$\alpha=\pi(a)$, $\beta=\pi(b)$. Then $\alpha$ and $\beta$ are
normal operators affiliated with $N$, and $(\ref{actionclassique})$ gives
\begin{equation}\label{Eqlambdagenerators}
\lambda_{z_{1},z_{2}}\alpha\lambda_{z_{1},z_{2}}^{*}=e^{z_{2}-z_{1}}\alpha,\quad
\lambda_{z_{1},z_{2}}\beta\lambda_{z_{1},z_{2}}^{*}=e^{-z_{1}}\beta.
\end{equation}
Now, using $(\ref{actiondualedeformee})$ and $(\ref{Eqlambdagenerators})$, we find
\begin{equation}\label{twisteddualgen}
\theta^{\Psi_{x}}_{z_{1},z_{2}}(\alpha)=e^{x(\overline{z}_{1}+\overline{z}_{2})}\alpha\, ,\quad
\theta^{\Psi_{x}}_{z_{1},z_{2}}(\beta)=e^{x\overline{z}_{1}}\beta.
\end{equation}
Let $T_{l}$ and $T_{r}$ be the infinitesimal generators of the
left and right translations, so $T_{l}$ and $T_{r}$ are
affiliated with $N$ and $\lambda_{z_{1},
z_{2}}=\exp\left(i\text{Im}(z_{1}T_{l})\right)\exp\left(i\text{Im}(z_{2}T_{r})\right).$
Then $\lambda(f)=f(T_{l},T_{r})$, for all $f\in L^{\infty}(\mathbb{C}^{2})$.

\begin{lemma}\label{Rbeta}
Let $L=e^{xT_{l}^{*}}$ and $R=e^{xT_{r}^{*}}$, then
\begin{itemize}
\item $(\beta,L)$ is a $e^{x}$-commuting pair.
\item $(\beta,R)$ is a $1$-commuting pair.
\item $(\alpha,R)$ is a $e^{-x}$-commuting pair.
\item $(\alpha,L)$ is a $e^{x}$-commuting pair.
\end{itemize}
\end{lemma}
\begin{proof}
Note that $\text{Ph}(L)=e^{-ix\text{Im}T_{l}}=\lambda_{-x,0}$ and $|L|^{is}=e^{isx\text{Re}T_{l}}=\lambda_{isx,0}$, so $(\ref{Eqlambdagenerators})$ gives $|L|^{is}\beta|L|^{-is}=e^{-isx}\beta$ and $\text{Ph}(L)\beta\text{Ph}(L)^{*}=e^{x}\beta$ which means that $(\beta,L)$ is a $e^{x}$-commuting pair. The proof of the other assertions is similar.
\end{proof}

Define $U=\lambda(\Psi_{x})$ and $\hat{\alpha}=U^{*}\alpha
U$, $\hat{v}=\text{Ph}(L)\text{Ph}(\beta)$ and $\hat{B}=|L||\beta|$. Then $\hat{\alpha}$ is normal, $\hat{B}$ is positive self adjoint, both
affiliated with $N$, and $\hat{v}\in N $ is unitary.

\begin{proposition}\label{generators2}
$\hat{\alpha}$ and $\hat{B}$ are affiliated with $\widehat{M}_{x}$ and $\hat{v}\in\widehat{M}_{x}$. Moreover,
$$\left\{f(\hat{\alpha})g(\hat{B})h(\hat{v}),\,\,f\in L^{\infty}(\mathbb{C}),g\in L^{\infty}(\mathbb{R}^{*}_{+}),h\in L^{\infty}(\mathbb{S}^{1})\right\}^{''}=\widehat{M}_{x}.$$
\end{proposition}
\begin{proof}
We have
\begin{eqnarray*}
\theta^{\Psi_{x}}_{z_{1},z_{2}}(U) &=& \lambda\left(\Psi_{x}(.-z_{1},.-z_{2})\right)\\
&=& U e^{ix\text{Im}(-\overline{z}_{2}T_{l})}e^{ix\text{Im}(\overline{z}_{1}T_{r})}\Psi_{x}(z_{1},z_{2})\\
&=&U \lambda_{-x\overline{z}_{2},x\overline{z}_{1}}\Psi_{x}(z_{1},z_{2}).
\end{eqnarray*}
This implies, using (\ref{twisteddualgen}) and $(\ref{Eqlambdagenerators})$:
$$
\theta^{\Psi_{x}}_{z_{1},z_{2}}(\hat{\alpha})
=e^{x(\bar{z}_{1}+\bar{z}_{2})}U^{*} \lambda_{x\bar{z}_{2},-x\bar{z}_{1}}\alpha\lambda_{x\bar{z}_{2},-x\bar{z}_{1}}^{*}U=\hat{\alpha}.
$$
Also,
$$
\theta^{\Psi_{x}}_{z_{1},z_{2}}(\hat{B})=
e^{x\text{Re}(T_{l}-\bar{z}_{1})}e^{x\text{Re}(\bar{z}_{1})}|\beta|=\hat{B}
$$
and
$$
\theta^{\Psi_{x}}_{z_{1},z_{2}}(\hat{v})=
e^{ix\text{Im}(T_{l}^{*}-\bar{z}_{1})}e^{ix\text{Im}(\bar{z}_{1})}\text{Ph}(\beta)=\hat{v}.
$$
Thus, $\hat{\alpha}$ and $\hat{B}$ are affiliated with $\widehat{M}_{x}$ and $\hat{v}\in\widehat{M}_{x}$. Let
$$
\mathcal{W}=\left\{zf(\hat{\alpha})g(\hat{B})h(\hat{v})y,\,\,f\in L^{\infty}(\mathbb{C}),
g\in L^{\infty}(\mathbb{R}^{*}_{+}),h\in
L^{\infty}(\mathbb{S}^{1}),\,\,
z,y\in\lambda(L^{\infty}(\mathbb{C}^{2}))\right\}^{''}.
$$
By Lemma \ref{FixSub}, it suffices to show that $\mathcal{W}=N$. Note that
$$
\left\{zf(\hat{\alpha}),\,\,f\in L^{\infty}(\mathbb{C}),z\in\lambda(L^{\infty}(\mathbb{C}^{2}))\right\}^{''}
=\left\{zU^{*}f(\alpha)U,\,\,f\in L^{\infty}(\mathbb{C}),z\in\lambda(L^{\infty}(\mathbb{C}^{2}))\right\}^{''}.
$$
Substituting $z\mapsto zU$, we get
$$
\left\{zf(\hat{\alpha}),\,\,f\in L^{\infty}(\mathbb{C}),z\in\lambda(L^{\infty}(\mathbb{C}^{2}))\right\}^{''}
=\left\{zf(\alpha)U,\,\,f\in L^{\infty}(\mathbb{C}),z\in\lambda(L^{\infty}(\mathbb{C}^{2}))\right\}^{''}.
$$
Observe that
$$\left\{f(\alpha)z\, ,\,\,f\in L^{\infty}(\mathbb{C}),\,z\in \lambda(L^{\infty}(\mathbb{C}^{2}))\right\}^{''}=\left\{zf(\alpha)\, ,\,\,f\in L^{\infty}(\mathbb{C}),\,z\in \lambda(L^{\infty}(\mathbb{C}^{2}))\right\}^{''},$$
so
$$\left\{zf(\hat{\alpha}),\,\,f\in L^{\infty}(\mathbb{C}),z\in\lambda(L^{\infty}(\mathbb{C}^{2}))\right\}^{''}
=\left\{f(\alpha)zU,\,\,f\in L^{\infty}(\mathbb{C}),z\in\lambda(L^{\infty}(\mathbb{C}^{2}))\right\}^{''}.$$
Substituting $z\mapsto zU^{*}$, we get
$$\left\{zf(\hat{\alpha}),\,\,f\in L^{\infty}(\mathbb{C}),z\in\lambda(L^{\infty}(\mathbb{C}^{2}))\right\}^{''}
=\left\{f(\alpha)z,\,\,f\in L^{\infty}(\mathbb{C}),z\in\lambda(L^{\infty}(\mathbb{C}^{2}))\right\}^{''},$$
so
$$\mathcal{W}=\left\{f(\alpha)z g(\hat{B})h(\hat{v})y,\,\,f\in L^{\infty}(\mathbb{C}),
g\in L^{\infty}(\mathbb{R}^{*}_{+}),h\in
L^{\infty}(\mathbb{S}^{1}),\,\,
z,y\in\lambda(L^{\infty}(\mathbb{C}^{2}))\right\}^{''}.$$
Note that
\begin{eqnarray*}
&&\left\{zg(\hat{B}),\,\, g\in L^{\infty}(\mathbb{R}^{*}_{+}),\,\,z\in\lambda(L^{\infty}(\mathbb{C}^{2}))\right\}^{''}
=\left\{z\hat{B}^{is},\,\, s\in\mathbb{R},\,\,z\in\lambda(L^{\infty}(\mathbb{C}^{2}))\right\}^{''}\\
&&\quad =\left\{ze^{ist\text{Re}T_{l}}|\beta|^{is},\,\, s\in\mathbb{R},\,\,x\in\lambda(L^{\infty}(\mathbb{C}^{2}))\right\}^{''}.\\
\end{eqnarray*}
Substitution $z\mapsto ze^{-ist\text{Re}T_{l}}$ gives
\begin{eqnarray*}
&&\left\{zg(\hat{B}),\,\, g\in L^{\infty}(\mathbb{R}^{*}_{+}),\,\,z\in\lambda(L^{\infty}(\mathbb{C}^{2}))\right\}^{''}
=\left\{z|\beta|^{is},\,\, s\in\mathbb{R},\,\,z\in\lambda(L^{\infty}(\mathbb{C}^{2}))\right\}^{''}\\
&&\quad =\left\{zg(|\beta|),\,\, g\in L^{\infty}(\mathbb{R}^{*}_{+}),\,\,z\in\lambda(L^{\infty}(\mathbb{C}^{2}))\right\}^{''}.
\end{eqnarray*}
Also, one can prove that
$$\left\{h(\hat{v})y,\,\, h\in L^{\infty}(\mathbb{S}^{1}),\,\,y\in\lambda(L^{\infty}(\mathbb{C}^{2}))\right\}^{''}=\left\{h(\text{Ph}\beta)y,\,\, h\in L^{\infty}(\mathbb{S}^{1}),\,\,y\in\lambda(L^{\infty}(\mathbb{C}^{2}))\right\}^{''}.$$
Thus,
\begin{eqnarray*}
\mathcal{W}&=&
\left\{f(\alpha)zg(|\beta|)h(\text{Ph}\beta)y,\,\,f\in L^{\infty}(\mathbb{C}),g\in L^{\infty}(\mathbb{R}^{*}_{+}),h\in L^{\infty}(\mathbb{S}^{1}),\,\, z,y\in\lambda(L^{\infty}(\mathbb{C}^{2}))\right\}^{''}\\
&=&\left\{f(\alpha)zg(\beta)y,\,\,f,g\in L^{\infty}(\mathbb{C}),\,\, z,y\in\lambda(L^{\infty}(\mathbb{C}^{2}))\right\}^{''}.
\end{eqnarray*}
Commuting back $f(\alpha)$ and $z$, we have the result.
\end{proof}

Let $\hat{\beta}=\hat{v}\hat{B}$. Then $\hat{\beta}$ is a closed (non normal) operator affiliated with $\widehat{M}_{x}$. Let us give now the commutation relations between $\hat{\alpha}$, $\hat{\beta}$.

\begin{proposition}\label{az+bRel}
$\alpha$ and $T_{l}^{*}+T_{r}^{*}$ strongly commute, and
$\hat{\alpha}=e^{x(T_{l}^{*}+T_{r}^{*})}$, so the polar
decomposition of $\hat{\alpha}$ is
$$
\text{Ph}(\hat{\alpha})=e^{-ix\text{Im}(T_{l}+T_{r})}\text{Ph}(\alpha)=
\text{Ph}(L)\text{Ph}(R)\text{Ph}(\alpha),\,\,
|\hat{\alpha}|=e^{x\text{Re}(T_{l}+T_{r})}|\alpha|=|L||R||\alpha|.
$$
Moreover, the following relations hold with $q=e^{2x}$:
\begin{itemize}
\item $\hat{\beta}\hat{\beta}^{*}=q\hat{\beta}^{*}\hat{\beta}$,
\item $(\hat{\alpha},\hat{\beta})$ is a $\sqrt{q}$-commuting pair.
\end{itemize}
\end{proposition}
\begin{proof}
Since
$$e^{i\text{Im}\left(z(T_{l}^{*}+T_{r}^{*})\right)}\alpha e^{-i\text{Im}\left(z(T_{l}^{*}+T_{r}^{*})\right)}=
\lambda_{-\overline{z},-\overline{z}}\alpha\lambda_{-\overline{z},-\overline{z}}^{*}
=e^{-\overline{z}+\overline{z}}\alpha=\alpha,$$
$T_{l}^{*}+T_{r}^{*}$ and $\alpha$ strongly commute. Moreover, since
$e^{ix\text{Im}T_{l}T_{l}^{*}}=1$,
$$\hat{\alpha}=e^{-ix\text{Im}T_{l}T_{r}^{*}}\alpha e^{ix\text{Im}T_{l}T_{r}^{*}}
=e^{-ix\text{Im}T_{l}(T_{l}+T_{r})^{*}}\alpha
e^{ix\text{Im}T_{l}(T_{l}+T_{r})^{*}}.$$ This
equality, the strong commutativity of $T_{l}^{*}+T_{r}^{*}$ with
$\alpha$, and the equality $e^{-ix\text{Im}T_{l}\omega}\alpha
e^{ix\text{Im}T_{l}\omega}=e^{x\omega}\alpha$ imply $\hat{\alpha}=
e^{x(T_{l}^{*}+T_{r}^{*})}$. The polar decomposition of
$\hat{\alpha}$ follows. All the relations can be checked using Lemma \ref{Rbeta}.
\end{proof}

We shall give now a nice formula for $\hat{\Delta}_{x}$. Let us define the following
(closed non-normal) operator affiliated with $\widehat{M}_{x}\otimes\widehat{M}_{x}$: $\hat{\Delta}_{x}(\hat{\beta})=\hat{\Delta}_{x}({\hat{v}})\hat{\Delta}_{x}({\hat{B}})$.

\begin{proposition}\label{azcom}
$$
\hat{\Delta}_{x}(\hat{\alpha})=\hat{\alpha}\otimes\hat{\alpha}\quad\text{and}\quad
\hat{\Delta}_{x}(\hat{\beta})=\hat{\alpha}\otimes\hat{\beta}\dot{+}\hat{\beta}\otimes 1.
$$
\end{proposition}
\begin{proof}
Proposition \ref{comultiplication} gives
$\hat{\Delta}_{x}=\Upsilon\Gamma(\cdot)\Upsilon^{*}$, where
$\Upsilon=e^{ix\text{Im}T_{r}\otimes T_{l}^{*}},$ and $\Gamma$ is
uniquely characterized by two properties:
\begin{itemize}
\item $\Gamma(T_{l})=T_{l}\otimes 1$,  $\Gamma(T_{r})=1\otimes
T_{r}$; \item $\Gamma$ restricted to $L^{\infty}(G)$ coincides
with the comultiplication $\Delta_G$.
\end{itemize}
With $V=\Upsilon\Gamma(U^{*})$, we have
$\hat{\Delta}_{x}(\hat{\alpha})=V(\alpha\otimes\alpha)V^{*}$, so
it suffices to show that $(U\otimes U)V$ commutes with
$\alpha\otimes\alpha$. Indeed in this case
$$\hat{\Delta}_{x}(\hat{\alpha})=V(\alpha\otimes\alpha)V^{*}
=(U^{*}\otimes U^{*})(U\otimes
U)V(\alpha\otimes\alpha)V^{*}(U^{*}\otimes U^{*})(U\otimes U)
=\hat{\alpha}\otimes\hat{\alpha}.$$ Let us show that $(U\otimes
U)V$ commutes with $\alpha\otimes\alpha$. From
$U=e^{ix\text{Im}T_{l}T_{r}^{*}}$ one has
$$\Gamma(U^{*})=e^{-ix\text{Im}T_{l}\otimes T_{r}^{*}},\quad U\otimes U=e^{ix\text{Im}(T_{l}T_{r}^{*}\otimes 1+1\otimes T_{l}T_{r}^{*})},$$
so $V=e^{-ix\text{Im}(T_{r}^{*}\otimes T_{l}+T_{l}\otimes
T_{r}^{*})}$ and
$$(U\otimes U)V=e^{ix\text{Im}(T_{l}T_{r}^{*}\otimes 1+1\otimes T_{l}T_{r}^{*}-T_{r}^{*}\otimes T_{l}-T_{l}\otimes T_{r}^{*})}.$$
Remark that
$$T_{l}T_{r}^{*}\otimes 1+1\otimes T_{l}T_{r}^{*}-T_{r}^{*}\otimes T_{l}-T_{l}\otimes T_{r}^{*}
=(T_{l}\otimes 1-1\otimes T_{l})(T_{r}^{*}\otimes 1-1\otimes
T_{r}^{*}),$$ so it is enough to show that $T_{l}\otimes
1-1\otimes T_{l}$ and $T_{r}^{*}\otimes 1-1\otimes T_{r}^{*}$
strongly commute with $\alpha\otimes\alpha$, which follows from

\begin{eqnarray*}
e^{i\text{Im}z(T_{r}^{*}\otimes 1-1\otimes T_{r}^{*})}(\alpha\otimes\alpha)e^{-i\text{Im}z(T_{r}^{*}\otimes 1-1\otimes T_{r}^{*})}
&=&
(\lambda_{0,-\overline{z}}\otimes\lambda_{0,\overline{z}})(\alpha\otimes\alpha)(\lambda_{0,-\overline{z}}\otimes\lambda_{0,\overline{z}})^{*}\\
&=&e^{-\overline{z}}e^{\overline{z}}\alpha\otimes\alpha=\alpha\otimes\alpha\, ,\\
e^{i\text{Im}z(T_{l}\otimes 1-1\otimes T_{l})}(\alpha\otimes\alpha)e^{-i\text{Im}z(T_{l}\otimes 1-1\otimes T_{l})}
&=&(\lambda_{z,0}\otimes\lambda_{-z,0})(\alpha\otimes\alpha)(\lambda_{z,0}\otimes\lambda_{-z,0})^{*}\\
&=&e^{-z}e^{z}\alpha\otimes\alpha=\alpha\otimes\alpha.\\
\end{eqnarray*}
By definition of $\hat{\Delta}_{x}$, we have
$$\hat{\Delta}_{x}(\hat{B})=\hat{\Delta}_{x}(e^{x\text{Re}T_{l}}|\beta|)
=(e^{x\text{Re}T_{l}}\otimes 1)\Upsilon|\alpha\otimes\beta
+\beta\otimes 1|\Upsilon^{*},$$
$$\hat{\Delta}_{x}(\hat{v})=\hat{\Delta}_{x}(e^{-ix\text{Im}T_{l}}\text{Ph}(\beta))
=(e^{-ix\text{Im}T_{l}}\otimes
1)\Upsilon\text{Ph}(\alpha\otimes\beta +\beta\otimes
1)\Upsilon^{*}.$$
A direct computation gives
$$\text{Ph}(\alpha\otimes\beta +\beta\otimes 1)(e^{x\text{Re}T_{l}}\otimes 1)=
e^{x}(e^{x\text{Re}T_{l}}\otimes 1)\text{Ph}(\alpha\otimes\beta
+\beta\otimes 1),$$
so
\begin{eqnarray*}
\hat{\Delta}_{x}(\hat{\beta}) & = & e^{x}(e^{xT_{l}^{*}}\otimes 1)\Upsilon(\alpha\otimes\beta +\beta\otimes 1)\Upsilon^{*}\\
 & = & e^{x}(e^{xT_{l}^{*}}\otimes 1)\Upsilon(\alpha\otimes\beta)\Upsilon^{*}\dot{+}e^{x}(e^{xT_{l}^{*}}\otimes 1)
 \Upsilon(\beta\otimes 1)\Upsilon^{*}.\\
\end{eqnarray*}
Thus, it suffices to show that
\begin{eqnarray}\label{eqcom1}
\hat{\alpha}\otimes\hat{\beta}&=&e^{x}(e^{xT_{l}^{*}}\otimes 1)\Upsilon(\alpha\otimes\beta)\Upsilon^{*}\\\label{eqcom2}
\hat{\beta}\otimes 1&=&e^{x}(e^{xT_{l}^{*}}\otimes 1)\Upsilon(\beta\otimes 1)\Upsilon^{*}.\\\notag
\end{eqnarray}
Let us prove $(\ref{eqcom1})$. Let us put $T=e^{x}e^{xT_{l}^{*}}\otimes 1=e^{x}L\otimes 1$ and $S=\Upsilon(\alpha\otimes\beta)\Upsilon^{*}$. We want to show that $\hat{\alpha}\otimes\hat{\beta}=TS$. For all $z\in\mathbb{C}$, we have
$$e^{ix\text{Im}z(T_{r}\otimes 1)}(\alpha\otimes 1)e^{-ix\text{Im}z(T_{r}\otimes 1)}
=(\lambda_{0,xz}\alpha\lambda_{0,-xz}^{*}\otimes
1)=e^{xz}(\alpha\otimes 1),$$ and, using the fact that $\alpha\otimes 1$
and $1\otimes T_{l}^{*}$ strongly commute, we obtain $
\Upsilon(\alpha\otimes 1)\Upsilon^{*}=\alpha\otimes
e^{xT_{l}^{*}}=\alpha\otimes L$. Similarly, $ \Upsilon(1\otimes
\beta)\Upsilon^{*}=R\otimes\beta$. Thus, using Lemma \ref{Rbeta}, we see that the polar decomposition of $S$ is
$$\text{Ph}(S)=\text{Ph}(\alpha)\text{Ph}(R)\otimes\text{Ph}(\beta)\text{Ph}(L),\quad
|S|=|\alpha||R|\otimes |\beta||L|.$$
Moreover, the polar decomposition of $T$ is given by $\text{Ph}(T)=\text{Ph}(L)\otimes 1,\quad |T|=e^{x}|L|\otimes 1$, so, using Lemma \ref{Rbeta}, one can see that $(T,S)$ is a $e^{x}$-commuting pair. In particular, the polar decomposition of $TS$ is
$$|TS|=e^{-x}|T||S|=|L||\alpha||R|\otimes |\beta||L|,\,\text{Ph}(TS)=\text{Ph}(L)\text{Ph}(\alpha)\text{Ph}(R)\otimes\text{Ph}(\beta)\text{Ph}(L).$$
But Proposition \ref{az+bRel} gives $\text{Ph}(\hat{\alpha})=\text{Ph}(L)\text{Ph}(R)\text{Ph}(\alpha)$ and $|\alpha|=|L||R||\alpha|$. Thus, we conclude that $\text{Ph}(\hat{\alpha}\otimes\hat{\beta})=\text{Ph}(TS)$ and $|\hat{\alpha}\otimes\hat{\beta}|=|TS|$ which concludes the proof of $(\ref{eqcom1})$. One can prove $(\ref{eqcom2})$ similarly.
\end{proof}

Now the proof of Theorem \ref{Theoremaz+bclassique} follows: Proposition \ref{generators2} says that $\hat{\alpha}$ and $\hat{\beta}$ generate $\widehat{M}_{x}$ and Proposition \ref{az+bRel} gives the commutation relations for $\hat{\alpha}$ and $\hat{\beta}$.

\subsection{Twisting of the quantum $az+b$ group}\label{Sectionaz+bquantique}

This Section is devoted to the proof of Theorem \ref{Theoremaz+bquantique}. Let $0<q<1$ and $(M,\Delta)$ be the $az+b$ Woronowicz' quantum group. Let $\alpha:\ L^{\infty}(\mathbb{C}^{q})\rightarrow M$ be defined by $\alpha(F)=F(a)$. Recall that (Section \ref{SectionResult}) $\widehat{\C^{q}}<(M,\Delta)$ is an abelian stable co-subgroup with the morphism $\gamma_{t}=q^{2it}\in\mathbb{C}^{q}$. Let us perform the twisting construction using the bicharacters
$$
\Psi_{x}(q^{k+i\varphi},q^{l+i\psi})=q^{ix(k\psi-l\varphi)}, \
\forall x\in\mathbb{Z},
$$
and let $(M_{x},\Delta_{x})$ be the twisted l.c. quantum group.

\begin{proposition}
$$
\Delta_{x}(a)=a\otimes a\quad\text{and}\quad \Delta_{x}(b)=u^{-x+1}|a|^{x+1}\otimes
b\dot{+}b\otimes u^{x}|a|^{-x},
$$
and $[D\varphi_{x}\,:\,D\varphi]_{t}=A_{x}^{it}=|a|^{-2ixt}$.
The modular element $\delta_{x}=|a|^{4x+2}$, the antipode is not
deformed. If $x,y\in\N$ and $x\neq y$, then $(M_{x},\Delta_{x})$ and
$(M_{y},\Delta_{y})$ are not isomorphic; if $x\neq 0$, then $(M_{x},\Delta_{x})$ and $(M_{-x},\Delta_{-x})$ are not isomorphic.
\end{proposition}
\begin{proof}
The relations of commutation from Preliminaries give
\begin{eqnarray*}
\Psi_{x}(a,q^{l+i\psi})b
&=&\Psi_{x}(u,q^{l+i\psi})\Psi_{x}(|a|,q^{l+i\psi})b\\
&=&u^{-xl}|a|^{ix\psi}v|b|\\
&=&q^{ix\psi -xl}v|b|u^{-xl}|a|^{ix\psi}\\
&=&q^{ix\psi -xl}b\Psi_{x}(a,q^{l+i\psi}).\\
\end{eqnarray*}
So, for any $\gamma\in\mathbb{C}^{q}$, one has
$$
\Psi_{x}(a\otimes 1,\gamma)(b\otimes 1)\Psi_{x}(a\otimes
1,\gamma)^{*}
=\left(\text{Phase}(\gamma)\right)^{x}|\gamma|^{-x}(b\otimes 1).
$$
Put $\Omega_{x}=(\alpha\otimes\alpha)(\Psi_{x})=\Psi_{x}(a\otimes
1,1\otimes a)$. Using the previous formula and the fact that $b\otimes 1$ and $1\otimes a$ strongly commute, one gets $\Omega_{x}(b\otimes 1)\Omega_{x}^{*} =b\otimes u^{x}|a|^{-x}.$ Similarly: $\Omega_{x}(1\otimes b)\Omega_{x}^{*} = u^{-x}|a|^{x}\otimes b.$ These formulas give the comultiplication on $b$. The comultiplication on $a$ is clear. We Since $\Psi_{s}(\gamma_{t},\gamma_{s})=1$, for all $s,t\in\R$, then $[D\varphi_{x}\,:\,D\varphi]_{t}=A_{x}^{it}=\Psi_{x}(a,\gamma_{t}^{-1})=|a|^{-2ixt}$.
Put $f_{t}^{x}=\Psi_{x} (\cdot,\gamma_{t})\Psi_{x}(\gamma_{t}^{-1},\cdot)$, then $f_{t}^{x}(q^{k+i\varphi}) = q^{4itxk}$ and $\alpha(f_{t}^{x})=|a|^{4itx}$. So, the modular element is $\delta_{x}=|a|^{2}|a|^{4x}$. The antipode is not deformed because $\Psi_{t}(x^{-1},x)=1$, for any $x$. The spectrum of the modular element is $\text{Sp}(\delta_{x})=q_{x}^{\Z}\cup\{0\},$ where $q_{x}=q^{4x+2}$, so, if $x\neq y$ are strictly positive, then $0<q_{x}\neq q_{y}<1$, so $q_{x}^{\Z}\neq q_{y}^{\Z}$, then $(M_{x},\Delta_{x})$ and  $(M_{y},\Delta_{y})$ are not isomorphic. Moreover, if $x>0$, then $(M_{x},\Delta_{x})$ is not isomorphic to $(M_{-x},\Delta_{-x})$ because in the opposite case we would have $q^{(4x+2)\Z}=q^{(4x-2)\Z}$, from where, as $x>0$, $4x+2=4x-2$ - contradiction.
\end{proof}

The group $\mathbb{C}^{q}$ is selfdual with the duality
$(q^{k+i\varphi},q^{l+i\psi})\mapsto q^{i(k\psi+l\varphi)},$ so
one can compute the representations $L$ and $R$ of
$\mathbb{C}^{q}$:
$$L_{q^{k+i\varphi}}=m^{i\varphi}\otimes s^{-k}\otimes 1\otimes s^{k},\quad
R_{q^{k+i\varphi}}=m^{-i\varphi}\otimes 1\otimes
m^{i\varphi}\otimes s^{-k}.$$ Then the left-right action of $\left(\mathbb{C}^{q}\right)^{2}$ on the
generators of $\widehat{M}$ is
\begin{equation}\label{leftright}
\alpha_{q^{k+i\varphi},q^{l+i\psi}}(\hat{a})=q^{l-k+i(\psi-\varphi)}\hat{a},\quad
\alpha_{q^{k+i\varphi},q^{l+i\psi}}(\hat{b})=q^{-k-i\varphi}\hat{b}.
\end{equation}
Let $N=\left(\mathbb{C}^{q}\right)^{2}\ltimes\hat{M}$, it is generated by the operators $\lambda_{q^{k+i\varphi},q^{l+i\psi}}$ and $\pi(x)$, for $x\in\widehat{M}$, and $\theta$ be the dual action of $\left(\mathbb{C}^{q}\right)^{2}$ on $N$. The deformed dual action is
$$\theta_{q^{k+i\varphi},q^{l+i\psi}}^{\Psi_{x}}=
\lambda_{q^{x(k-i\varphi)},q^{x(-l+i\psi)}}
\theta_{q^{k+i\varphi},q^{l+i\psi}}(.)\lambda^{*}_{q^{x(k-i\varphi)},q^{x(-l+i\psi)}}.$$
Let $\widehat{M}_{x}$ be the fixed point algebra. The left-right action is very similar to the one for the classical $az+b$. Define $\alpha=\pi(\hat{a})$, $\beta=\pi(\hat{b})$. Then $\alpha$ and $\beta$ are normal operators affiliated with $N$ and one can see that
\begin{equation}\label{actionalbe}
\theta_{q^{k+i\varphi},q^{l+i\psi}}^{\Psi_{x}}(\alpha)=q^{-x(l+k)+ix(\varphi+\psi)}\alpha\,\, ,\,\,\,\,\theta_{q^{k+i\varphi},q^{l+i\psi}}^{\Psi_{x}}(\beta)=q^{-xk+ix\varphi}\beta.
\end{equation}
Let $T_{l}$ and $T_{r}$ be the "infinitesimal generators" of the
left and right translations, so $T_{l}$ and $T_{r}$ are
affiliated with $N$ and
\begin{equation}\label{infgen}
\lambda_{q^{k+i\varphi},q^{l+i\psi}}=\left(\text{Ph}T_{l}\right)^{k}|T_{l}|^{i\varphi}\left(\text{Ph}T_{r}\right)^{l}|T_{r}|^{i\psi}.
\end{equation}
Then $\lambda(f)=f(T_{l},T_{r})\,\,\,\forall f\in
L^{\infty}\left(\left(\mathbb{C}^{q}\right)^{2}\right).$ Let
$U=\lambda(\Psi_{x})$ and $\hat{\alpha}=U^{*}\alpha U$.

\begin{proposition}\label{generators3}
$\left(T^{*}_{l} T^{*}_{r}\right)^{-x}$ and $\alpha$ strongly
commute and $\hat{\alpha}=\left(T^{*}_{l}
T^{*}_{r}\right)^{-x}\alpha$. The polar decomposition of
$\hat{\alpha}$ is
$\hat{u}:=\text{Ph}\hat{\alpha}=\left(\text{Ph}T_{l}T_{r}\right)^{x},\quad
\hat{A}:=|\hat{\alpha}|=|T_{l}T_{r}|^{-x}|\alpha|.$ Also, $|T_{l}|$
and $|\beta|$ strongly commute, so we can define a positive
operator $\hat{B}=|T_{l}|^{-x}|\beta|$. Let
$\hat{v}=\text{Ph}(T_{l})^{x}\text{Ph}(\beta)$. Then
$\hat{\alpha}$ and $\hat{B}$ are affiliated with $\widehat{M}_{x}$,
$\hat{v}\in \widehat{M}_{x}$, and we have the following relations of commutation:
\begin{itemize}
\item $\hat{u}\hat{v}=\hat{v}\hat{u}$, $\hat{A}\hat{B}=\hat{B}\hat{A}$;
\item $\hat{v}\hat{B}\hat{v}^{*}=q^{-2x}\hat{B}$, $\hat{u}\hat{B}\hat{u}^{*}=q^{-2x+1}\hat{B}$
 and, $\hat{v}\hat{A}\hat{v}^{*}=q^{-2x-1}\hat{A}$.
\end{itemize}
Moreover, these three operators generate $\widehat{M}_{x}$ in the sense that
$$\widehat{M}_{x}=\left\{f(\hat{\alpha})g(\hat{v})h(\hat{B}),\,\,f\in
L^{\infty}(\mathbb{C}^{q}),\,\,g\in L^{\infty}
(\mathbb{S}^{1}),\,\,h\in L^{\infty}(q^{\mathbb{Z}})\right\}^{''}.$$
\end{proposition}
\begin{proof}
Using $(\ref{leftright})$ and $(\ref{infgen})$, we find:
\begin{eqnarray}\label{qazrel1}
&&|T_{l}T_{r}|^{is}\alpha
|T_{l}T_{r}|^{-is}=\alpha,\\\label{qazrel2}
&&\text{Ph}(T_{l}T_{r})\alpha\text{Ph}(T_{l}T_{r})^{*}=\alpha,\\\label{qazrel3}
&&|T_{l}|^{is}\beta |T_{l}|^{-is}=q^{-is}\beta,\\\label{qazrel4}
&&\text{Ph}(T_{l})\beta\text{Ph}(T_{l})^{*}=q^{-1}\beta. \\\notag
\end{eqnarray}
Due to $(\ref{qazrel1})$ and $(\ref{qazrel2})$, $\alpha$ and
$T^{*}_{l}T^{*}_{r}$ strongly commute. Because
$\Psi_{x}(T_{r},T_{r})=1$, we have
$\hat{\alpha}=\Psi_{x}(T_{l}T_{r},T_{r})^{*}\alpha\Psi_{x}(T_{l}T_{r},T_{r}).$
Next, using
$\Psi_{x}(q^{k+i\varphi},T_{r})^{*}\alpha\Psi_{x}(q^{k+i\varphi},T_{r})
=\lambda_{1,q^{-xk+ix\varphi}}\alpha\lambda_{1,q^{-xk+ix\varphi}}^{*}=q^{-xk+ix\varphi}\alpha,$
and because $T_{l}T_{r}$ and $\alpha$ strongly commute, we have
$$\hat{\alpha}=|T_{l}T_{r}|^{-x}\left(\text{Ph}T_{l}T_{r}\right)^{x}\alpha
=\left(T^{*}_{l} T^{*}_{r}\right)^{-x}\alpha.$$ The polar
decomposition of $\hat{\alpha}$ follows. Equality
$(\ref{qazrel3})$ implies that $|T_{l}|$ and $|\beta|$ strongly
commute. Note that
\begin{eqnarray*}
\theta^{\Psi_{x}}_{q^{k+i\varphi},q^{l+i\psi}}(U)
&=&\Psi_{x}(T_{l}q^{-k-i\varphi},T_{r}q^{-l-i\psi})\\
&=&U\lambda_{q^{xl-ix\psi},q^{-xk+ix\varphi}}\Psi_{x}(q^{k+i\varphi},q^{l+i\psi}).
\end{eqnarray*}
Then, it follows from $(\ref{leftright})$ and $(\ref{actionalbe})$
that $\hat{\alpha}$ is affiliated with $\widehat{M}_{x}$. Also, using
$(\ref{actionalbe})$ we find
$\theta^{\Psi_{x}}_{q^{k+i\varphi},q^{l+i\psi}}(\hat{v})=
\left(\text{Ph}(T_{l}q^{-k-i\varphi})\right)^{x}q^{ix\varphi}\text{Ph}\beta=\hat{v},$
so $\hat{v}\in \widehat{M}_{x}$. In the same way we prove that $\hat{B}$
is affiliated with $\widehat{M}_{x}$. It is easy to see that
$\text{Ph}T_{l}$ and $\text{Ph}T_{r}$ commute with
$\text{Ph}\alpha$ and $\text{Ph}\beta$, and because
$\text{Ph}\alpha$ and $\text{Ph}\beta$ commute, it follows that
$\hat{u}\hat{v}=\hat{v}\hat{u}$. Also, $|T_{l}|$ and $|T_{r}|$
strongly commute with $|\alpha|$ and $|\beta|$, and because
$|\alpha|$ and $|\beta|$ strongly commute, it follows that
$\hat{A}\hat{B}=\hat{B}\hat{A}$. The relation
$\hat{v}\hat{B}\hat{v}^{*}=q^{-2x}\hat{B}$ follows from
$(\ref{qazrel3})$ and $(\ref{qazrel4})$. Remark that
$$\text{Ph}\alpha |T_{l}|^{-x}\text{Ph}\alpha ^{*}=
q^{-x}|T_{l}|^{-x},\quad\text{Ph}\beta
|T_{l}T_{r}|^{-x}\text{Ph}\beta ^{*}=q^{-x}|T_{l}T_{r}|^{-x},$$
and the two last relations follow from $\text{Ph}\alpha |\beta
|\text{Ph}\alpha ^{*}=q|\beta |$ and $\text{Ph}\beta |\alpha
|\text{Ph}\beta ^{*}=q^{-1}|\beta |$. The generating property is
proved as in Proposition \ref{generators2}.
\end{proof}

Let $\hat{\Delta}_{x}$ be the comultiplication on $\widehat{M}_{x}$ and
$\hat{\beta}=\hat{v}\hat{B}$. Then $\hat{\beta}$ is a closed (non-normal) operator
affiliated with $\widehat{M}_{x}$. As before, we define
$\hat{\Delta}_{x}(\hat{\beta})=\hat{\Delta}_{x}(\hat{v})\hat{\Delta}_{x}(\hat{B})$ which
is closed, non-normal and affiliated with $\widehat{M}_{x}\otimes\widehat{M}_{x}$. The
proof of the following Proposition is similar to the one of Proposition \ref{azcom}.

\begin{proposition}
$$
\hat{\Delta}_{x}(\hat{\alpha})=\hat{\alpha}\otimes\hat{\alpha}\quad\text{and}\quad
\hat{\Delta}_{x}(\hat{\beta})=\hat{\alpha}\otimes\hat{\beta}\dot{+}\hat{\beta}\otimes
1.
$$
\end{proposition}

The proof of Theorem \ref{Theoremaz+bquantique} follows from the
results of this section.

\section{Appendix}

Let $\alpha$ be an action of a \lc{} quantum group $(M,\Delta)$ on the von Neumann algebra $N$.
Let $\theta$ be a \nsf{} weight on $N$ and suppose that $N$ acts on a Hilbert space $K$ such that $(K,\iota,\Lambda_{\theta})$ is the \GNS{} construction for $\theta$. We define
$$\mathcal{D}_{0}=\text{span}\{(a\otimes 1)\alpha(x)\,|\,a\in\mathcal{N}_{\hat{\varphi}},x\in\mathcal{N}_{\theta}\}.$$
Let $(H,\iota,\Lambda)$ be the \GNS{} construction for the left invariant weight $\varphi$ of $(M,\Delta)$, $\hat{\varphi}$ the dual weight, and $\hat{\Lambda}$ its canonical G.N.S.-map. We
recall that the \GNS{} construction for the dual weight $\tilde{\theta}$ is given by $(H\otimes K,\iota,\tilde{\Lambda})$, where $\tilde{\Lambda}_{\theta}$ is the $\sigma$-strong*-norm closure
of the map
$$\mathcal{D}_{0}\rightarrow H\otimes K\,:\,
 (a\otimes1)\alpha(x)\mapsto\hat{\Lambda}(a)\otimes\Lambda_{\theta}(x).$$

\begin{proposition}\label{PropCor}
Let $C_{1}$ be a $\sigma$-strong*-norm core for $\hat{\Lambda}$ and $C_{2}$ a $\sigma$-strong*-norm core for $\Lambda_{\theta}$. Then the set $\mathcal{C}=\text{span}\{(a\otimes1)\alpha(x)\,|\,a\in C_{1},x\in C_{2}\}$ is a $\sigma$-weak-weak core for $\tilde{\Lambda}_{\theta}$.
\end{proposition}
\begin{proof}
Let $a\in\mathcal{N}_{\hat{\varphi}}$ and $x\in\mathcal{N}_{\theta}$. There exists two nets $(a_{i})$
and $(x_{i})$, with $a_{i}\in C_{1}$ and $x_{i}\in C_{2}$, such that
$$a_{i}\rightarrow a,\,\, x_{i}\rightarrow x\quad\sigma- \text{strongly}* \quad\text{and}\,\,\,\hat{\Lambda}(a_{i})\rightarrow \hat{\Lambda}(a),\,\,\Lambda_{\theta}(x_{i})\rightarrow \Lambda_{\theta}(x).$$
Thus, $(a_{i}\otimes 1)\alpha(x_{i})\rightarrow (a\otimes 1)\alpha(x)$ $\sigma$-weakly and
$$\tilde{\Lambda}_{\theta}\left((a_{i}\otimes 1)\alpha(x_{i})\right)=\hat{\Lambda}(a_{i})\otimes\Lambda_{\theta}(x_{i})\rightarrow \hat{\Lambda}(a)\otimes \Lambda_{\theta}(x)=\tilde{\Lambda}_{\theta}\left((a\otimes 1)\alpha(x)\right).$$
\end{proof}

\begin{proposition}\label{PropTx}
Let $M$ be a von Neumann algebra with a \nsf{} weight $\varphi$, $(H,\iota,\Lambda)$ the \GNS{} construction for $\varphi$, and $T$ a positive self-adjoint operator affiliated with $M$. Then
$\mathcal{C}=\{x\in\mathcal{N}_{\varphi}\,|\, Tx\,\,\,\text{ is bounded and }\,\,\,\Lambda(x)\in\mathcal{D}(T)\}$ is a $\sigma$-strong*-norm core for $\Lambda$ and, if $x\in\mathcal{C}$, then $\overline{Tx}\in\mathcal{N}_{\varphi}$ and
$\Lambda(\overline{Tx})=T\Lambda(x).$
\end{proposition}

\begin{proof}
Let $T=\int_{0}^{+\infty}\lambda de_{\lambda}$ be the spectral decomposition of $T$. Let $e_{n}=\int_{0}^{n}de_{\lambda}$. Then $e_{n}\rightarrow 1$ $\sigma$-strongly*, $Te_{n}$ is bounded with domain $H$. Let $x\in\mathcal{N}_{\varphi}$ and put $x_{n}=e_{n}x$. We have $x_{n}\rightarrow x$ $\sigma$-strongly* and $\Lambda(x_{n})=e_{n}\Lambda(x)\rightarrow\Lambda(x)$ in norm. Moreover, $Tx_{n}=Te_{n}x$ is bounded and $\Lambda(x_{n})=e_{n}\Lambda(x)\in\mathcal{D}(T)$, so $x_{n}\in\mathcal{C}$, and it follows that $\mathcal{C}$ is a $\sigma$-strong*-norm core for $\Lambda$. Now let $x\in\mathcal{C}$. Note that $e_{n}\overline{Tx}=Te_{n}x=\overline{e_{n}T}x$ is in $\mathcal{N}_{\varphi}$ and it converges $\sigma$-strongly* to $\overline{Tx}$. Moreover,
$$
\Lambda(e_{n}\overline{Tx})=
\overline{e_{n}T}\Lambda(x)=e_{n}T\Lambda(x)\rightarrow T\Lambda(x).
$$
Because $\Lambda$ is $\sigma$-strong*-norm closed, we have $\overline{Tx}\in\mathcal{N}_{\varphi}$ and $\Lambda(\overline{Tx})=T\Lambda(x)$.
\end{proof}

\begin{proposition}\label{Equality}
Let $M$ be a von Neumann algebra, $\varphi_{1}$ and $\varphi_{2}$ two \nsf{} weights on $M$ having the same modular group. Let $(H_{i},\pi_{i},\Lambda_{i})$ be the \GNS{} construction for $\varphi_{i}$ ($i=1,2$). Suppose that there exist a $\sigma$-weak-weak core $\mathcal{C}$ for $\Lambda_{1}$ such that $\mathcal{C}\subset \mathcal{N}_{\varphi_{1}}\cap\mathcal{N}_{\varphi_{2}}$ and a unitary $Z\,:\,H_{1}\rightarrow H_{2}$ such that
$\Lambda_{2}(x)=Z\Lambda_{1}(x),\quad\text{for all}\,\, x\in\mathcal{C}.$
Then $\varphi_{1}=\varphi_{2}$.
\end{proposition}

\begin{proof}
Because $\mathcal{C}$ is a $\sigma$-weak-weak core for $\Lambda_{1}$ and because $\Lambda_{2}$ is $\sigma$-weak-weak closed, we have $\mathcal{N}_{\varphi_{1}}\subset\mathcal{N}_{\varphi_{2}}$ and, for all $x\in\mathcal{N}_{\varphi_{1}}$ we have $\Lambda_{1}(x)=Z\Lambda_{2}(x).$
Thus, $\varphi_{1}(y^{*}x)=\varphi_{2}(y^{*}x)$, for all $x,y\in\mathcal{N}_{\varphi_{1}}$. Let $\mathcal{B}=\mathcal{N}^{*}_{\varphi_{1}}\mathcal{N}_{\varphi_{1}}$. This is a dense *-subalgebra of $\mathcal{M}_{\varphi_{1}}\cap\mathcal{M}_{\varphi_{2}}$ and, for all $x\in\mathcal{B}$, we have $\varphi_{1}(x)=\varphi_{2}(x)$. Because $\varphi_{1}$ and $\varphi_{2}$ have the same modular group, we can use the Pedersen-Takesaki Theorem \cite{Tak2} to conclude the proof.
\end{proof}

Let $M$ be a von Neumann algebra, $\varphi$ a \nfs{} weight on $M$,  $(H,\iota,\Lambda)$ the \gns{} construction for $\varphi$, and $\sigma$ the modular group of $\varphi$. Let $\delta$ be a positive self-adjoint operator affiliated with $M$, $\lambda>0$ such that $\sigma_{t}(\delta^{is})= \lambda^{ist}\delta^{is}$, and $\Lambda_{\delta}$ the canonical \gns{} map of the Vaes' weight $\varphi_{\delta}$. One can consider on $M\otimes M$ two \nfs{} weights: $\varphi_{\delta}\otimes\varphi_{\delta}$, with the canonical \gns{} map $\Lambda_{\delta}\otimes\Lambda_{\delta}$, and the Vaes' weight $(\varphi\otimes\varphi)_{\delta\otimes\delta}$ associated with $\varphi\otimes\varphi$, $\delta\otimes\delta$ and $\lambda^{2}$. Let $\Lambda\otimes\Lambda$ be the \gns{} map for $\varphi\otimes\varphi$, and $(\Lambda\otimes\Lambda)_{\delta\otimes\delta}$ the \gns{} map for $(\varphi\otimes\varphi)_{\delta\otimes\delta}$ (see Section \ref{sectionVaes}).

\begin{proposition}\label{propVaes}
$\varphi_{\delta}\otimes\varphi_{\delta}=
(\varphi\otimes\varphi)_{\delta\otimes\delta}$ and
$\Lambda_{\delta}\otimes\Lambda_{\delta}=
(\Lambda\otimes\Lambda)_{\delta\otimes\delta}.$
\end{proposition}
\begin{proof} Let us apply the Pedersen-Takesaki theorem to the weights
$\varphi_{1}:=(\varphi_{\delta}\otimes\varphi_{\delta})$ and $\varphi_{2}:=
(\varphi\otimes\varphi)_{\delta\otimes\delta}$ which have the same modular group
and are equal on the dense *-subalgebra $\mathcal{B}=N\odot N$ of $\mathcal{M}_{\varphi_{1}}\cap\mathcal{M}_{\varphi_{2}}$, where
$$
N :=\left\{x\in M\,|\, x\delta^{\frac{1}{2}}\,\,\text{is bounded and}\,\,\overline{x\delta^{\frac{1}{2}}}\in\mathcal{N}_{\varphi}\right\}.
$$
Let $\Lambda_{i}$ be the canonical \GNS{} map of $\varphi_{i}$. By definition, $N\odot N$ is a $\sigma$-strong*-norm core for $\Lambda_{1}$, and $\Lambda_{1}\vert_N=\Lambda_{2}\vert_N$.
Since $\Lambda_{1}$ and $\Lambda_{2}$ are $\sigma$-strongly*-norm closed, then $\Lambda_{1}\subset\Lambda_{2}.$ And $\Lambda_{1}=\Lambda_{2}$ since $\mathcal{D}(\Lambda_{1})=\mathcal{N}_{\varphi_{1}}=\mathcal{N}_{\varphi_{2}}=
\mathcal{D}(\Lambda_{2})$.
\end{proof}

Finally, let us formulate the von Neumann algebraic version of \cite{Kas}, Lemma
3.6. Let $N$ be a von Neumann algebra, $G$ a \lc{} abelian group,
$u\,:\,G\rightarrow N$ a unitary representation of $G$ and
$\theta\,:\,\hat{G}\rightarrow\text{Aut}(N)$ an action of
$\hat{G}$ on $N$ such that
$$\theta_{\gamma}(u(g))=\overline{<\gamma,g>}u(g).$$
Let $\alpha$ be the action of $G$ on $N$ implemented by $u$. The unitary representation $u$ of $G$ gives a *-homomorphism $\pi\,:\,L^{\infty}(\hat{G})\rightarrow N$.

\begin{lemma}\label{FixSub}
Let $V$ be a linear subspace of $N^{\theta}$ invariant under the
action $\alpha$ and such that
$\left(\pi(L^{\infty}(\hat{G}))V\pi(L^{\infty}(\hat{G}))\right)^{''}=N.$
Then $V^{''}=N^{\theta}.$
\end{lemma}

\end{document}